\documentclass[11pt]{article}

\usepackage[latin1]{inputenc} 
\usepackage{amsthm,amsmath,amssymb} 
\usepackage{graphicx}
\usepackage{color}
\usepackage{verbatim}

\usepackage{enumitem} 
\setlist[enumerate]{topsep=0pt,itemsep=-1ex,partopsep=1ex,parsep=1ex}

\theoremstyle{plain}
\newtheorem{theo}{Theorem}[section]

\newtheorem{lemma}[theo]{Lemma}

\theoremstyle{definition}
\newtheorem{defn}[theo]{Definition}
\newtheorem{rem}[theo]{Remark}

\newtheorem{alg}[theo]{Algorithm}

\newcommand{\mc}[1]{\mathcal{#1}}
\newcommand{\mb}[1]{\mathbb{#1}}
\newcommand{\nib}[1]{\noindent {\bf #1}}

\newcommand{\bsize}[1]{\left| #1 \right|}
\newcommand{\brak}[1]{\left[ #1 \right]}
\newcommand{\bgen}[1]{\left\langle #1 \right\rangle}
\newcommand{\sgen}[1]{\langle #1 \rangle}
\newcommand{\bfl}[1]{\left\lfloor #1 \right\rfloor}

\newcommand{\sub}{\subseteq}
\newcommand{\subn}{\subsetneq}

\newcommand{\lra}{\leftrightarrow}
\newcommand{\Lra}{\Leftrightarrow}
\newcommand{\Ra}{\Rightarrow}
\newcommand{\sm}{\setminus}
\newcommand{\ov}{\overline}
\newcommand{\ul}{\underline}

\newcommand{\eps}{\varepsilon}
\newcommand{\es}{\emptyset}
\newcommand{\pl}{\partial}

\newcommand{\ua}{\uparrow}

\newcommand{\ova}{\overrightarrow}

\newcommand{\aA}{\alpha}

\newcommand{\gG}{\gamma}
\newcommand{\dD}{\delta}

\newcommand{\lL}{\lambda}
\newcommand{\tT}{\theta}
\newcommand{\sS}{\sigma}
\newcommand{\oO}{\omega}
\newcommand{\ups}{\upsilon}

\newcommand{\GG}{\Gamma}

\newcommand{\OO}{\Omega}
\newcommand{\TT}{\Theta}
\newcommand{\Ss}{\Sigma}
\newcommand{\Ups}{\Upsilon}

\def\qed{\hfill $\Box$}

\title{The existence of designs II}

\author{Peter Keevash\thanks{Mathematical Institute,
University of Oxford, Oxford, UK. Email: keevash@maths.ox.ac.uk.
\newline \hspace*{1.8em}Research supported
in part by ERC Consolidator Grant 647678.}}

\begin{document}

\maketitle

\begin{abstract}
We generalise the existence of combinatorial designs
to the setting of subset sums in lattices with coordinates
indexed by labelled faces of simplicial complexes.
This general framework includes the problem of decomposing
hypergraphs with extra edge data, such as colours and orders, 
and so incorporates a wide range of variations on the basic design problem, 
notably Baranyai-type generalisations, such as resolvable hypergraph designs,
large sets of hypergraph designs and decompositions of designs by designs.
Our method also gives approximate counting results,
which is new for many structures whose existence was previously known,
such as high dimensional permutations or Sudoku squares.
\end{abstract}

\section{Introduction}

The existence of combinatorial designs was proved in \cite{Kexist}, to which 
we refer the reader for an introduction to and some history of the problem.
There we obtained a more general result on clique decompositions of hypergraphs,
that can be roughly understood as saying that under certain extendability conditions,
the obstructions to decomposition can already be seen 
in two natural relaxations of the problem:
the fractional relaxation (where we see geometric obstructions)
and the integer relaxation (where we see arithmetic obstructions).
The main theorem of this paper is an analogous result 
in a more general setting of lattices with coordinates
indexed by labelled faces of simplicial complexes.
There are many prerequisites for the statement of this result,
so in this introduction we will first discuss several
applications to longstanding open problems in Design Theory,
which illustrate various aspects of the general picture,
and give some indication of why it is more complicated
than one might have expected given the results of \cite{Kexist}.

\subsection{Resolvable designs}

In 1850, Kirkman formulated his famous `schoolgirls problem':

\medskip

{\em Fifteen young ladies in a school walk out three abreast for seven days in succession:
it is required to arrange them daily, so that no two shall walk twice abreast.}

\medskip

The general problem is to determine when one can 
find designs that are `resolvable':
the set of blocks can be partitioned into perfect matchings.
Recall from \cite{Kexist} that a set $S$ of $q$-subsets of an 
$n$-set $X$ is a \emph{design} with parameters $(n,q,r,\lL)$ if 
every $r$-subset of $X$ belongs to exactly $\lL$ elements of $S$,
and that such $S$ exists for any $n>n_0(q,r,\lL)$ satisfying the
necessary divisibility conditions that $\tbinom{q-i}{r-i}$ 
divides $\lL \tbinom{n-i}{r-i}$ for $0 \le i \le r-1$. 
For a resolvable design, a further necessary condition is that
$q \mid n$, otherwise there is no perfect matching on $X$,
let alone a partition of $S$ into perfect matchings!
We will show that these necessary conditions suffice for large $n$;
the case $r=2$ of this result was proved in 1973 
by Ray-Chaudhuri and Wilson \cite{RW,RW2}.

\begin{theo} \label{resolvable}
Suppose $q \ge r \ge 1$ and $\lL$ are fixed and $n>n_0(q,r,\lL)$ is large
with $q \mid n$ and $\tbinom{q-i}{r-i} \mid \lL \tbinom{n-i}{r-i}$ 
for $0 \le i \le r-1$. Then there is a resolvable $(n,q,r,\lL)$ design.
\end{theo}

It is convenient to generalise the problem and then
translate the generalisation into an equivalent 
hypergraph decomposition problem.
Suppose $G \in \mb{N}^{X_r}$ is a $r$-multigraph
supported in an $n$-set $X$, where $q \mid n$,
satisfying the necessary divisibility conditions
for a $K^r_q$-decomposition, i.e.\
$\tbinom{q-i}{r-i}$ divides $\sum \{ G_e: f \sub e\}$ 
for any $i$-set $f$ with $0 \le i \le r$
(we say $G$ is $K^r_q$-divisible).
We also suppose that $G$ is `vertex-regular'
in that $|G(x)|$ is the same for all $x \in X$
(this is clearly necessary for a resolvable decomposition).
Under certain conditions (extendability and regularity) 
to be defined later, we will show that $G$ has a
$K^r_q$-decomposition that is resolvable,
i.e.\ its $q$-sets can be partitioned into perfect matchings,
each of which can be viewed as a $K^r_q$-tiling,
i.e.\ $n/q$ vertex-disjoint $K^r_q$'s.

Now we set up an equivalent hypergraph decomposition problem.
Let $Y$ be a set of $m$ vertices disjoint from $X$, 
where $m$ is the least integer with\footnote{
We identify $G$ with its edge set,
so $|G|$ denotes the number of edges in $G$.} 
$\tbinom{m}{r-1} \ge q|G|/Qn$, with $Q=\tbinom{q}{r}$.
Let $J$ be an $(r-1)$-graph on $Y$ with 
$|J|=q|G|/Qn = \tbinom{q-1}{r-1}^{-1}|G(x)|$ for all $x \in X$.
Let $G'$ be the $r$-multigraph obtained from $G$
by adding as edges (with multiplicity one)
all $r$-sets of the form $f \cup \{x\}$
where $f \in J$ and $x \in X$. 
Let $H$ be the $r$-graph whose vertex set is the disjoint
union of a $q$-set $A$ and an $(r-1)$-set $B$,
and whose edges consist of all $r$-sets in $A \cup B$
that are contained in $A$ or have exactly one vertex in $A$.

We claim that a resolvable $K^r_q$-decomposition of $G$
is equivalent to an $H$-decomposition of $G'$.
To see this, suppose that an $H$-decomposition $\mc{H}$ 
of $G'$ is given. Let $\mc{A}$ be the set of the restrictions 
to $A$ of the copies of $H$ in $\mc{H}$.
Then $\mc{A}$ is a $K^r_q$-decomposition of $G$.
Furthermore, for each $f \in J$ and $x \in X$
there is a unique copy of $H$ in $\mc{H}$
containing $f \cup \{x\}$, so the elements of $\mc{A}$
corresponding to copies of $H$ containing $f$ 
form a perfect matching of $X$, and every element of $\mc{A}$ 
is thus obtained from a unique such $f$, so $\mc{A}$ is resolvable.
Conversely, any resolvable $K^r_q$-decomposition of $G$
can be converted into an $H$-decomposition of $G'$
by assigning an edge of $J$ to each perfect matching in
the resolution and forming copies of $H$ in the obvious way.

For general $r$-graphs $H$, Glock, K\"uhn, Lo and Osthus \cite{GKLO2}
solved the $H$-decomposition problem for certain structures
that they call `supercomplexes', which in particular solves the
problem for $r$-graphs $G$ that are `typical' or have large minimum degree:
they show that in this setting there is an $H$-decomposition if $G$ is
$H$-divisible, i.e.\ every $i$-set degree is divisible by the
greatest common divisor of the $i$-set degrees in $H$.
However, the bipartite form of the problem considered here is
not covered by their framework; indeed, we will see later in more
generality that there are additional complications in partite settings.
Our general result on $H$-decompositions will be of the form discussed
above, i.e.\ that under certain extendability conditions, 
the obstructions to decomposition appear in the 
fractional or integer relaxation. However, one should note that
there can be additional obstructions in the integer relaxation 
besides the divisibility conditions mentioned above.

\subsection{Baranyai-type designs}

Next we develop the theme suggested by the construction
in the previous section, namely that of obtaining variations
on the basic design problem that are equivalent to certain partite
hypergraph decomposition problems. We will call these Baranyai-type designs,
after the classical result of Baranyai \cite{B} that any complete $r$-graph
$K^r_n$ with $r \mid n$ can be partitioned into perfect matchings.

One natural question of this type is whether $K^q_n$ can be decomposed into 
$(n,q,r,\lL)$-designs; such a decomposition is known as a `large set' of designs. 
Besides the necessary divisibility conditions discussed above 
for the existence of one such design, another obvious necessary condition
is that the size $\lL \tbinom{q}{r}^{-1} \tbinom{n}{r}$ of each design 
in the decomposition should divide the size $\tbinom{n}{q}$ of $K^q_n$;
equivalently, we need $\lL \mid \tbinom{n-r}{q-r}$.
It is natural to conjecture that these necessary divisibility conditions 
should be sufficient, apart from a finite number of exceptions.
Even in the very special case of Steiner Triple Systems,
this was a longstanding open problem, settled in 1991 by Teirlinck \cite{T}.
Lovett, Rao and Vardy \cite{LRV} extended the method of \cite{KLP} 
to show that if $q>9r$, $\ell \in \mb{N}$
and $n>n_0(q,r,\ell)$ satisfies the divisibility conditions 
then there is a large set of $(n,q,r,\lL)$-designs, 
where $\tbinom{n}{q} = \ell \lL \tbinom{q}{r}^{-1} \tbinom{n}{r}$.
This settles the existence conjecture when the edge multiplicity $\lL$
is within a constant factor of the maximum possible multiplicity, 
but leaves it open otherwise (for example, it does not include
the case of large sets of Steiner systems).
We will prove the general form of the 
existence conjecture for large sets of designs.

\begin{theo} \label{large}
Suppose $q \ge r \ge 1$ are fixed, $n>n_0(q,r)$ is large
and $\lL \mid \tbinom{n-r}{q-r}$ with all
$\tbinom{q-i}{r-i} \mid \lL \tbinom{n-i}{r-i}$.
Then there is a large set of $(n,q,r,\lL)$ designs.
\end{theo}

As for resolvable designs, we can consider the more general problem of 
decomposing any $q$-multigraph $G$ on an $n$-set $X$ into $(n,q,r,\lL)$-designs.
This clarifies the general form of the divisibility conditions,
as there are several conditions that collapse into one
in the case that $G=K^q_n$. Indeed, for each $0 \le i \le r$
and $i$-set $I \sub [n]$ we need the degree $|G(I)|$ of $I$
to be divisible by the number 
$Z_i := \lL \tbinom{q-i}{r-i}^{-1} \tbinom{n-i}{r-i}$
of $q$-sets containing $I$ in any $(n,q,r,\lL)$-design.
Furthermore, we clearly need $G$ to be an `$r$-multidesign',
meaning that all $|G(e)|$ with $e \in [n]_r$ are equal.

Again we formulate an equivalent hypergraph decomposition problem.
Let $Y$ be a set of $m$ vertices disjoint from $X$,
where $m$ is the least integer with 
$\tbinom{m}{q-r} \ge |G| / Z_0$.
Let $J$ be an $(q-r)$-graph on $Y$ with $|J|=|G|/Z_0$.
Let $G'$ be the $q$-multigraph obtained from $G$
by adding as edges with multiplicity $\lL$
all $q$-sets of the form $e \cup f$
with $e \sub X$ and $f \in J$.
Let $H$ be the $q$-graph whose vertex set is the disjoint
union of a $q$-set $A$ and a $(q-r)$-set $B$,
and whose edges consist of $A$ 
and all $q$-sets in $A \cup B$ that contain $B$.
Then a decomposition of $G$ into $(n,q,r,\lL)$-designs
is equivalent to an $H$-decomposition of $G$.
Indeed, given an $H$-decomposition $\mc{H}$ of $G$,
each edge of $G$ appears as exactly one copy of $A$ in $\mc{H}$,
and for each $f \in J$ the copies of $A$ within the copies of $H$
that contain $f$ form an $(n,q,r,\lL)$-design.
Conversely, any decomposition of $G$ into $(n,q,r,\lL)$-designs
can be converted into an $H$-decomposition of $G$
by assigning an edge of $J$ to each design in the decomposition.

Another natural example of a Baranyai-type design is what
we will call a `complete resolution' of $K^q_n$: 
we partition $K^q_n$ into Steiner $(n,q,q-1)$ systems, 
each of which is partitioned into Steiner $(n,q,q-2)$ systems,
and so on, down to Steiner $(n,q,1)$ systems
(which are perfect matchings).
Again we show that this exists for $n>n_0(q)$ 
under the necessary divisibility conditions,
which take the simple form $q-j \mid n-j$ for $0 \le j < q$,
i.e.\ $n = q$ mod $lcm([q])$.

\begin{theo} \label{resolution}
Suppose $q$ is fixed and $n>n_0(q)$ is large
with $n = q$ mod $lcm([q])$. 
Then there is a complete resolution of $K^q_n$.
\end{theo}

To formulate an equivalent hypergraph decomposition problem,
we consider disjoint sets of vertices $X$ and $Y$ where $|X|=n$ 
and $Y$ is partitioned into $Y_j$, $0 \le j < q$ 
with $|Y_j| = \tfrac{n-j}{q-j}$. We let $G'$ be the $q$-graph 
whose edges are all $q$-sets $e \sub X \cup Y$ 
such that $|e \cap Y_j| \le 1$ for all $0 \le j < q$,
and if $e \cap Y_j \ne\es$ then $e \cap Y_i \ne\es$ 
for all $i>j$. Let $H$ be the $q$-graph whose vertex set 
is the disjoint union of two $q$-sets $A$ and 
$B = \{b_0,\dots,b_{q-1}\}$, whose edges are all $q$-sets 
$e \sub A \cup B$ such that if $b_j \in e$ 
then $b_i \in e$ for all $i>j$. 

Then a complete resolution of $K^q_n$
is equivalent to an $H$-decomposition of $G'$.
Indeed, given an $H$-decomposition $\mc{H}$ of $G'$,
we note that for any $y_i \in Y_i$ for $j \le i \le q$
the set of copies of $A$ in the copies of $H$ in $\mc{H}$
that contain $\{y_j,\dots,y_q\}$ form a Steiner $(n,q,j-1)$ system,
and as $y_j$ ranges over $Y_j$ we obtain a partition
of the Steiner $(n,q,j)$ system corresponding
to the copies of $H$ in $\mc{H}$
that contain $\{y_{j+1},\dots,y_q\}$.
Conversely, a complete resolution of $K^q_n$
can be converted into an $H$-decomposition of $G'$
by iteratively assigning vertices of $Y_j$
to the Steiner $(n,q,j-1)$ systems that
decompose each Steiner $(n,q,j)$ system.

\subsection{Partite decompositions}

The above applications demonstrate the need for 
hypergraph decomposition in various partite settings.
We defer our general statement
and just give here some easily stated particular cases.
First we consider the nonpartite setting and
the typicality condition from \cite{Kexist}.

\begin{defn} \label{def:typ}
Suppose $G$ is an $r$-graph on $[n]$.
The density of $G$ is $d(G) = |G| \tbinom{n}{r}^{-1}$.
We say that $G$ is \emph{$(c,s)$-typical} if for any set $A$ 
of $(r-1)$-subsets of $V(G)$ with $|A| \le s$ we have 
$\bsize{\cap_{f \in A} G(f)} = (1 \pm |A|c) d(G)^{|A|} n$.
\end{defn}

We show that any typical $r$-graph has an $H$-decomposition
provided that it satisfies the necessary divisibility condition
discussed above (this result was also proved in \cite{GKLO2}).

\begin{theo} \label{Hdecomp:typ}
Let $H$ be an $r$-graph on $[q]$
and $G$ be an $H$-divisible
$(c,h^q)$-typical $r$-graph on $[n]$,
where $n=|V(G)|>n_0(q)$ is large,
$h=2^{50q^3}$, $\dD = 2^{-10^3 q^5}$,
$d(G) > 2n^{-\dD/h^q}$, $c < c_0 d(G)^{h^{30q}}$
where $c_0=c_0(q)$ is small.
Then $G$ has an $H$-decomposition.
\end{theo}
 
Next we consider the other extreme in terms of partite settings.

\begin{defn} \label{Hblowup}
Let $H$ be an $r$-graph. 
We call an $r$-graph $G$ an $H$-blowup
if $V(G)$ is partitioned as $(V_x: x \in V(H))$
and each $e \in G$ is $f$-partite for some $f \in H$,
i.e.\ $f = \{x: e \cap V_x \ne \es \}$.

We write $G_f$ for the set of $f$-partite $e \in G$.
For $f \in H$ let $d_f(G) = |G_f| \prod_{x \in f} |V_x|^{-1}$.
We call $G$ a $(c,s)$-typical $H$-blowup if 
for any $s' \le s$ and distinct $e_1,\dots,e_{s'}$ where
each $e_j$ is $f_j$-partite for some $f_j \in V(H)_{r-1}$,
and any $x \in \cap_{j=1}^{s'} H(f_j)$ we have 
\[ \bsize{V_x \cap \bigcap_{j=1}^{s'} G(e_j)}
= (1 \pm s'c) |V_x| \prod_{j=1}^{s'} d_{f_j+x}(G).\]

We say $G$ has a partite $H$-decomposition if 
it has an $H$-decomposition using copies of $H$ 
with one vertex in each part $V_x$.

We say $G$ is $H$-balanced if for every $f \sub V(H)$
and $f$-partite $e \sub V(G)$ there is some $n_e$
such that $|G_{f'}(e)|=n_e$ for all $f \sub f' \in H$.
\end{defn}

Note that if $G$ has a partite $H$-decomposition
then $G$ is $H$-balanced; we establish the converse 
for typical $H$-blowups.

\begin{theo} \label{Hdecomp:partite1}
Let $H$ be an $r$-graph on $[q]$
and $G$ be an $H$-balanced $(c,h^q)$-typical
$H$-blowup on $(V_x: x \in V(H))$,
where each $n/h \le |V_x| \le n$
for some large $n > n_0(q)$ and $h=2^{50q^3}$,
$\dD = 2^{-10^3 q^5}$, $d_f(G) > d > 2n^{-\dD/h^q}$ 
for all $f \in H$ and $c < c_0 d^{h^{30q}}$,
where $c_0=c_0(q)$ is small.
Then $G$ has a partite $H$-decomposition.
\end{theo}

For example, if $H=K^r_{r+1}$ and $G=K^r_{r+1}(n)$
is the complete $(r+1)$-partite $r$-graph
with $n$ vertices in each part then Theorem \ref{Hdecomp:partite1}
shows the existence of an object known variously as
an $r$-dimensional permutation or latin hypercube.
(It can be viewed as an assignment of $0$ or $1$
to the elements of $[n]^{r+1}$ so that every line 
has a unique $1$, or as an assignment of $[n]$
to the elements of $[n]^r$ so that each line
contains every element of $[n]$ exactly once.)
The result for general $G$
implies a lower bound on the number of
$r$-dimensional permutations:
we can estimate the number
of choices for an almost $H$-decomposition by
analysing a random greedy algorithm, and show
that almost all of these can be completed
to an (actual) $H$-decomposition
(we omit the details of the proof, 
which are similar to those in \cite{Kcount}).
In combination with the upper bound of 
Linial and Luria \cite{LL} we obtain the
following answer to an open problem from \cite{LL}.

\begin{theo}
The number of $r$-dimensional permutations 
of order $n$ is $(n/e^r + o(n))^{n^r}$.
\end{theo}

More generally, many applications of our main theorem
can be similarly converted to an approximate counting result,
where the upper bound comes from a general bound by Luria \cite{L}
on the number of perfect matchings in a uniform hypergraph 
with small codegrees (for example, we could give such estimates
for the number of resolvable designs or large sets of designs).
Another example\footnote{
Alexey Pokrovskiy drew this to my attention.
Our theorem does not apply to the construction
given by Luria \cite{L}, but it is not hard
to give a suitable alternative construction.
For example, let $H$ be the $4$-graph with
$V(H) = \{x_1,x_2,y_1,y_2,z_1,z_2\}$ and
$E(H) = \{ x_1x_2y_1y_2, x_1x_2z_1z_2, 
y_1y_2z_1z_2, x_1y_1z_1z_2 \}$.
Then an $H$-decomposition of the complete $n$-blowup of $H$
can be viewed as a Sudoku square, 
where we represent rows by pairs $(a_1,a_2)$,
columns by $(b_1,b_2)$, symbols by $(c_1,c_2)$
and boxes by $(a_1,b_1)$; a copy of $H$
with vertices $\{a_1,a_2,b_1,b_2,c_1,c_2\}$
represents a cell in row $(a_1,a_2)$
and column $(b_1,b_2)$ with symbol $(c_1,c_2)$.}
is the following estimate
for the number of (generalised) Sudoku squares
(the theorem says nothing about the squares 
of the popular puzzle, in which $n=3$).

\begin{theo}
The number of Sudoku squares with $n^2$ boxes 
of order $n$ is $(n^2/e^3 + o(n^2))^{n^4}$.
\end{theo}

\subsection{Colours and labels}

There are many questions in design theory that are naturally expressed 
as a decomposition problem for hypergraphs with extra data associated
to edges, such as colours or vertex labels. A decomposition theorem
for coloured multidigraphs with several such applications was
given by Lamken and Wilson \cite{LW}. Here we illustrate one such 
application and an example of a hypergraph generalisation. 
(There are many other such applications, but for the sake of 
brevity we leave a detailed study for future research.)

The Whist Tournament Problem (posed in 1896 by Moore \cite{M}) 
is to find a schedule of games for $4n$ players, 
where in each game two players oppose two others, such that 
(1) the games are arranged into rounds, where
each player plays in exactly one game in each of the rounds,
(2) each player partners every other player exactly once
and opposes every other player exactly twice.
(There is also a similar problem for $4n+1$ players
in which one player sits out in each round.)
Whist Tournaments exist for all $n$ (see \cite[Chapter VI.64]{CD}).
If we remove condition (1) we obtain the Whist Table Problem.
As observed in \cite{LW}, we obtain an equivalent form of the latter 
by considering a red/blue coloured multigraph on the set of players,
where between each pair of players there is one red edge (`partner')
and two blue edges (`oppose'), and we seek a decomposition into
copies of $K_4$ coloured as a blue $C_4$ with two red diagonals.
The Whist Tournament Problem is equivalent to a partite 
decomposition problem that fits into our framework, 
but not that of \cite{LW} (which only covers the
case of $4n+1$ players).

There are many ways to formulate similar problems 
with more complexities, such as larger teams
and particular roles for players within teams.
Here we describe a fictional illustration of this idea,
which we may call a `tryst tournament'
(sports aficionados will no doubt
be able to provide real examples).
A tryst team consists of three players,
one of whom is designated the captain.
A tryst game is played by nine players
divided into three tryst teams.
The Tryst Table Problem is to find a schedule
of tryst games for $n$ players, such that
(1) for every triple $T$ of players and every $x \in T$
there is exactly one game in which $T$ is a team 
and $x$ is the captain,
(2) for every triple $T$ of players
there is exactly one game in which $T$ is the set 
of captains of the three teams in that game.\footnote{
We choose these simple rules for simplicity of exposition,
and there is no doubt a direct proof
of Theorem \ref{tryst} not using our main theorem.
The point is that one can use the same method to analyse 
variations with more rules, such as 
a Tryst Tournament Problem (arranging the games into rounds)
and/or constraining more triples, e.g.\
we could also ask for every triple $T$ of players 
and every $x \in T$ to have exactly two games 
in which $x$ captains a team and $T \sm \{x\}$
is the set of non-captains in a different team.}

\begin{theo} \label{tryst}
The Tryst Table Problem has a solution
for all sufficiently large $n$.
\end{theo}

We reformulate the Tryst Table Problem 
(somewhat vaguely at first) as follows.
Form a `structure' $G$ on the set $V$ of players by
including a red triple (`captains') for each triple
and a blue `pointed' triple (`teams') for each
triple $T$ and $x \in T$. We want to decompose $G$
by copies of a `structure' $H$ on $9$ vertices,
with $3$ vertex-disjoint blue pointed triples, and 
a red triple consisting of the points of the blue triples.

To make sense of the undefined terms just used
we now switch to a setting in which all edges 
come with labels on their vertices, 
so our fundamental object becomes a set of functions
(instead of a hypergraph, which is a set of sets).
For the Tryst Table Problem, we let $G^*$
contain a red copy and a blue copy of
each injection from $[3]$ to $V$.
We define a set $H^*$ of red and blue injections
from $[3]$ to $[9]$ as follows,
in which we imagine that three teams are labelled
$123$, $456$ and $789$ with captains $1$, $4$ and $7$.
The red functions of $H^*$ consist of 
all bijections from $[3]$ to $147$.
The blue functions of $H^*$ consist of all bijections
from $[3]$ to one of the teams $123$, $456$ or $789$, 
such that $1$ is mapped to the captain. A copy of $H^*$ in $G^*$ 
is defined by fixing any injection $\phi:[9] \to V$
and composing all functions in $H^*$ with $\phi$;
the interpretation of this copy is a tryst game
between teams $\phi(123)$, $\phi(456)$ and $\phi(789)$ 
with captains $\phi(1)$, $\phi(4)$ and $\phi(7)$.
It is clear that the Tryst Table Problem is equivalent
to finding an $H^*$-decomposition of $G^*$.

Our main theorem is a decomposition result for vectors
where coordinates are indexed by functions and
take values in some lattice $\mb{Z}^D$.
The `subcoordinates' in $\mb{Z}^D$ may be interpreted
as colours, so e.g.\ we may think of 
$J_\psi = (2,3) \in \mb{Z}^2$ as saying that
$J$ has $2$ red copies and $3$ blue copies 
of some function $\psi$. This general framework
includes all of the problems discussed above
and many other variations thereupon
(see subsection \ref{sec:vvd} for more examples).

One consequence of our main theorem is a generalisation
of the hypergraph decomposition result alluded to above
to decompositions of coloured multihypergraphs
by coloured hypergraphs. It seems hard to describe
the divisibility conditions in general, so here we
will specialise to the setting of rainbow clique decompositions,
for which the divisibility conditions are quite simple.
We write $[\tbinom{q}{r}]K^r_n$ for the $r$-multigraph on $[n]$
in which there are $\tbinom{q}{r}$ copies of each $r$-set
coloured by $[\tbinom{q}{r}]=\{1,\dots,\tbinom{q}{r}\}$.
We ask when $[\tbinom{q}{r}]K^r_n$ can be decomposed
into rainbow copies of $K^r_q$, i.e.\ copies of $K^r_q$
in $[\tbinom{q}{r}]K^r_n$ in which the colours of edges
are all distinct. A stricter version of the question is
to fix some rainbow colouring of $K^r_q$ and only allow
the decomposition to use copies of $K^r_q$ that
are isomorphic to the fixed rainbow colouring.
We will answer both versions of the question.

First we consider the question in which
we allow any rainbow $K^r_q$. Ignoring colours,
we have the same necessary divisibility condition as before
for the multigraph $\tbinom{q}{r} K^r_n$ 
to have a $K^r_q$-decomposition, namely
$\tbinom{q-i}{r-i} \mid \tbinom{q}{r} \tbinom{n-i}{r-i}$
for $0 \le i \le r-1$. We will show that under the
same conditions we even have a rainbow $K^r_q$-decomposition.

\begin{theo} \label{rainbow:all}
Suppose $n>n_0(q)$ is large and 
$\tbinom{q-i}{r-i} \mid \tbinom{q}{r} \tbinom{n-i}{r-i}$
for $0 \le i \le r-1$. Then $[\tbinom{q}{r}]K^r_n$ has 
a rainbow $K^r_q$-decomposition.
\end{theo}

Now suppose that we only allow
copies of some fixed rainbow colouring.
For convenient notation we
identify the set of colours with 
$[q]_r := \{B \sub [q]: |B|=r \}$
and suppose that in the fixed colouring
of $[q]_r$ we colour each set by itself.
We write $[q]_rK^r_n$ for the corresponding
relabelling of $[\tbinom{q}{r}]K^r_n$.
Any injection $\phi:[q] \to [n]$ defines
a copy of $[q]_r$ where for each $B \in [q]_r$
we use the colour $B$ copy of $\phi(B)$.
We say $[q]_rK^r_n$ has a $[q]_r$-decomposition
if it can be decomposed into such copies.

\begin{theo} \label{rainbow:fixed}
Suppose $n>n_0(q)$ is large and 
$\tbinom{r}{i} \mid \tbinom{n-i}{r-i}$
for $0 \le i \le r-1$. 
Then $[q]_rK^r_n$ has a $[q]_r$-decomposition.
\end{theo}

The divisibility conditions in Theorem \ref{rainbow:fixed} 
are necessary for $r \le q/2$ but not in general.\footnote{
The precise divisibility conditions can
be extracted from our characterisation of 
the divisibility lattice, but we omit
this for the sake of a simple illustration.}
To see necessity, suppose $\mc{D}$ is a
$[q]_r$-decomposition of $[q]_rK^r_n$.
Identify each copy $\phi([q]_r)$
with a vector $v^\phi \in (\mb{Z}^{[q]_r})^{K^r_n}$
where each $v^\phi_{\phi(B)} = e_B$
(the standard basis vector for $B \in [q]_r$).
For any $f \in [n]_i$ we have
$\sum \{ v^\phi_e: f \sub e \in K^r_n, \phi \in \mc{D} \}
= \tbinom{n-i}{r-i} \ul{1} \in \mb{Z}^{[q]_r}$ 
equal to $\tbinom{n-i}{r-i}$ in each coordinate. 
On the other hand, the contribution of any given
$\phi([q]_r)$ to this sum is
$\sum \{ e_B: \phi^{-1}(f) \sub B \}$,
which is a row of the inclusion matrix $M^r_i(q)$:
this has rows indexed by $[q]_i$, columns by $[q]_r$
and each $M^r_i(q)_{fe} = 1_{f \sub e}$.
As $M^r_i(q)$ has full row rank
(by Gottlieb's Theorem \cite{Go}, using $r \le q/2$),
each row must appear the same number of times,
say $m$, and then any $B \in [q]_r$ contributes
$m \tbinom{r}{i}$ times, so
$\tbinom{r}{i} \mid \tbinom{n-i}{r-i}$.

\subsection{Decomposition lattices}

Here we will give some indication of what new ideas 
are needed in the general setting besides those in \cite{Kexist}.
We start by recalling the proof strategy in \cite{Kexist}
for clique decompositions of hypergraphs that are extendable
and have no obstruction (fractional or integer) to decomposition.
The first step is a Randomised Algebraic Construction,
which results in a partial decomposition (the `template')
that covers a constant fraction of the edge set,
and carries a rich structure of possible local modifications.
By the nibble and other random greedy algorithms, 
we can choose another partial decomposition that covers 
all edges not in the template, which also spills over 
slightly into the template, so that every edge is covered
once or twice, and very few edges (the `spill') are covered twice.

Next we find an `integral decomposition' of the spill,
and apply a `clique exchange algorithm'
that replaces the integral decomposition 
by a `signed decomposition',
i.e.\ two partial decompositions,
called `positive' and `negative',
such that the underlying hypergraph 
of the negative decomposition is contained 
in that of the positive decomposition,
and the difference forms a `hole'
that is precisely equal to the spill.
We also ensure that for each positive clique $Q^+$
there is a modification (a `cascade') to the clique decomposition 
of the template so that it contains $Q^+$.
Then deleting the positive cliques and replacing them by the
negative cliques eliminates one of the two uses of each edge in
the spill, and we end up with a perfect decomposition.

In broad terms, the strategy of the proof in this paper is similar. 
Furthermore, much of the proof for designs can be adapted to 
the general setting (we give the details of this in section 3).
However, the `integral decomposition' step becomes much more difficult
(and it is necessary to overcome these difficulties, as this
is a relaxation of the problem we are trying to solve).
There are in fact two aspects of this step,
which both become much more difficult:
(i) characterising the decomposition lattice
(i.e.\ the set of $\mb{Z}$-linear combinations
of the vectors that we allow in our decomposition),
(ii) finding bounded integral decompositions.
Regarding (ii), there is a `local decoding' trick for designs
that greatly simplifies the proof of \cite[Lemma 5.12]{Kexist} (version 2),
but in the general setting we do not have local decodability,
so we revert to the original randomised rounding 
method of \cite{Kexist} (version 1).

As for (i), we give a cautionary example here to show that 
the decomposition lattice is in general not what one might guess
given its simple structure for designs. Let us first recall the
characterisation by Graver and Jurkat \cite{GJ} and Wilson \cite{W4}
of the set of $J \in \mb{Z}^{K^r_n}$ with an integral $K^r_q$-decomposition,
i.e.\ the $\mb{Z}$-linear combinations of (characteristic vectors of)
copies of $K^r_q$. Clearly, any such $J$ is $K^r_q$-divisible
(as defined above), and the converse is also true for $n \ge q+r$.
For partite problems there may be additional `balance' constraints
(as in Theorem \ref{Hdecomp:partite1}) 
but there may be further more subtle constraints.

Let us consider the decomposition lattice of 
the triangles of a rainbow $K_4$, defined as follows.
Fix any bijection $b:E(K_4) \to [6]$ 
and colouring $c:E(K_n) \to [6]$.
Let $\mc{B}$ be the set of 
all $b$-coloured copies of $K^3_4$,
i.e.\ for each injection $\phi:[4] \to [n]$
such that all $c(\phi(i)\phi(j)) = b(ij)$
we include $\{\phi([4] \sm \{i\}): i \in [4]\}$.
We wish to characterise the lattice
$\bgen{\mc{B}}$ generated by $\mc{B}$.
Certainly any $J \in \bgen{\mc{B}}$
must be $K^3_4$-divisible and supported on the
set $\mc{T}$ of triangles that appear in $\mc{B}$.
One might guess by analogy with the lattice of $K^3_4$'s
in a random $3$-graph (see \cite{Kexist})
that if $c$ is random then whp
there would be no further condition.

However, we will now describe a $K^3_4$-divisible vector
(a `twisted octahedron') that is not in $\bgen{\mc{B}}$.
Suppose that $b(12)=3$, $b(13)=2$, $b(23)=1$,
$b(14)=4$, $b(24)=5$, $b(34)=6$.
Consider any octahedron, i.e.\ complete tripartite graph,
with parts $\{x_0,x_1\}$, $\{y_0,y_1\}$, $\{z_0,z_1\}$
coloured so that all $c(y_i z_j) = 1$,
$c(x_0 y_0) = c(x_0 y_1) = c(x_1 z_0) = c(x_1 z_1) = 2$,
$c(x_1 y_0) = c(x_1 y_1) = c(x_0 z_0) = c(x_0 z_1) = 3$.
Then every triangle $x_i y_j z_k$ is rainbow.
Let $J \in \mb{Z}^{\mc{T}}$ be $\pm 1$ on these triangles
and $0$ otherwise, with $J_{x_i y_j z_k} = (-1)^{i+j+k}$.
Then $J$ is null ($\sum \{ J_e: f \sub e \in \mc{T}\} = 0$
whenever $|f| \le 2$) so is $K^3_4$-divisible.

To see $J \notin \bgen{\mc{B}}$ we use the
colouring to define an algebraic invariant.
For each triangle $T = xy_0z_0 \in \mc{T}$ 
containing $y_0z_0$ we let $f(T)$ be one of
the standard basis vectors of $\mb{Z}^4$
according to the colouring of $T$:
we let $f(T)$ be $e_1$, $e_2$, $e_3$ or $e_4$
according to whether $(c(xy_0),c(xz_0))$ is
$(2,3)$, $(3,2)$, $(5,6)$ or $(6,5)$.
We extend $f$ to a $\mb{Z}$-linear map
from $\mb{Z}^{\mc{T}}$ to $\mb{Z}^6$, where 
$f(T')=0$ if $T'$ does not contain $y_0z_0$.
If $J' \in \bgen{\mc{B}}$ then $f(J')$ 
lies in the lattice generated by
$(1,0,1,0)$ and $(0,1,0,1)$.
However, we have $f(J)=(1,-1,0,0)$,
so $J \notin \bgen{\mc{B}}$.
This example hints at the importance of vertex labels
(even when not explicitly presented in a problem)
when characterising decomposition lattices.

\subsection{Organisation}

The organisation of this paper is as follows.
In the next section we set up our general framework
of labelled complexes and vector-valued decompositions,
and develop some basic theory of these definitions
that will be used throughout the paper.
In section 3 we state our main theorem and
present those parts of the proof that are 
somewhat similar to those in \cite{Kexist}.
We define and analyse the Clique Exchange Algorithm
in section 4. We complete the proof
of our main theorem by solving the problems of
integral decomposition (section 5) and
bounded integral decomposition (section 6).
Section 7 gives several applications of our main theorem.

\subsection{Notation}

Most of the following notation is as in \cite{Kexist}.
Let $[n] = \{1,\dots,n\}$.
Let $\tbinom{S}{r}$ denote the set of $r$-subsets of $S$.
We write $Q = \tbinom{[q]}{r}$ and also $Q = \tbinom{q}{r}$
(the use will be clear from the context).
We identify $Q = \tbinom{[q]}{r}$ with the
edge set of $K^r_q$ (the complete $r$-graph on $[q]$).
We write $K^r_q(S)$ for the complete
$q$-partite $r$-graph with parts of size $|S|$ 
where each part is identified with $S$.
If $S=[s]$ we write $K^r_q(S)=K^r_q(s)$.

We use `concatenation notation' for sets
($xyz$ may denote $\{x,y,z\}$)
and for function composition
($fg$ may denote $f \circ g$).

We say $E$ holds with high probability (whp) 
if $\mb{P}(E) = 1-e^{-\Omega(n^c)}$ for some $c>0$ as $n \to \infty$.

We write $Y^X$ for the set of vectors with entries in $Y$
and coordinates indexed by $X$, which we also identify
with the set of functions $f:X \to Y$.
For example, we may consider $v \in \mb{F}_p^q$ 
as an element of a vector space over $\mb{F}_p$
or as a function from $[q]$ to $\mb{F}_p$.

We identify $v \in \{0,1\}^X$ with the set $\{x \in X: v_x=1\}$,
and $v \in \mb{N}^X$ with the multiset in $X$ 
in which $x$ has multiplicity $v_x$
(for our purposes $0 \in \mb{N}$). 
We often consider algorithms with input $v \in \mb{Z}^X$,
where each $x \in X$ is considered $|v_x|$ times,
with a sign attached to it (the same as that of $v_x$);
then we refer to $x$ as a `signed element' of $v$.

If $G$ is a hypergraph, $v \in \mb{Z}^G$ and $e \in G$
we define $v(e) \in \mb{Z}^{G(e)}$ 
by $v(e)_f = v_{e \cup f}$ for $f \in G(e)$.

We denote the standard basis vectors in $\mb{R}^d$
by $e_1,\dots,e_d$. Given $I \sub [d]$, 
we let $e_I$ denote the $I$ by $[d]$ matrix 
in which the row indexed by $i \in I$ is $e_i$.

We write $M \in \mb{F}_p^{q \times r}$ to mean that $M$ 
is a matrix with $q$ rows and $r$ columns
having entries in $\mb{F}_p$.
For $I \in Q = \tbinom{[q]}{r}$ we let $M_I$ be the 
square submatrix with rows indexed by $I$.
Note that $M_I = e_I M$.

We will regard $\mb{F}_{p^a}$ as a vector space 
over $\mb{F}_p$. For $e \sub \mb{F}_{p^a}$
we write $\dim(e)$ for the dimension of
the subspace spanned by the elements of $e$.
For $e \in \mb{F}_{p^a}^d$ we write $\dim(e)$
for the dimension of the set of coordinates of $e$.

When we use `big-O' notation, 
the implicit constant will depend only on $q$.

We write $a = b \pm c$ to mean $b - c \leq a \leq b + c$.

Throughout the paper we omit floor and ceiling symbols 
where they do not affect the argument. 

We also use the following notation (not from \cite{Kexist}).

Let $Bij(B,B')$ denote the set of bijections from $B$ to $B'$.

Let $Inj(B,V)$ denote the set of injections from $B$ to $V$.

We extend our concatenation notation to sets of functions,
e.g.\ $\phi\Ups = \{ \phi \circ \psi: \psi \in \Ups \}$.

We let $\es$ denote the empty set and
also the unique function with empty domain.

For any set $X$ write $X_j = \tbinom{X}{j}$ 
(convenient notation for use in exponents).

We write $[q](S)$ for the set of partite maps
$f:[q] \to [q] \times S$.

We write $Im(\phi)$ and $Dom(\phi)$
for the image and domain of a function $\phi$.

We write $\psi \sub \phi$ for functions $\psi$ and $\phi$
if $\psi$ is a restriction of $\phi$.
Then $\phi \sm \psi$ denotes the restriction
of $\phi$ to $Dom(\phi) \sm Dom(\psi)$.
Given functions $\phi_j$ on $A_j$ for $j=1,2$ that agree on $A_1 \cap A_2$
we write $\phi_1 \cup \phi_2$ for the function on $A_1 \cup A_2$
that restricts to $\phi_j$ on $A_j$ for $j=1,2$.

Given $\GG \sub \mb{R}$, the $\GG$-span of $S \sub \mb{R}^d$ is 
$\bgen{S}_\GG = \{ \sum_{x \in S} \Phi_x x : \Phi \in \GG^S \}$.
We write $\bgen{S} = \bgen{S}_{\mb{Z}}$.

If $u \in (\mb{Z}^D)^X$ we write
$|u|=\sum_{x \in X} \sum_{d \in [D]} |(u_x)_d|$.

\section{Basic structures}

In this section we define the basic objects
needed for the statement of our main theorem
and record some simple properties of them 
that will be used throughout the paper.
We start in the first subsection with
the labelled complex structure that is
the functional analogue of a simplicial complex.
The second subsection concerns embeddings and extensions
of labelled complexes, and defines the extendability
property mentioned in the introduction.
In the third subsection we consider adapted complexes,
in which we add the structure of a permutation group
that acts on labellings: the orbits of this action 
play the role of edges in the example of hypergraph decomposition.
Then in the fourth subsection we formalise our general
decomposition problem with respect to a superimposed 
system of vector values on functions;
here we also show how to realise several
concrete examples within this general framework. 
We describe some basic properties of vector-valued
decompositions in the fifth subsection and introduce
some terminology (atoms and types) for them;
we also define the regularity property that
formalises the `no fractional obstacle' 
assumption discussed above.

\subsection{Complexes}

We start by defining a structure 
that we call a labelled complex.\footnote{
We suppress the term `labelled'
in our terminology, as the labels are 
indicated by the labelling set $R$.}

\begin{defn} \label{def:complex} (labelled complexes)
We call $\Phi=(\Phi_B: B \sub R)$ an $R$-system on $V$ 
if $\phi:B \to V$ is injective for each $\phi \in \Phi_B$.
We call an $R$-system $\Phi$ an $R$-complex on $V$
if whenever $\phi \in \Phi_B$ and $B' \sub B$ 
we have $\phi\mid_{B'} \in \Phi_{B'}$.
Let $\Phi^\circ_B = \{ \phi(B): \phi \in \Phi_B \}$,
$\Phi^\circ_j = \bigcup \{ \Phi^\circ_B: B \in \tbinom{R}{j} \}$
and $\Phi^\circ = \bigcup \{ \Phi^\circ_B: B \sub R\}$.
We write $V(\Phi) = \Phi^\circ_1$.
\end{defn}

Note that if $A \sub A' \in \Phi^\circ_{B'}$ then $A \in \Phi^\circ_B$ 
for some $B \sub B'$ (not necessarily unique);
thus $\Phi^\circ$ is a (simplicial) complex.
We will now define some basic operations
(forming restrictions and neighbourhoods)
for working with labelled complexes.

\begin{defn} \label{def:restrict} (restriction)
Let $\Phi$ be an $R$-complex and $\Phi'$ an $R$-system.
We let $\Phi[\Phi']$ be the $R$-system
where $\Phi[\Phi']_B$ is the set of $\phi \in \Phi_B$
such that $\phi\mid_{B'} \in \Phi'_{B'}$ 
for all $B' \sub B$ such that $\Phi'_{B'}$ is defined
(we allow some $\Phi'_{B'}$ to be undefined).
If $\Phi'_x = U \sub V(\Phi)$ for all $x \in R$ 
and $\Phi'_B$ is undefined otherwise
then we also write $\Phi[\Phi']=\Phi[U]
= \{ \phi \in \Phi: Im(\phi) \sub U \}$.  
\end{defn}

\begin{lemma} \label{restrict}
$\Phi[\Phi']$ is an $R$-complex.
\end{lemma}

\nib{Proof.}
Consider $\phi \in \Phi[\Phi']_B$ and
$B^* \sub B$ such that $\Phi'_{B^*}$ is defined.
We need to show $\phi\mid_{B^*} \in \Phi[\Phi']_{B^*}$.
To see this, note that $\phi\mid_{B^*} \in \Phi_{B^*}$,
and if $B' \sub B^*$ is such that $\Phi'_{B'}$ is defined 
then $(\phi\mid_{B^*})\mid_{B'}=\phi\mid_{B'} \in \Phi'_{B'}$
as $\phi \in \Phi[\Phi']$. \qed

\medskip

\begin{defn}
Let $\Phi$ be an $R$-complex and $\phi^* \in \Phi$.
We write $\Phi\mid_{\phi^*} 
= \{ \phi \in \Phi: \phi^* \sub \phi \}$.
\end{defn}

\begin{defn} \label{def:nhood} (neighbourhoods)
Let $\Phi$ be an $R$-complex and $\phi^* \in \Phi_{B^*}$.
For $\phi \in (\Phi\mid_{\phi^*})_{B \cup B^*}$
with $B \sub R \sm B^*$ let $\phi/\phi^* = \phi\mid_B$.
We define an $(R \sm B^*)$-system $\Phi/\phi^*$ 
where each $(\Phi/\phi^*)_B$ consists of all
$\phi/\phi^*$ with $\phi \in (\Phi\mid_{\phi^*})_{B \cup B^*}$.
For $J \in \GG^\Phi$ we define $J/\phi^* \in \GG^{\Phi/\phi^*}$
by $(J/\phi^*)_{\phi/\phi^*}=J_\phi$ whenever $\phi^* \sub \phi$. 
\end{defn}

\begin{lemma} \label{nhood}
$\Phi/\phi^*$ is an $R \sm B^*$-complex.
\end{lemma}

\nib{Proof.}
Consider $B' \sub B \sub R \sm B^*$
and $\phi/\phi^* \in (\Phi/\phi^*)_B$.
We need to show 
$(\phi/\phi^*)\mid_{B'} \in (\Phi/\phi^*)_{B'}$.
This holds as 
$\phi':=\phi\mid_{B^* \cup B'} \in \Phi_{B^* \cup B'}$
with $\phi'\mid_{B^*}=\phi^*$. \qed

\subsection{Embeddings and extensions}

Here we formulate our extendability property
and show that it is maintained 
under taking neighbourhoods.
We start by defining embeddings
of labelled complexes.

\begin{defn}
Let $H$ and $\Phi$ be $R$-complexes.
Suppose $\phi:V(H) \to V(\Phi)$ is injective.
We call $\phi$ a $\Phi$-embedding of $H$
if $\phi \circ \psi \in \Phi$ for all $\psi \in H$.
\end{defn}

We will define extendability using
the following labelled complex of partite maps.

\begin{defn}
Let $R(S)$ be the $R$-complex of all 
partite maps from $R$ to $R \times S$, i.e.\ 
whenever $i \in B \sub R$ and $\psi \in R(S)_B$ 
we have $\psi(i)=(i,x)$ for some $x \in S$.
If $S=[s]$ we write $R(S)=R(s)$.
\end{defn}

The following extendability property 
can be viewed as a labelled analogue
of that in \cite{Kexist}.\footnote{
The unlabelled analogue of our assumption 
here is weaker than that in \cite{Kexist},
as we only consider partite extensions.}

\begin{defn} \label{def:ext}
Suppose $H \sub R(S)$ is an $R$-complex
and $F \sub V(H)$. Define $H[F] \sub R(S)$ by
$H[F] = \{ \psi \in H: Im(\psi) \sub F \}$.
Suppose $\phi$ is a $\Phi$-embedding of $H[F]$.
We call $E=(H,F,\phi)$ a $\Phi$-extension of rank $s=|S|$.
We say $E$ is simple if $|V(H) \sm F|=1$.

We write $X_E(\Phi)$ for the set or number of 
$\Phi$-embeddings of $H$ that restrict to $\phi$ on $F$.
We say $E$ is $\oO$-dense (in $\Phi$)
if $X_E(\Phi) \ge \oO |V(\Phi)|^{v_E}$,
where $v_E := |V(H) \sm F|$.
We say $\Phi$ is $(\oO,s)$-extendable
if all $\Phi$-extensions of rank $s$ are $\oO$-dense.
\end{defn}
 
We will also require the following extension
of the previous definition that allows 
for a system of extra restrictions.

\begin{defn} \label{def:ext+}
Let $\Phi$ be an $R$-complex and
$\Phi' = (\Phi^t: t \in T)$ 
with each $\Phi^t \sub \Phi$.
Let $E=(H,F,\phi)$ be a $\Phi$-extension 
and $H' = (H^t: t \in T)$ for some
mutually disjoint $H^t \sub H \sm H[F]$;
we call $(E,H')$ a $(\Phi,\Phi')$-extension.

We write $X_{E,H'}(\Phi,\Phi')$ for the set or number 
of $\phi^* \in X_E(\Phi)$ with $\phi^* \circ \psi \in \Phi^t_B$ 
whenever $\psi \in H^t_B$ and $\Phi^t_B$ is defined.
We say $(E,H')$ is $\oO$-dense in $(\Phi,\Phi')$
if $X_{E,H'}(\Phi,\Phi') \ge \oO |V(\Phi)|^{v_E}$.
We say $(\Phi,\Phi')$ is $(\oO,s)$-extendable if all 
$(\Phi,\Phi')$-extensions of rank $s$ are $\oO$-dense in $(\Phi,\Phi')$.

When $|T|=1$ we identify $\Phi' \sub \Phi$ with $(\Phi')$.
We also write $X_E(\Phi,\Phi')=X_{E,H \sm H[F]}(\Phi,\Phi')$.
For $L \sub \Phi^\circ$ we write 
$X_E(\Phi,L)=X_E(\Phi,\Phi')$ where
$\Phi'=\{\phi \in \Phi: Im(\phi) \in L\}$;
we also write $\Phi[L]=\Phi[\Phi']$,
and say that $(\Phi,L)$ is $(\oO,s)$-extendable
if $(\Phi,\Phi')$ is $(\oO,s)$-extendable.

If $L \sub V(\Phi)$ we also say that
$(\Phi,L)$ is $(\oO,s)$-extendable wrt $L$ 
if $X_{E,H'}(\Phi,\Phi') \ge \oO |L|^{v_E}$
for all $(\Phi,\Phi')$-extensions $(E,H')$ of rank $s$.
\end{defn}

Note that if $\Phi' \sub \Phi$ and
$(\Phi,\Phi')$ is $(\oO,s)$-extendable
then $\Phi[\Phi']$ is $(\oO,s)$-extendable.
In the next definition we combine the operations
of taking neighbourhoods and restriction to an
(unordered) hypergraph; the accompanying lemma
shows that under the generalised extendability
condition of the previous condition
the resulting labelled complex is extendable.
We note that if $L=\Phi^\circ_r$ 
the restriction has no effect,
so $\Phi/\phi^*L=\Phi/\phi^*$,
and in this case Lemma \ref{nhood:ext}
states that if $\Phi$ is $(\oO,s)$-extendable 
then $\Phi/\phi^*$ is $(\oO,s)$-extendable.
A less trivial example is when 
$L \sub V(\Phi) = \Phi^\circ_1$;
then $\Phi/\phi^*L = (\Phi/\phi^*)[L]$
is obtained by restricting $\Phi/\phi^*$ to $L$.

\begin{defn}
Suppose $\Phi$ is an $R$-complex, 
$L \sub \Phi^\circ_r$ and $\phi^* \in \Phi_{B^*}$.
Let $\Phi/\phi^*L$ be the set of all
$\phi/\phi^* \in \Phi/\phi^*$ such that
$e \in L$ for all $e \in Im(\phi)_r$
with $e \sm Im(\phi^*) \ne \es$.
\end{defn}

\begin{lemma} \label{nhood:ext}
If $(\Phi,L)$ is $(\oO,s)$-extendable (wrt $L$)
then $\Phi/\phi^*L$ is $(\oO,s)$-extendable.
\end{lemma}

\nib{Proof.}
Consider any $\Phi/\phi^*L$-extension 
$E = (H,F,\phi)$ where $H \sub (R \sm B^*)(s)$.
We need to show that $E$ is $\oO$-dense in $\Phi/\phi^*L$.
To see this, we consider $H^+ \sub R(s)$ where for each 
$B \sub R \sm B^*$, $B' \sub B^*$, $\psi \in H_B$ 
we include $\psi \cup B'$ in $H^+_{B \cup B'}$
defined by $(\psi \cup B')\mid_B = \psi$
and $(\psi \cup B')(i)=i=(i,1)$ for $i \in B'$.
Consider the $\Phi$-extension 
$E^+ = (H^+,F^+,\phi \cup \phi^*)$
with $F^+ = F \cup B^*$.

As $(\Phi,L)$ is $(\oO,s)$-extendable, we have 
$X_{E^+}(\Phi,L) > \oO |V(\Phi)|^{v_{E^+}}$,
or $X_{E^+}(\Phi,L) > \oO |L|^{v_{E^+}}$
if $(\Phi,L)$ is $(\oO,s)$-extendable wrt $L$.
It remains to show that if $\phi^+ \in X_{E^+}(\Phi,L)$ 
then $\phi^+\mid_{V(H)} \in X_E(\Phi/\phi^*L)$.
For any $B \sub R \sm B^*$, $\psi \in H_B$
as $\phi^+ \in X_{E^+}(\Phi)$ we have 
$\phi' := \phi^+ \circ (\psi \cup B^*) \in \Phi_{B \cup B^*}$,
and as $\phi' \mid_{B^*} = \phi^*$ 
we have $\phi^+\mid_{V(H)} \circ \psi
= \phi' \mid_B = \phi'/\phi^* \in (\Phi/\phi^*)_B$. 
Furthermore, as $\phi^+ \in X_{E^+}(\Phi,L)$
we have $e \in L$ for any $e \in Im(\phi^+)_r$ 
with $e \sm Im(\phi^*) \ne \es$.
Therefore $\phi^+\mid_{V(H)} \in X_E(\Phi/\phi^*L)$. \qed

\subsection{Adapted complexes}

Next we introduce the setting of adapted complexes,
where we have a permutation group acting 
on the functions in a labelled complex.
We start with some notation for permutation groups;
in particular, given a permutation group $\Ss$ on $R$
we define an $R$-complex $\Ss^\le$ that consists
of all restrictions of elements of $\Ss$.

\begin{defn} \label{def:perm} 
Suppose $\Ss$ is a permutation group on $R$.
For $B,B' \sub R$ we write 
$\Ss^{B'}_B = \{ \sS\mid_B: \sS \in \Ss, \sS(B)=B' \}$,
$\Ss_B = \cup_{B'} \Ss^{B'}_B$,
$\Ss^{B'} = \cup_B \Ss^{B'}_B$,
$\Ss[B] = \cup_{B' \sub B} \Ss^{B'}$,
$\Ss^\le = \cup_{B,B'} \Ss^{B'}_B$.

We let $\mc{P}^\Ss$ be the equivalence classes of the
relation $B \sim B' \lra \Ss^B_{B'} \ne \es$.
Note that $B \sim B'$ implies $|B|=|B'|$.
We write $\mc{P}^\Ss_j = 
\{ C \in \mc{P}^\Ss: B \in C \Ra |B| = j \}$.
\end{defn}

We will restrict attention to labelled complexes 
in which any function can be relabelled under the
group action, as follows.

\begin{defn} \label{def:adapt} (adapted)
Suppose $\Phi$ is an $R$-complex
and $\Ss$ is a permutation group on $R$.
For $\sS \in \Ss$ and $\phi \in \Phi_{\sS(B)}$
let $\phi \sS = \phi \circ \sS\mid_B$.
We say $\Phi$ is $\Ss$-adapted if 
$\phi \sS \in \Phi$ for any $\phi \in \Phi$, $\sS \in \Ss$.
\end{defn}

Next we introduce some notation for the orbits
of the action implicit in the previous definition;
these will play the role of edges
in hypergraph decompositions.

\begin{defn} \label{def:orbit} (orbits)

For $\psi \in \Phi_B$ with $B \sub R$ 
we define the orbit of $\psi$ by
$\psi\Ss := \psi\Ss^B = \{ \psi\sS: \sS \in \Ss^B \}$.
We denote the set of orbits by $\Phi/\Ss$.
We write $\Phi_C = \cup_{B \in C} \Phi_B$ for $C \in \mc{P}^\Ss$.
We write $Im(O)=Im(\psi)$ for $\psi \in O \in \Phi/\Ss$.
For $O,O' \in \Phi/\Ss$ we write $O \sub O'$ if there are 
$\psi \in O$, $\psi' \in O'$ with $\psi \sub \psi'$.
\end{defn}

Note that the orbits partition $\Phi$ and
$\Phi_C = \bigcup \{\psi\Ss: \psi \in \Phi_B\}$ for any $B \in C$.
When we later consider functions on $\Phi$
we will decompose them by orbits as follows.

\begin{defn} (orbit decomposition)
Let $\GG$ be an abelian group.
For $J \in \GG^{\Phi_r}$ and $O \in \Phi_r/\Ss$
we define $J^O$ by $J^O_\psi = J_\psi 1_{\psi \in O}$.
The orbit decomposition of $J$ is
$J = \sum_{O \in \Phi_r/\Ss} J^O$.
\end{defn}

Now we will illustrate the role of orbits with the two
most obvious examples (see also subsection 
\ref{sec:vvd} for more examples).

\nib{Examples.}\
\begin{enumerate}
\item If $\Ss=\{id_R\}$ is the trivial group then
each equivalence class and orbit has size $1$,
and we can identify $\Phi$ with $\Phi/\Ss$.
This choice of $\Ss$ is suitable for `fully partite'
hypergraph decompositions, in which every edge 
is uniquely labelled by the set of parts that it meets.
We also denote $\Ss^\le$ by $\ova{R}$,
or by $\ova{q}$ when $R=[q]$.
Then $\ova{q}_B=\{id_B\}$ for all $B \sub [q]$.
\item If $\Ss$ is the symmetric group $S_R$ on $R$
then the equivalence classes of $\mc{P}^\Ss$ 
are $\tbinom{R}{r}$ for $0 \le r \le |R|$.
We also denote $\Ss^\le$ by $R^\le$.
Then each $R^\le_B = Inj(B,R)$ consists of
all injections from $B$ to $R$.
We can identify $\Phi^\circ$ with $\Phi/\Ss$,
where $e \in \Phi^\circ$ is identified
with $\{ \psi \in \Phi: Im(\psi)=e \}$.
This choice of $\Ss$ is suitable for nonpartite
hypergraph decompositions, in which
the labels play no essential role.
\end{enumerate}

\medskip

Next we show that adapted complexes 
have neighbourhoods that are also adapted complexes.

\begin{defn} \label{nhood:perm}
Let $\Ss$ be a permutation group on $R$.
For $B^* \sub R$ and $\sS \in \Ss$ with
$\sS\mid_{B^*}=id_{B^*}$ we write
$\sS/B^* = \sS\mid_{R \sm B^*}$.
We let $\Ss/B^*$ be the set of all such $\sS/B^*$.
\end{defn}

Note that $\Ss/B^*$ is a permutation group on $R \sm B^*$.

\begin{lemma} \label{nhood:adapt}
Let $\Phi$ be a $\Ss$-adapted $R$-complex and $\phi^* \in \Phi_{B^*}$.
Then $\Phi/\phi^*$ is a $\Ss/B^*$-adapted $(R \sm B^*)$-complex. 
\end{lemma}

\nib{Proof.} 
Suppose $B \sub R \sm B^*$, $\sS=\sS'/B^* \in \Ss/B^*$
and $\psi = \psi'/\phi^* \in (\Phi/\phi^*)_B$.
As $\Phi$ is $\Ss$-adapted, $\psi'\sS' \in \Phi$, 
so $\psi\sS = \psi'\sS'/\phi^* \in \Phi/\phi^*$. \qed

\medskip

Next we introduce the labelled complex structure
defined by embeddings of one labelled complex in another.

\begin{defn} \label{def:APhi}
Given $R$-complexes $\Phi$ and $A$ we let $A(\Phi)$ 
denote the set of $\Phi$-embeddings of $A$.
We let $A(\Phi)^\le$ denote the $V(A)$-complex
where each $A(\Phi)^\le_F$ for $F \sub V(A)$
is the set of $\Phi$-embeddings of $A[F]$.
\end{defn}

In the next subsection we will apply 
Definition \ref{def:APhi} with $A = \Ss^\le$;
we conclude this subsection by showing that
if $\Phi$ is $\Ss$-adapted then we can identify
the resulting complex of embeddings with $\Phi$ itself.

\begin{lemma} \label{Phi=APhi}
If $\Phi$ is $\Ss$-adapted and $B \sub [q]$
then $\Ss[B](\Phi) = \Phi_B$.
\end{lemma}

\nib{Proof.}
Consider any $\phi \in \Phi_B$.
As $\Phi$ is $\Ss$-adapted,
for any $\sS \in \Ss^B$ we have
$\phi\sS \in \Phi$.
As $\Phi$ is a $[q]$-complex
we deduce $\phi\sS \in \Phi$
for any $\sS \in \Ss[B]$,
so $\phi \in \Ss[B](\Phi)$.
Conversely, if $\phi \in \Ss[B](\Phi)$
then $\phi = \phi\ id_B \in \Phi_B$. \qed

\subsection{Vector-valued decompositions} \label{sec:vvd}

Now we introduce
our general framework for decomposing vectors 
with coordinates indexed by the functions
of a labelled complex and entries in some abelian group.
We follow the definition with several examples
to show how it captures hypergraph decompositions
and other related problems. 

\begin{defn} \label{def:vsys} 
Let $\mc{A}$ be a set of $R$-complexes;
we call $\mc{A}$ an $R$-complex family.
If each $A \in \mc{A}$ is a copy of $\Ss^\le$
we call $\mc{A}$ a $\Ss^\le$-family.
For $r \in \mb{N}$ we write 
$A_r = \bigcup \{ A_B: B \in \tbinom{R}{r} \}$
and $\mc{A}_r = \cup_{A \in \mc{A}} A_r$.
We let $\mc{A}(\Phi)^\le$ denote the 
$V(A)$-complex family $(A(\Phi)^\le: A \in \mc{A})$.

Let $\gG \in \GG^{\mc{A}_r}$ for some abelian group $\GG$;
we call $\gG$ a $\GG$-system for $\mc{A}_r$.

Let $\Phi$ be an $R$-complex.
For $\phi \in A(\Phi)^\le$ with $A \in \mc{A}$
we define $\gG(\phi) \in \GG^{\Phi_r}$ 
by $\gG(\phi)_{\phi \circ \tT} = \gG_\tT$ for $\tT \in A_r$
(zero otherwise).
We call $\gG(\phi)$ a $\gG$-molecule
and let $\gG(\Phi)$ be the set of $\gG$-molecules.

Given $\Psi \in \mb{Z}^{\mc{A}(\Phi)}$ 
we define $\pl \Psi = \pl^\gG \Psi
= \sum_\phi \Psi_\phi \gG(\phi) \in \GG^{\Phi_r}$.
We also call $\Psi$ an integral 
$\gG(\Phi)$-decomposition of $G = \pl \Psi$
and call $\bgen{\gG(\Phi)}$ the decomposition lattice.
If furthermore $\Psi \in \{0,1\}^{\mc{A}(\Phi)}$
(i.e.\ $\Psi \sub \mc{A}(\Phi)$) we call $\Psi$ 
a $\gG(\Phi)$-decomposition of $G$.
\end{defn}

\nib{Examples.}
\begin{enumerate}
\item Suppose $H$ and $G$ are $r$-graphs with $V(H)=[q]$.
Let $\Phi$ be the complete $[q]$-complex on $V(G)$, 
i.e.\ all $\Phi_B = Inj(B,V(G))$.
Let $G^* = \{ \psi \in \Phi_r: Im(\psi) \in G\}$.
Let $\Ss=S_{[q]}$, $\mc{A}=\{A\}$ with $A=\Ss^\le$,
and $\gG \in \{0,1\}^{A_r}$ with
each $\gG_\tT = 1_{Im(\tT) \in H}$.
Note that $\Phi$ is $\Ss$-adapted.
For any $\phi \in A(\Phi) = \Phi_q$ 
and $\tT \in A_r$ we have
$\gG(\phi)_{\phi\tT} = \gG_\tT = 1_{Im(\tT) \in H}$,
so an (integral) $H$-decomposition of $G$ is equivalent 
to a (an integral) $\gG(\Phi)$-decomposition of $G^*$.
Similarly, if $\mc{H}$ is a family of $r$-graphs,
by introducing isolated vertices we may assume
they all have vertex set $[q]$.
Then an $\mc{H}$-decomposition of $G$ is equivalent 
to a $\gG(\Phi)$-decomposition of $G^*$,
where now $\mc{A}$ contains $A^H$ defined
as above for each $H \in \mc{H}$,
and $\gG_\tT$ is as defined above
whenever $\tT \in A^H_r$.
\item We generalise the previous example
(for simplicity we revert to one $r$-graph $H$).
Now suppose $H$ and $G$ have coloured edges.
Let the set of colours be $[D]$,
let $H^d$ and $G^d$ be the edges 
in $H$ and $G$ of colour $d$.
Let $e_1,\dots,e_D$ be the 
standard basis of $\mb{Z}^D$.
Let $\Phi$, $\Ss$, $\mc{A}$ be as above.
Define $G^* \in (\mb{N}^D)^{\Phi_r}$
by $G^*_\psi = e_d$ for all $\psi$
with $Im(\psi) \in G^d$
and $G^*_\psi = 0$ otherwise.
Define $\gG \in (\mb{N}^D)^{A_r}$ by
$\gG_\tT = e_d$ for all $\tT \in A_r$
with $Im(\psi) \in H^d$
and $\gG_\tT = 0$ otherwise.
Then an $H$-decomposition of $G$ 
that respects colours is equivalent 
to a $\gG(\Phi)$-decomposition of $G^*$.
\item Now suppose that $G$ is $q$-partite,
say with parts $V_1,\dots,V_q$.
The previous examples will not be useful
for finding an $H$-decomposition of $G$,
as our main theorem requires $(\Phi,G)$
to be extendable, but if we allow $\Phi$
to disrespect the partition then we cannot
extend all partial embeddings within $G$.
Instead, we define $\Phi_B$ for $B \sub [q]$
to consist of all partite $\psi \in Inj(B,V(G))$,
i.e.\ $\psi(i) \in V_i$ for all $i \in B$.
We let $\Ss = \{id\}$ be trivial,
$\mc{A}=\{A\}$ with $A=\Ss^\le$,
i.e.\ all $A_B = \{id_B\}$.
Defining $G^*$ and $\gG$ as in the first example,
we again see that an $H$-decomposition of $G$ 
is equivalent to a $\gG(\Phi)$-decomposition of $G^*$.
\item Next we consider the $H$-decomposition problem
for $G$ when we are given bipartitions $(X,Y)$ of $V(G)$
and $(A,B)$ of $V(H)=[q]$, and we only allow copies of $H$
in which $A$ maps into $X$ and $B$ into $Y$
(recall that the problems of resolvable designs
and large sets of designs are equivalent
to such bipartite decomposition problems).
We let $\Phi_F$ for $F \sub [q]$ consist
of all $\psi \in Inj(F,V(G))$ such that
$\psi(F \cap A) \sub X$ and $\psi(F \cap B) \sub Y$.
As usual, we let 
$G^* = \{ \psi \in \Phi_r: Im(\psi) \in G\}$.
We let $\Ss$ be the group of all $\sS \in S_q$
such that $\sS(A)=A$ and $\sS(B)=B$.
Then $\Phi$ is $\Ss$-adapted.
As usual, we let
$\mc{A} = \{A'\}$ with $A'=\Ss^\le$
and $\gG \in \{0,1\}^{A'_r}$ with
each $\gG_\tT = 1_{Im(\tT) \in H}$.
Then an $H$-decomposition of $G$ is equivalent 
to a $\gG(\Phi)$-decomposition of $G^*$.
\item 
In the above examples we were decomposing
hypergraphs (sets of sets) and treating the labellings
(sets of functions) as a convenient device,
but many applications explicitly require labellings.
An example that may have some topological motivation
is that of decomposing the set of top-dimensional cells 
of an oriented simplicial complex.
The standard definition of orientations
fits very well with our framework: 
for $r$-graphs $H$ and $G$, an orientation is defined
by a bijective labelling of each edge by $[r]$,
where two labellings are considered equivalent
if they differ by an even permutation in $S_r$.
Then we wish to decompose $G$ by copies of $H$,
where we only allow copies $\phi(H)$ such that
for each edge $e$ of $H$ composing the labelling
of $e$ with $\phi$ gives a labelling of $\phi(e)$
equivalent to that in $G$.
To realise this problem in our framework
(consider for simplicity the nonpartite setting 
where $\Phi$ is the complete $[q]$-complex on $V(G)$
and $\Ss=S_q$), for each $B \in Q=[q]_r$
we let $\pi_B \in Bij([r],B)$ be order preserving
and let $G^*_B$ consist of all $\psi \in \Phi_B$ 
such that $\psi' = \psi \circ \pi_B$ with
$Im(\psi') \in G$ is correctly oriented.
Similarly, we let $\mc{A} = \{A\}$ with $A=\Ss^\le$
and $\gG \in \{0,1\}^{A_r}$ where for $\tT \in A_B$
we let $\gG_\tT$ be $1$ if $\tT' = \tT \circ \pi_B$
with $Im(\tT') \in H$ is correctly oriented,
otherwise $\gG_\tT=0$. Then an oriented $H$-decomposition 
of $G$ is equivalent to a $\gG(\Phi)$-decomposition of $G^*$.
\end{enumerate}

\subsection{Atoms and types}

In this subsection we introduce some structures
and terminology for working with vector-valued 
decompositions, and make some preliminary observations
regarding the decomposition lattice $\bgen{\gG(\Phi)}$.
We also define the regularity property referred to above.
Throughout we let $\Ss \le S_q$ be a permutation group,
$\Phi$ be a $\Ss$-adapted $[q]$-complex,
$\mc{A}$ be a $\Ss^\le$-family and
$\gG \in \GG^{\mc{A}_r}$.

\begin{defn} \label{def:atom} (atoms)
For any $\phi \in \mc{A}(\Phi)$ and $O \in \Phi_r/\Ss$
such that $\gG(\phi)^O \ne 0$ we call
$\gG(\phi)^O$ a $\gG$-atom at $O$.
We write $\gG[O]$ for the set of $\gG$-atoms at $O$.
We say $\gG$ is elementary if
all $\gG$-atoms are linearly independent.
We define a partial order $\le_\gG$
on $\GG^{\Phi_r}$ where $H \le_\gG G$
iff $G-H$ can be expressed as the sum
of a multiset of $\gG$-atoms.
\end{defn}

Note that if $J \in \bgen{\gG(\Phi)}$ 
then each $J^O$ can be expressed as a 
$\mb{Z}$-linear combination of $\gG$-atoms at $O$.
Furthermore, if $\gG$ is elementary
then this expression is unique, 
so if $J$ is the sum of a multiset $Z$ of $\gG$-atoms
then a $\gG(\Phi)$-decomposition of $J$
may be thought of as a partition of $Z$,
where each part is the set of $\gG$-atoms
contained in some molecule $\gG(\phi)$.
This is a combinatorially natural condition,
as it avoids arithmetic issues that arise 
e.g.\ for decompositions of integers
(the Frobenius coin problem).
In our main theorem we will assume 
that $\gG$ is elementary, but the proof also uses other vector systems 
derived from $\gG$ that are not necessarily elementary.

The following definition and accompanying lemma
give various equivalent ways to represent atoms.
The notation $\gG(\phi)$ matches the notation for molecules 
in Definition \ref{def:vsys} when $\phi \in \mc{A}(\Phi)$.

\begin{defn} \label{def:atom=}
For $\psi \in \Phi_B$ and $\tT \in \mc{A}_B$
we define $\gG[\psi]^\tT \in \GG^{\psi\Ss}$
by $\gG[\psi]^\tT_{\psi\sS} = \gG_{\tT\sS}$.

For $\phi \in A(\Phi)^\le = \Phi$ we define
$\gG(\phi) \in \GG^{\Phi_r}$ by
$\gG(\phi)_{\phi \tT} = \gG_\tT$ 
whenever $\tT \in A_r$ with $Im(\tT) \sub Dom(\phi)$.
\end{defn}

\begin{lemma} \label{atom=} 
Suppose $\phi \in A(\Phi)$ and
$\psi \in A(\Phi)^\le_B = \Phi_B$.
\begin{enumerate}
\item 
If $\psi = \phi \tT$ with $\tT \in A_r$
then $\gG(\phi)^{\psi\Ss} = \gG[\psi]^\tT$. \\
Furthermore, if $\tT \in A_B$ and $\sS \in \Ss^B$ 
then $\gG[\psi]^\tT = \gG[\psi\sS]^{\tT\sS}$.
\item 
If $\psi \sub \phi$ then 
$\gG(\phi)^{\psi\Ss} = \gG[\psi]^{id_B} = \gG(\psi)$.
\end{enumerate}
\end{lemma}

\nib{Proof.}
For (i), by Definitions
\ref{def:vsys} and \ref{def:atom=},
for any $\sS \in \Ss^B$ we have 
$\gG(\phi)_{\psi\sS} = \gG(\phi)_{\phi\tT\sS}
= \gG_{\tT\sS} = \gG[\psi]^\tT_{\psi\sS}$,
i.e.\ $\gG(\phi)^{\psi\Ss} = \gG[\psi]^\tT$.
Furthermore, if $\tT \in A_B$, $\sS \in \Ss^B_{B'}$,
$\sS' \in \Ss^{B'}$ then 
$\gG[\psi]^\tT_{\psi\sS\sS'} = \gG_{\tT\sS\sS'}
= \gG[\psi\sS]^{\tT\sS}_{\psi\sS\sS'}$.
For (ii), we have $\psi = \phi\ id_B$, so 
$\gG(\phi)^{\psi\Ss} = \gG[\psi]^{id_B}$ by (i).
Also, for any $\sS \in \Ss^B$ we have 
$\gG(\psi)_{\psi\sS} = \gG_\sS = \gG[\psi]^{id_B}_{\psi\sS}$,
so $\gG[\psi]^{id_B} = \gG(\psi)$. \qed

\medskip

The following definition will be used for the 
extendability assumption on $(\Phi,\gG[G])$
in our main theorem, which gives a lower bound
on extensions such that all atoms belong to $G$.

\begin{defn} \label{Gatoms}
For $G \in \GG^{\Phi_r}$ we let
$\gG[G]=(\gG[G]^A: A \in \mc{A})$
where each $\gG[G]^A$ is the set of
$\psi \in A(\Phi)^\le_r = \Phi_r$ 
such that $\gG(\psi) \le_\gG G$.
\end{defn}

When using the notation $\gG[\psi]^\tT$ for an atom,
there may be several choices of $\tT$ that give rise to
the same atom; this defines an equivalence relation 
that we will call a type. To illustrate the following definition,
we recall example $i$ (nonpartite $H$-decomposition)
from subsection \ref{sec:vvd}. In this case, there are
two types for each $r$-set $B$ of labels:
for any $\tT \in A_B$, if $Im(\tT) \in H$
then $\gG^\tT$ is the all-1 vector
(we think of this type as an edge),
whereas $Im(\tT) \notin H$
then $\gG^\tT$ is the all-0 vector
(the zero type, which we think of as a `non-edge').

\begin{defn} \label{def:type} (types)
For $\tT \in \mc{A}_B$ with $B \in Q$ we define 
$\gG^\tT \in \GG^{\Ss^B}$ by $\gG^\tT_\sS = \gG_{\tT\sS}$.

A type $t=[\tT]$ in $\gG$ is an equivalence class 
of the relation $\sim$ on any $\mc{A}_B$ with $B \in Q$
where $\tT \sim \tT'$ iff $\gG^\tT = \gG^{\tT'}$.
We write $T_B$ for the set of types in $\mc{A}_B$.

For $\tT \in t \in T_B$ and $\psi \in \Phi_B$ we write 
$\gG^t = \gG^\tT$ and $\gG[\psi]^t = \gG[\psi]^\tT$.

If $\gG^t=0$ call $t$ a zero type and write $t=0$.

If $\phi \in \mc{A}(\Phi)$ with
$\gG(\phi)^{\psi\Ss} = \gG[\psi]^t$
we write $t_\phi(\psi) = t$.
\end{defn}

The next lemma shows that $\gG[\psi]^t$ is well-defined.

\begin{defn} \label{def:fB}
For $B \in C \in \mc{P}^\Ss_r$ and $J \in \GG^{\Phi_C}$
we define $f_B(J) \in (\GG^{\Ss^B})^{\Phi_B}$ 
by $(f_B(J)_\psi)_\sS = J_{\psi\sS}$.
\end{defn}

\begin{lemma} \label{fB}
If $J = \gG[\psi]^t$ for some 
$\psi \in \Phi_B$, $t \in T_B$
then $f_B(J)_\psi = \gG^t$.
\end{lemma}

\nib{Proof.}
For any $\sS \in \Ss^B$ and $\tT \in t$ we have
$(f_B(J)_\psi)_\sS = J_{\psi\sS} = \gG[\psi]^\tT_{\psi\sS}
= \gG_{\tT\sS} = \gG^\tT_\sS = \gG^t_\sS$. \qed

\medskip

We also see from Lemma \ref{fB} that $\gG$ is elementary 
iff for any $B \in Q$ the set of nonzero $\gG^t$ 
with $t \in T_B$ is linearly independent.
Next we introduce certain group actions
that will be important in section \ref{sec:int};
the following lemma records their effect on types.

\begin{defn} \label{def:action} 
For any set $X$ we define a right $\Ss^B_B$ action on $X^{\Ss^B}$ 
by $(v \tau)_\sS = v_{\tau\sS}$ whenever 
$v \in X^{\Ss^B}$, $\tau \in \Ss^B_B$, $\sS \in \Ss^B$.
\end{defn}

Note that Definition \ref{def:action} is indeed a right action,
as for $\tau_1,\tau_2 \in \Ss^B_B$ we have
$((v \tau_1)\tau_2 )_\sS = (v \tau_1)_{\tau_2 \sS}
= v_{\tau_1\tau_2\sS} = (v (\tau_1\tau_2))_\sS$.
For future reference we also note the linearity
$(v+v') \tau = v \tau + v' \tau$; 
indeed, for $\sS \in \Ss^B$ we have
$((v+v') \tau)_\sS = (v+v')_{\tau\sS}
= v_{\tau\sS} + v'_{\tau\sS}
= (v \tau)_\sS + (v' \tau)_\sS$.

\begin{lemma} \label{taut}
If $\tT \in \mc{A}_B$ and $\tau \in \Ss^B_B$
then $\gG^\tT \tau = \gG^{\tT\tau}$.
\end{lemma}

\nib{Proof.}
For any $\sS \in \Ss^B$ we have
$(\gG^\tT \tau)_\sS = \gG^\tT_{\tau \sS}
  = \gG_{\tT \tau \sS} = \gG^{\tT \tau}_\sS$. \qed

\medskip

The next definition and accompanying lemma
restate and provide notation for the earlier
observation that any vector in the decomposition lattice
can be expressed as a $\mb{Z}$-linear combination of atoms
(we omit the trivial proof).

\begin{defn} \label{def:L-} 
Let $\mc{L}^-_\gG(\Phi)$ be the set of $J \in \GG^{\Phi_r}$ 
such that $J^O \in \sgen{\gG[O]}$ for all $O \in \Phi_r/\Ss$.
\end{defn}

\begin{lemma} \label{LsubL-}
$\bgen{\gG(\Phi)} \sub \mc{L}^-_\gG(\Phi)$.
\end{lemma}

Next we define two notions of symmetry,
one for vectors and the other for subsets.

\begin{defn} \label{def:symm}
We call $v \in (\GG^{\Ss^B})^{\Phi_B}$ symmetric
if $v_\psi \tau = v_{\psi\tau}$ whenever
$\psi \in \Phi_B$, $\tau \in \Ss^B_B$. 

We call $H \sub \GG^{\Ss^B}$ symmetric
if $g \tau \in H$ whenever $g \in H$, $\tau \in \Ss^B_B$.
\end{defn}

Note that 
\[ G^B := \{\gG^t: t \in T_B \} \ \text{ and } \
 \gG^B := \sgen{G^B} \le \GG^{\Ss^B} \]
are symmetric by Lemma \ref{taut}.
Now we use types to give an alternative description 
of the lattice from Definition \ref{def:L-}.

\begin{lemma} \label{gB}
Let $B \in C \in \mc{P}^\Ss_r$ and 
$J \in \GG^{\Phi_C} \cap \mc{L}^-_\gG(\Phi)$.
Then $f_B(J) \in (\gG^B)^{\Phi_B}$ is symmetric.
\end{lemma}

\nib{Proof.}
By linearity, we can assume $J$ is a $\gG$-atom, say 
$J=\gG[\psi]^\tT$ with $\psi \in \Phi_B$, $\tT \in \mc{A}_B$.
For any $\tau \in \Ss^B_B$ we have
$J = \gG[\psi\tau]^{\tT\tau}$ by Lemma \ref{atom=}.i 
and $f_B(J)_{\psi\tau} = \gG^{\tT\tau} = \gG^\tT \tau
= f_B(J)_\psi \tau$ by Lemmas \ref{fB},
\ref{taut} and \ref{fB} again. \qed

\medskip

For future reference (in section \ref{sec:cea})
we also note the following lemma which will allow
us to split any linear dependence of $\gG$-atoms
into constant sized pieces. If $\gG$ is elementary
then $Z_B(\gG)=\{0\}$ so there is nothing to prove;
in this case we let $C_0=1$.
We call $C_0$ the lattice constant.

\begin{lemma} \label{AAspan}
There is $C_0=C_0(\gG)$ such that for any $n \in Z_B(\gG) 
:= \{ n \in \mb{Z}^{\mc{A}_r}: \sum_\tT n_\tT \gG^\tT = 0 \}$,
there are $n^i \in Z_B(\gG)$ for $i \in [t]$ 
for some $t \le C_0 |n|$ with each $|n^i| \le C_0$
and $n = \sum_{i \in [t]} n^i$.
\end{lemma}

\nib{Proof.}
Let $\mc{X}$ be an integral basis for $Z_B(\gG)$.
Let $Z$ be the matrix with columns $\mc{X}$.
We will find an integral solution $v$ of $n=Zv$ 
and then for each $X \in \mc{X}$ take $|v_X|$ of 
the $n^i$ equal to $\pm X$ (with the sign of $v_X$).
The following explicit construction implies
the required bound for $|v|$ in terms of $|n|$.
We can put $Z$ in `diagonal form' via
elementary row and column operations:
there are unimodular
(integral and having integral inverses)
matrices $P$ and $Q$ such that $D=PZQ$
has $D_{ij} \ne 0 \Lra i=j$.
To solve $n = Zv$ we need to solve $Pn = DQ^{-1}v$.
Let $R$ be the set of nonzero rows of $D$
and let $(Pn)_R$ and $D_R$ denote the
corresponding restrictions of $Pn$ and $D$.
Then $D_R^{-1}(Pn)_R$ is integral 
(as $n \in \bgen{\mc{X}}$)
so $v = Q D_R^{-1}(Pn)_R$ 
is an integral solution of $n=Zv$. \qed

\medskip

Next we introduce some notation for the coefficients 
that arise from decomposing a vector into atoms.

\begin{defn} \label{def:adecomp} (atom decomposition)

Suppose $\gG$ is elementary and $J \in \mc{L}^-_\gG(\Phi)$.
For $\psi \in \Phi_B$ with $|B|=r$ we define 
integers $J^t_\psi$ for all nonzero $t \in T_B$ by
$J^{\psi\Ss} = \sum_{0 \ne t \in T_B} J^t_\psi \gG[\psi]^t$.
Any choice of orbit representatives $\psi^O \in \Phi_{B^O}$ 
for each orbit $O \in \Phi_r/\Ss$ defines an atom decomposition
$J = \sum_{O \in \Phi_r/\Ss} \sum_{0 \ne t \in T_{B^O}} J^t_{\psi^O} \gG[\psi^O]^t$.
\end{defn}

We need one final definition before
stating our main theorem in the next section;
note that the coefficients $G^t_\psi$ are 
as in the previous definition.

\begin{defn} \label{def:reg} (regularity)
Suppose $\gG \in (\mb{Z}^D)^{\mc{A}_r}$
and $G \in (\mb{Z}^D)^{\Phi_r}$. Let 
\[ \mc{A}(\Phi,G) 
= \{ \phi \in \mc{A}(\Phi): \gG(\phi) \le_\gG G \}.\]
We say $G$ is $(\gG,c,\oO)$-regular (in $\Phi$) if there is 
$y \in [\oO n^{r-q},\oO^{-1} n^{r-q}]^{\mc{A}(\Phi,G)}$ such that for all 
$B \in [q]_r$, $\psi \in \Phi_B$, $0 \ne t \in T_B$ we have 
\[\pl^t y_\psi := \sum_{\phi: t_\phi(\psi)=t} y_\phi = (1 \pm c)G^t_\psi.\]
\end{defn}

Note that if $G$ and $y$ are as in the previous definition
and $\psi \in O \in \Phi_r/\Ss$ then
\[ (\pl^\gG y)^O = \sum_\phi y_\phi \gG(\phi)^O 
= \sum_{0 \ne t \in T_B} \sum_{\phi: t_\phi(\psi)=t} y_\phi \gG[\psi]^t
= \sum_{0 \ne t \in T_B} (1 \pm c)G^t_\psi \gG[\psi]^t = (1 \pm c)G^O.\]
We note for future reference that this implies an upper bound
on the use (see Definition \ref{def:use+bdd} below)
of any orbit $O \in \Phi_r/\Ss$, namely $U(G)_O < 2|\mc{A}|\oO^{-1}$
(if $c<1/2$). Also, summing over $O$ we obtain
$\pl^\gG y = (1 \pm c)G$, i.e.\
$\sum_\phi y_\phi \gG(\phi)_{\psi,d} = (1 \pm c)G_{\psi,d}$ 
for all $\psi \in \Phi_r$, $d \in [D]$.

\section{Main theorem} 

Now we can state our main theorem.
We will give the proof in this section,
assuming Lemmas \ref{hole/bddint} and \ref{bddint},
which will be proved in sections \ref{sec:cea} 
and \ref{sec:bddint}. The parts of the proof
given in this section are those that are 
somewhat similar to the proof in \cite{Kexist},
so we will be quite concise in places
where they are similar, and give more
details at points of significant difference.
To apply Theorem \ref{main}
we also need a concrete description of 
the decomposition lattice $\bgen{\gG(\Phi)}$;
this will be given in section \ref{sec:int}.

\begin{theo} \label{main}
For any $q \ge r$ and $D$ there are $\oO_0$ and $n_0$ such that 
the following holds for $n > n_0$, $h=2^{50q^3}$, $\dD = 2^{-10^3 q^5}$,
$n^{-\dD}<\oO<\oO_0$ and $c \le \oO^{h^{20}}$.
Let $\mc{A}$ be a $\Ss^\le$-family with $\Ss \le S_q$.
Suppose $\gG \in (\mb{Z}^D)^{\mc{A}_r}$ is elementary.
Let $\Phi$ be a $\Ss$-adapted $[q]$-complex on $[n]$. 
Let $G \in \bgen{\gG(\Phi)}$ be $(\gG,c,\oO)$-regular 
in $\Phi$ such that $(\Phi,\gG[G]^A)$ 
is $(\oO,h)$-extendable for each $A \in \mc{A}$.
Then $G$ has a $\gG(\Phi)$-decomposition.
\end{theo}

Throughout this section we let
$\Ss$, $\mc{A}$, $\gG$, $\Phi$ and $G$ 
be as in the statement of Theorem \ref{main}.
We note that the assumption that $\gG$ is elementary
bounds $|\mc{A}|$ as a function of $q$ and $D$,
say $|\mc{A}| < (Dq)^{q^q}$.
For convenient reference, we list here
several parameters used throughout the paper.
\begin{gather*}
Q = \tbinom{q}{r}, \quad z=h=2^{50q^3}, \quad \dD = 2^{-10^3 q^5},
\quad n^{-\dD}<\oO<\oO_0(q,D), \quad \oO_q:=\oO^{(9q)^{q+5}}, \\
\quad p \text{ is a prime with } 2^{8q}<p<2^{9q}, \quad
a \in \mb{N} \text{ with } p^{a-2} < n \le p^{a-1}, \\
\gG = np^{-a}, \quad
\rho = \oO z^{-Q} |\mc{A}|^{-1} (q)_r^{-Q} \gG^{q-r}, 
\text{ where } (q)_r=q!/(q-r)!,\\
c=\oO^{h^{20}}, \quad c_1 = (2Qc)^{1/2Q}, \quad
c_{i+1} = \oO^{-h^3} c_i \ \text{ for } \ i \in [4].
\end{gather*}

\subsection{Probabilistic methods} 

\def\olddom{2.4}
\def\oldlip3{2.11}

We briefly recall two concentration 
inequalities 
(see \cite[Lemmas \olddom \ and \oldlip3]{Kexist}).

\begin{defn} 
Suppose $Y$ is a random variable and 
$\mc{F} = (\mc{F}_0,\dots,\mc{F}_n)$ is a filtration.
We say that $Y$ is \emph{$(C,\mu)$-dominated (wrt $\mc{F}$)} 
if we can write $Y = \sum_{i=1}^n Y_i$, 
where $Y_i$ is $\mc{F}_i$-measurable, $|Y_i| \le C$
and $\mb{E}[|Y_i| \mid \mc{F}_{i-1}] < \mu_i$ for $i \in [n]$, 
where $\sum_{i=1}^n \mu_i < \mu$.
\end{defn}

\begin{lemma} \label{dom}
If $Y$ is $(C,\mu)$-dominated then 
$\mb{P}(|Y|>(1+c)\mu) < 2e^{-\mu c^2/2(1+2c)C}$.
\end{lemma}

\begin{defn}
Let $a = (a_1,\dots,a_n)$ and $a' = (a'_1,\dots,a'_n)$, 
where $a_i \in \mb{N}$ and $a'_i \in [a_i]$ for $i \in [n]$,
and $\Pi(a,a')$ be the set of $\pi=(\pi_1,\dots,\pi_n)$
where $\pi_i:[a'_i] \to [a_i]$ is injective. 
Suppose $f:\Pi(a,a') \to \mb{R}$ and $b = (b_1,\dots,b_n)$ 
with $b_i \ge 0$ for $i \in [n]$.
We say that $f$ is \emph{$b$-Lipschitz} if for any $i \in [n]$
and $\pi,\pi' \in \Pi(a,a')$ such that $\pi_j = \pi'_j$ for $j \ne i$
and $\pi_i = \tau \circ \pi'_i$ for some transposition $\tau \in S_{a_i}$
we have $|f(s)-f(s')| \le b_i$. 
We also say that $f$ is \emph{$B$-varying} 
where $B=\sum_{i=1}^n a'_i b_i^2$.
\end{defn}

\begin{lemma} \label{lip3} 
Suppose $f:\Pi(a,a') \to \mb{R}$ is $B$-varying 
and $X=f(\pi)$, where $\pi=(\pi_i) \in \Pi(a,a')$ 
is random with $\{\pi_i: i \in [n]\}$ independent
and $\pi_i$ uniform whenever $a'_i>1$.
Then $\mb{P}(|X-\mb{E}X|>t) \le 2e^{-t^2/2B}$.
\end{lemma}

The following lemma will be used to pass from
fractional matchings to almost perfect matchings.
The statement and proof are similar\footnote{
The constraint set $P$ acts on the edges of $H$ in \cite{KaLP},
whereas here it is more convenient to use vertices.
It is assumed to be of constant size in \cite{KaLP}, 
which allows for a simple second moment argument,
but we need to allow $P$ to grow polynomially in $\aA^{-1}$,
which can be achieved by proving exponential tails on
the failure probabilities (which follows from the same proof
by applying standard concentration inequalities).
Also, the error term in the conclusion is not
explicitly given in \cite{KaLP}, whereas we state a polynomial
dependence on $\aA$ (the proof gives $\aA^{1/2k}$), 
which is needed if one desires counting versions of our results, 
as in \cite{Kcount}. A final comment is that 
it is not essential to consider edge weights, as we anyway reduce 
to the case that all $w_e$ are equal, but it is convenient
to leave the weights in the statement, and this also
facilitates comparison with the statement in \cite{KaLP}.} 
to those given by Kahn \cite{KaLP}, so we omit the details.
Call a hypergraph $H$ a $k^\le$-graph
if all edges have size at most $k$.

\begin{lemma} \label{kahn+}
Suppose $H$ is a $k^\le$-graph
and $w$ is a fractional matching in $H$
with $\sum_{\{x,y\} \sub e \in H} w_e 
< \aA < \aA_0(k)$ sufficiently small
for all $\{x,y\} \sub V(H)$. 
Let $P \sub \mb{R}^{V(H)}$ with $|P| < \aA^{-k}$ 
and $\max_v p_v < (\log \aA)^{-2} \sum_v p_v$ 
for all $p \in P$. 
Then there is a matching $M$ of $H$ such that
\[\sum_{v \in \bigcup M} p_v = (1 \pm \aA^{1/2k}) 
\sum_{e \in H} w_e \sum_{v \in e} p_v 
\quad \text{ for all } p \in P.\]
\end{lemma}

We will apply Lemma \ref{kahn+} to a hypergraph
whose vertices can be identified with $\gG$-atoms,
we have $\aA=O(n^{-1})$, and elements of $P$
indicate atoms that `use' a given 
ordered $(r-1)$-tuple from $[n]$;
the conclusion will be that there is 
a matching with `bounded leave'
(see Definition \ref{def:use+bdd}
and Lemma \ref{nibble} below).

\subsection{Template} 

Recalling the proof strategy discussed in the introduction,
we start by describing the template.
This will be determined by some $M^* \sub \mc{A}(\Phi)$
such that $G^* := \sum_{\phi \in M^*} \phi(Q) \sub \Phi^\circ_r$,
i.e.\ $G^*$ is an $r$-graph (with no multiple edges)
contained in $\Phi^\circ_r$ and $\{ \phi(Q): \phi \in M^* \}$
is a $K^r_q$-decomposition of $G^*$. 
Thus for each $e \in \Phi^\circ_r$ there is at most 
one orbit $O \in \Phi_r/\Ss$ with $Im(O)=e$
and $O \sub \phi\Ss$ for some $\phi \in M^*$, 
and given such $O$ with representative $\psi^O \in \Phi_B$ 
the use of $O$ by $\phi$ has
a unique type $t_\phi(\psi^O)=t \in T_B$
(which may be the zero type).

As in \cite{Kexist}, we fix $M \in \mb{F}_p^{q \times r}$ 
as a $q \times r$ matrix over $\mb{F}_p$ that is generic,
in that every square submatrix of $M$ is nonsingular.

As $G$ is $(\gG,c,\oO)$-regular, there is 
$y \in [\oO n^{r-q},\oO^{-1 }n^{r-q}]^{\mc{A}(\Phi,G)}$ with 
$\pl^t y_\psi = (1 \pm c)G^t_\psi$ 
for all $B \in [q]_r$, $\psi \in \Phi_B$, $0 \ne t \in T_B$.
We activate each $\phi \in \mc{A}(\Phi)$ independently
with probability $y_\phi \oO n^{q-r}$.

Let $f = (f_j: j \in [z])$, with $z=h=2^{50q^3}$,
where we choose independent uniformly random
injections $f_j: [n] \to \mb{F}_{p^a}$.
Given $f$, for each $e \in \Phi^\circ_r$ we let
\[\mc{T}_e = \{j \in [z]: \dim(f_j(e))=r\}.\]
We abort if any $|\mc{T}_e| \le z-2r$,
which occurs with probability $O(n^{-r})$.
We assume without further comment 
that the template does not abort.

We choose $T_e \in [z]$ for all $e \in \Phi^\circ_r$ 
independently and uniformly at random. 
We say $\phi \in \mc{A}(\Phi)$ is compatible with $j$
if $T_e=j \in \mc{T}_e$ for all $e \in \phi(Q)$
and for some $y \in \mb{F}_{p^a}^r$ we have
$f_j(\phi(i)) = (My)_i$ for all $i \in [q]$.

Let $\pi = (\pi_e: e \in \Phi^\circ_r)$ 
where we choose independent uniformly random
injections $\pi_e:e \to [q]$.
We say $\phi \in \Phi$ is compatible with $\pi$
if $\pi_e\phi(i)=i$ whenever
$\phi(i) \in e \in \phi(Q)$
(for brevity we write this as $\pi_e\phi=id$).

We choose independent uniformly random $A_\phi \in \mc{A}$
for each injection $\phi: [q] \to [n]$.

\begin{defn} \label{def:template}
Let $M^*_j$ be the set of all activated
$\phi \in A_\phi(\Phi)$ compatible with $j$ and $\pi$
such that $\gG(\phi) \le_\gG G$. 
The template is $M^* = \cup_{j \in [z]} M^*_j$.

The underlying $r$-graph of the template
is $G^* = \cup_{j \in [z]} G^*_j$, where
each $G^*_j = \cup_{\phi \in M^*_j} \phi(Q)$.
\end{defn}

\def\oldtdec{3.3}

Note that $G^*$ is the edge-disjoint union of the $G^*_j$,
and each $\{ \phi(Q): \phi \in M^*_j \}$
is a $K^r_q$-decomposition of $G^*_j$
(see \cite[Lemma \oldtdec]{Kexist}).

We introduce the following further notation that
will be used in the analysis of the template.

\begin{defn} \label{def:further}
For $e \in G^*$ let $\phi^e \in M^*$ be such 
that $e \in \phi^e(Q)$. We write $M^*(e)=\phi^e(Q)$.

For $J \sub G^*$ let $M^*(J) = \sum_{e \in J} M^*(e) \in \mb{N}^{G^*}$.

For $e \in G^*$ we write $\gG(e) = \gG(\phi^e)^{O_e}$
where $O_e = \pi_e^{-1} \Ss$. We call $\gG(e)$ an $M^*$-atom.
\end{defn}

Note that $\gG(e)$ may be zero in Definition \ref{def:further}.
If $\gG(e) \ne 0$ then $\gG(e)$ is also a $\gG$-atom. 
Furthermore, for any $\gG$-atom $\gG(\psi) \le_\gG \pl^\gG M^*$ 
we have $\gG(\psi)=\gG(e)$ where $Im(\psi)=e$.

If $0 \le_\gG J \le_\gG \pl^\gG M^*$
with $J = \sum Z$ for some set $Z$ of $\gG$-atoms
we write $J^\circ = \{ e \in G^*: \gG(e) \in Z \}$.
For example, $G^* = (\pl^\gG M^*)^\circ$.

\subsection{Extensions}

Next we give estimates on the probability
that certain $\gG$-atoms appear in the template
and deduce that the template is whp extendable.
Our estimates are conditional on the following local events 
(defined similarly to \cite{Kexist}).

\begin{defn} (local events)
Suppose $e \in \Phi^\circ_r$.
We reveal $T_e=j$ and $f_j\!\mid_e=\aA$.
If $\dim(\aA)<r$ then $\mc{E}^e$ is the event
that $T_e=j$ and $f_j\!\mid_e=\aA$,
which witnesses $e \notin G^*$.

Now suppose $\dim(\aA)=r$, reveal $\pi_e$,
and let $y \in \mb{F}_{p^a}^r$ 
with $f_j(x) = (My)_i$ for all $x \in e$, $\pi_e(x)=i$.
We reveal $f_j^{-1}((My)_i)$ for all $i \in [q] \sm \pi_e(e)$,
and let $\phi:[q] \to [n]$ be such that $f_j\phi=My$,
and reveal $A_\phi$. If we do not have $\gG(\phi) \le_\gG G$
with $\phi \in A_\phi(\Phi)$ then $\mc{E}^e$ is the event 
that $T_e=j$ and $f_j\phi=My$, which witnesses $e \notin G^*$.

Finally, if $\gG(\phi) \le_\gG G$ we reveal whether $\phi$ is activated,
and reveal $(T_{e'},\pi_{e'})$ for all $e' \in \phi(Q) \sm \{e\}$.
Then $\mc{E}^e$ is defined by all the above information,
which determines whether $e \in G^*$: given 
$T_e=j$, $f_j\phi=My$, $\phi \in A_\phi(\Phi)$, $\gG(\phi) \le_\gG G$ 
we have $e \in G^*$ iff $\phi$ is activated and 
$T_{e'}=j$ and $\pi_{e'}\phi=id$ for all $e' \in \phi(Q)$.

We say that a vertex $x$ is touched by $\mc{E}^e$
if $f_j(x)$ is revealed by $\mc{E}^e$.

We say that an edge $e'$ is touched by $\mc{E}^e$
if $T_{e'}$ is revealed by $\mc{E}^e$.
\end{defn}

\def\olde|E{3.6}

The following lemma is analogous
to \cite[Lemma \olde|E]{Kexist}. Let 
\[ \rho := \oO z^{-Q} |\mc{A}|^{-1} (q)_r^{-Q} \gG^{q-r}. \]

\begin{lemma} \label{e|E}
Let $S \sub \Phi^\circ_r$ with $|S|<h=z$
and $\mc{E} = \cap_{f \in S} \mc{E}^f$.
Let $\psi \in \Phi_B$ and 
$t \in T_B$ with $t \ne 0$ and
$\gG[\psi]^t \le_\gG G$ (i.e.\ $G^t_\psi > 0$).
Suppose $e:=Im(\psi)$ is not touched by $\mc{E}$
and $j \in [z] \sm \{T_f: f \in S\}$.
Then $\mb{P}( \gG[\psi]^t \le_\gG \pl^\gG M^*_j \mid \mc{E})
= (1 \pm 1.1c) \rho G^t_\psi$.
\end{lemma}

\nib{Proof.}
We fix any $\phi \in A(\Phi)$ with
$\gG(\phi) \le_\gG G$ and $t_\phi(\psi)=t$,
and estimate the probability that $\phi \in M^*_j$.
We have $\mb{P}(A_\phi=A) = |\mc{A}|^{-1}$.
We activate $\phi$ with probability $y_\phi \oO n^{q-r}$.
We can assume every $e' \in \phi(Q)$ is not touched
by $\mc{E}$, as this excludes $O(n^{q-r-1})$ choices.
Then all $T_{e'}=j$ with probability $z^{-Q}$.
With probability $(q)_r^{-Q}$ all $\pi_{e'}\phi=id$.
We condition on $f_j\!\mid_e$ such that $\dim(f_j(e))=r$;
this occurs with probability $1-O(n^{-1})$.
There is a unique $y \in \mb{F}_{p^a}^r$ 
such that $(My)_i = f_j(x)$ for all $x \in e$, $i=\pi_e(x)$.
With probability $(1+O(n^{-1}))(p^{-a})^{q-r}$
we have $f_j(\phi(i))=(My)_i$ 
for all $i \in [q] \sm \pi_e(e)$.
Therefore $\mb{P}(\phi \in M^*_j \mid \mc{E})
= (1+O(n^{-1})) |\mc{A}|^{-1} y_\phi \oO n^{q-r} z^{-Q} 
 (q)_r^{-Q} (p^{-a})^{q-r}$.
The lemma follows by summing over $\phi$, 
using $\pl^t y_\psi = (1 \pm c)G^t_\psi$. \qed

\begin{rem} \label{rem:e|E}
The same proof shows 
\begin{enumerate}
\item $\mb{P}( \gG[\psi]^t \le_\gG \pl^\gG M^*_j 
\mid \mc{E} \cap \{T_e=j\})
= (1 \pm 1.1c) z\rho G^t_\psi$
for any $j \in [z]\sm \{T_f: f \in S\}$, 
\item $\mb{P}( \{\gG[\psi]^t \le_\gG \pl^\gG M^*_j\}
\cap \{\pi_e=\pi\} \cap \{A_{\phi^e}=A\} \mid \mc{E} ) 
> \oO^2 \rho G^t_\psi$, for any $A \in \mc{A}$
and injection $\pi: e \to [q]$
such that $\psi' = \pi^{-1} \in A(\Phi)^\le$
satisfies $\gG[\psi]^t=\gG(\psi')$.
\end{enumerate}
Note that (ii) is weaker than the corresponding
bound in \cite{Kexist} as we cannot permute $\phi$ 
(this may change $\gG(\phi)$).
Instead, we use extendability to see that
there are at least $\oO n^{q-r}$ choices
of $\phi$ with $\gG(\phi) \le_\gG G$ containing $\psi'$,
and each has $y_\phi > \oO n^{r-q}$.
\end{rem}

\def\oldext*{3.8}

The following lemma is analogous
to \cite[Lemma \oldext*]{Kexist}.

\begin{lemma} \label{ext*}
Suppose $E=(\phi,F,H)$ is a $\Phi$-extension with $|H| \le h/3$.
Let $A \in \mc{A}$ and $H' \sub H_r \sm H[F]$.
Then whp $X_{E,H'}(\Phi,\gG[\pl^\gG M^*]^A) 
> \oO n^{v_E} (z\rho/2)^{|H'|}$.
\end{lemma}

\nib{Proof.} 
As $(\Phi,\gG[G]^A)$ is $(\oO,h)$-extendable,
there are at least $\oO n^{v_E}$
choices of $\phi^+ \in X_{E,H'}(\Phi,\gG[G]^A)$,
i.e.\ $\phi^+ \in X_E(\Phi)$ with 
$\phi^+ \psi' \in \gG[G]^A_B$ 
for all $\psi' \in H'_B$.
We fix any such $\phi^+$ and estimate 
$\mb{P}( \phi^+ \in X_{E,H'}(\Phi,\gG[\pl^\gG M^*]^A) )$
by repeated application of Lemma \ref{e|E}.
For any $\psi=\phi^+\psi'$ with $\psi' \in H'_B$,
we condition on the intersection $\mc{E}$ of
all previously considered local events, and estimate
$p^j_\psi = \mb{P}( \gG[\psi]^t \le_\gG \pl^\gG M^*_j \mid \mc{E})$,
where $t \in T_B$ contains $id_B \in A$ and
we can assume $Im(\psi)$ is not touched by $\mc{E}$.
If $t=0$ then $p^j_\psi=1$; otherwise, there are 
at least $2z/3$ choices of $j \in [z]$ not used by 
any previous edge such that Lemma \ref{e|E} applies to give
$p^j_\psi = (1 \pm 1.1c) \rho G^t_\psi$. Multiplying all conditional 
probabilities and summing over $\phi^+$ gives 
$\mb{E}X_{E,H'}(\Phi,\gG[\pl^\gG M^*]^A) 
> \oO n^{v_E} (0.6 z\rho)^{|H'|}$.
The proof of concentration is similar
to that in \cite[Lemma \oldext*]{Kexist}, noting that the effect
of changing any $A_\phi$ has a similar effect
to that of changing whether $\phi$ is activated. \qed

\subsection{Approximate decomposition}

Similarly to \cite[section 4]{Kexist}
we will now complete the template
to an approximate decomposition,
namely $M' \sub \mc{A}(\Phi)$
such that $\pl^\gG M'$ is almost equal to $G$,
except that some (suitably bounded) set of $M^*$-atoms
are each covered one time too many.
First we introduce some notation and terminology
that will be used throughout the rest of the paper.

\begin{defn} \label{def:use+bdd}
For $J \in (\mb{Z}^D)^{\Phi_r}$ and $\psi \in \Phi_r$, 
we define the {\em use} $U(J)_\psi$ of $\psi$ by $J$ 
as the minimum possible value of $\sum_{w \in W} |x_w|$
where $W$ is the set of $\gG$-atoms at $O=\psi\Ss$
and $x \in \mb{Z}^W$ with $J^O = \sum x_w w$.
If there is no such $x$
then $U(J)_\psi$ is undefined.
For $\psi' \in \Phi$ we let $U(J)_{\psi'} 
= \sum \{ U(J)_\psi: \psi' \sub \psi \in \Phi_r \}$.
We note that use is a property of orbits,
so $U(J)_{\psi\Ss} = U(J)_\psi$ is well-defined.
We say $J$ is $\tT$-bounded
if $U(J)_\psi < \tT |V(\Phi)|$ 
whenever $\psi \in \Phi_{r-1}$. 
\end{defn}

Note that as $\gG$ is elementary the use of `minimum'
in Definition \ref{def:use+bdd} is redundant,
as $J^O$ has a unique atom decomposition;
however, we will also need this definition
for other $\gG$ that are not necessarily elementary.
We will also sometimes use the following definition that
ignores the edge labellings and atom structure
(so is analogous to that used in \cite{Kexist}).

\begin{defn} \label{def:use+bdd(unordered)}
Suppose $J \in (\mb{Z}^D)^{\Phi^\circ_r}$.
We let $U(J)_e = \sum_{d \in [D]} |(J_e)_d|$
for $e \in \Phi^\circ_r$ and
$U(J)_f = \sum \{ U(J)_e: f \sub e \in \Phi^\circ_r\}$
for $f \in \Phi^\circ$.
We say $J$ is $\tT$-bounded if
$U(J)_f < \tT |V(\Phi)|$ 
for all $f \in \Phi^\circ_{r-1}$.
\end{defn}

\def\oldnibble{4.1}

Now we find a $\gG(\Phi)$-decomposition
of almost all of $G - \pl^\gG M^*$,
such that the leave is bounded 
in the sense of Definition \ref{def:use+bdd}.
The following lemma is analogous
to \cite[Lemma \oldnibble]{Kexist}. 

\begin{lemma} \label{nibble}
There is $M^n \sub \mc{A}(\Phi)$ with $c_1$-bounded  
{\em leave} $L := G - \pl^\gG M^* - \pl^\gG M^n \ge_\gG 0$. 
\end{lemma}

\nib{Proof.}
We define $\Phi' \sub \mc{A}(\Phi)$ randomly as follows.
Consider any $\phi \in \mc{A}(\Phi)$ such that 
$\gG(\phi) \le_\gG G$ and reveal the local events 
$\mc{E}^e$ for each $e \in Q':=\phi(Q)$.
If $\phi$ is not activated or
$T_e=T_{e'}$ for any $e \ne e'$ in $Q'$
then we do not include $\phi$ in $\Phi'$.
For each $e \in Q'$ we fix $\psi_e$ with image $e$
such that $\gG(\psi_e) \le_\gG \gG(\phi)$,
and let $t_e = t_\phi(\psi_e)$,
so $\gG(\psi_e) = \gG[\psi_e]^{t_e}$.
For $v \in \{0,1\}^{Q'}$ let $\mc{E}^\phi_v$ be the event that 
all $v_e = 1_{\gG(\psi_e) \le_\gG \pl^\gG M^*} $. If $\phi$ is activated,
all $T_e$ for $e \in Q'$ are distinct and $\mc{E}^\phi_v$ holds 
then we include $\phi$ in $\Phi'$ independently with probability 
$p^v_\phi = \prod_{e \in Q'} (1-v_eq_e)$,
where $q_e = 1/G^{t_e}_{\psi_e}$ if $t_e \ne 0$
or $q_e = z\rho$ if $t_e=0$.

Now we fix any $\psi$ and
$t \ne 0$ with $G^t_\psi>0$ and estimate 
the number $X$ of $\phi \in \Phi'$ with $t_\phi(\psi)=t$.
We consider any activated $\phi$ with $t_\phi(\psi)=t$,
let $e'=Im(\psi)$, $Q'=\phi(Q)$ and condition on 
the local event $\mc{E}^{e'}$ and any event 
$\mc{C} = \cap_{e \in Q'} \{T_e=j_e\}$ 
such that all $j_e$ are distinct
(the latter occurs with probability $(z)_Q z^{-Q}$).

\def\oldnibble{4.1}

For any $v  \in \{0,1\}^{Q'}$ with $v_{e'}=(\pl^\gG M^*)^t_\psi$, 
by repeated application of Remark \ref{rem:e|E}.i we have
$\mb{P}(\mc{E}^\phi_v \mid \mc{E}^{e'} \cap \mc{C} )
= ( 1 \pm Qc) \prod_{e \in Q' \sm \{e'\}} p^{v_e}_e$, 
where if $t_e=0$ we let $p^1_e=1$ and $p^0_e=0$,
and otherwise $p^1_e = z\rho G^{t_e}_{\psi_e}$
and $p^0_e=1-z\rho G^{t_e}_{\psi_e}$.
Then $p^0_e + p^1_e (1-q_e) = 1-z\rho$,
so as in the proof of \cite[Lemma \oldnibble]{Kexist},
$\mb{P}[ \phi \in \Phi' \mid \mc{E}^{e'} \cap \mc{C} ]
= ( 1 \pm Qc) (1-(\pl^\gG M^*)^t_\psi/G^t_\psi) (1-z\rho)^{Q-1}$.

We activate each $\phi$ with probability $y_\phi \oO n^{q-r}$,
so as $\pl^t y_\psi = (1 \pm c)G^t_\psi$ we deduce
$\mb{E}X = (z)_Q z^{-Q} \cdot
( 1 \pm Qc) (1-(\pl^\gG M^*)^t_\psi/G^t_\psi) (1-z\rho)^{Q-1} 
\cdot (1 \pm c) \oO n^{q-r} G^t_\psi 
= (1 \pm 1.1Qc) d' n^{q-r} (G-\pl^\gG M^*)^t_\psi$,
where $d' = (z)_Q z^{-Q} \oO (1-z\rho)^{Q-1}$.
As in the proof of \cite[Lemma \oldnibble]{Kexist},
by Lemma \ref{lip3} whp
$X = (1 \pm 1.2Qc) d' n^{q-r} (G-\pl^\gG M^*)^t_\psi$.

Finally, we consider the following hypergraph $H$,
where $V(H)$ is the disjoint union of sets $V^t_\psi$
of size $(G-\pl^\gG M^*)^t_\psi$ corresponding to the 
$\gG$-atoms of $G-\pl^\gG M^*$ counted with multiplicity. 
For each $\phi \in \Phi'$ we let $\mc{V}_\phi$ be the set 
of all $V^t_\psi$ with $\gG[\psi]^t \le_\gG \gG(\phi)$,
and include as an edge a uniformly random set $e^\phi$
with one vertex in each $V \in \mc{V}_\phi$.
Then whp every $v \in V(H)$ has degree
$(1 \pm 2Qc) d' n^{q-r}$. Also any $\{u,v\} \sub V(H)$
is contained in $O(n^{q-r-1})$ edges.
Let $P \sub \mb{R}^{V(H)}$ where for each
$\psi' \in \Phi_{r-1}$ and we include $p^{\psi'}$
where $p^{\psi'}_v$ is $1$ if $v$ is in some $V^t_\psi$
with $\psi'\Ss \sub \psi\Ss$, otherwise $0$. 
Then $\sum_v p^{\psi'}_v$ counts the number of
$\gG$-atoms of $G-\pl^\gG M^*$ on orbits
containing $\psi'\Ss$. 

By Lemma \ref{kahn+}, applied with uniform weights 
$w_e=((1 + 2Qc) d' n^{q-r})^{-1}$,
there is a matching $M^n$ in $H$ with
$\sum_{v \in \bigcup M^n} p_v = (1 \pm (1.1Qc)^{1/2Q})
\sum_{e \in H} w_e \sum_{v \in e} p_v
= (1 \pm (1.2Qc)^{1/2Q}) \sum_{v \in V(H)} p_v$ 
for all $p \in P$. We can also view $M^n$ 
as a subset of $\mc{A}(\Phi)$. Then
$L := G - \pl^\gG M^* - \pl^\gG M^n$
contains at most $(1.2Qc)^{1/2Q}$ proportion
of the $\gG$-atoms of $G-\pl^\gG M^*$ on orbits
containing $\psi'\Ss$, for any $\psi' \in \Phi_{r-1}$.
Recalling that $U(G)_O < 2|\mc{A}|\oO^{-1}$
and $|\mc{A}| < D^{q^{2q}}$ we see that
$L$ is $c_1$-bounded. \qed
 
\medskip

\def\oldcover{4.2}

To complete the approximate decomposition,
we choose a partial $\gG(\Phi)$-decomposition
that exactly matches the leave on $\Phi^\circ_r \sm G^*$,
but has some `spill' $S$ in $G^*$ that we will
need to correct for later. 
The following lemma is analogous
to \cite[Lemma \oldcover]{Kexist}. 

\begin{lemma} \label{cover}
Suppose $0 \le_\gG L \le_\gG G -  \pl^\gG M^*$ is $c_1$-bounded.
Then there is $M^c \sub \mc{A}(\Phi)$ such that 
$\gG(\phi) \le_\gG L + \pl^\gG M^*$ for all $\phi \in M^c$
and $\pl^\gG M^c_\psi = L_\psi$ for all 
$\psi \in \Phi_r$ with $Im(\psi) \notin G^*$,
with {\em spill} $S := G^* \cap \sum_{\phi \in M^c} \phi(Q)$, 
such that $M^*(S)$ is a set and $c_2$-bounded.
\end{lemma}

\nib{Proof.}
We order the $\gG$-atoms of $L$ as $(\gG(\psi_i): i \in [n_L])$,
where each $\psi_i \in A^i(\Phi)^\le_{B^i} = \Phi_{B^i}$.
For each $i$ we consider the $\Phi$-extension 
$E_i=(\ova{q},B^i,\psi_i)$ and let $H' = \ova{q}_r \sm \{B^i\}$.
We apply a random greedy algorithm to select 
$\phi_i \in X_{E_i,H'_i}(\Phi,\gG[\pl^\gG M^*]^{A^i})$,
where we write $S_i = G^* \cap \cup_{i'<i} \phi_{i'}(Q)$
and choose $\phi_i$ uniformly at random
such that $M^*(\phi_i(Q))$ is a set disjoint from $M^*(S_i)$.
For each $i$ we add $\phi_i \in A^i(\Phi)$ to $M^c$.
The remainder of the proof is very similar
to that of \cite[Lemma \oldcover]{Kexist}.
We show that whp $M^c$ has the stated properties.
At any step $i$ before $M^*(S)$ fails to be $c_2$-bounded
at most half of the choices of $\phi_i$ are forbidden, as
$X_{E_i,H'_i}(\Phi,\gG[\pl^\gG M^*]^{A^i}) > \oO (z\rho/2)^Q n^{q-r}$
by Lemma \ref{ext*}. Then for all $e \in G^*$ we estimate
$r_e := \sum_{i \in [n_L]} \mb{P}'(e \in M^*(\phi_i(Q)))
< 2 (2q)^{2q} \oO^{-1} (z\rho/2)^{-Q} c_1$,
so by Lemma \ref{dom} whp $M^*(S)$ is $c_2$-bounded. \qed

\subsection{Proof modulo lemmas}

In this subsection we give the proof of Theorem \ref{main},
assuming two lemmas that will be proved later.
First we need to define use and boundedness
for vectors indexed by $\mc{A}(\Phi)$.

\begin{defn} \label{bddPsi}

For $\Psi \in \mb{Z}^{\mc{A}(\Phi)}$ and $\psi \in \Phi$
the use of $\psi$ by $\Psi$ is 
$U(\Psi)_\psi = \sum \{ |\Psi_\phi| : \psi\Ss \sub \phi\Ss \}$.
We also write $U(\Psi)_{\psi\Ss} = U(\Psi)_\psi$.
We say $\Psi$ is $\tT$-bounded if
$\sum \{ U(\Psi)_\psi: \psi' \sub \psi \in \Phi_r \} 
< \tT |V(\Phi)|$ whenever $\psi' \in \Phi_{r-1}$. 
\end{defn}

\def\oldbddint{5.1}

The following is a bounded integral decomposition lemma,
analogous to \cite[Lemma \oldbddint]{Kexist},
which will be proved in section \ref{sec:bddint}.
Recall  $\oO_q:=\oO^{(9q)^{q+5}}$.

\begin{lemma} \label{bddint}
Let $\mc{A}$ be a $\Ss^\le$-family 
with $\Ss \le S_q$ and $|\mc{A}| \le K$
and suppose $\gG \in (\mb{Z}^D)^{\mc{A}_r}$.
Let $\Phi$ be an $(\oO,h)$-extendable $\Ss$-adapted $[q]$-complex
on $[n]$, where $n^{-h^{-3q}}<\oO<\oO_0(q,D,K)$ and $n>n_0(q,D,K)$.
Suppose $J \in \bgen{\gG(\Phi)}$ is $\tT$-bounded, 
with $n^{-(5hq)^{-r}} < \tT < 1$.
Then there is some $\oO_q^{-2h} \tT$-bounded
$\Psi \in \mb{Z}^{\mc{A}(\Phi)}$ with $\pl^\gG \Psi = J$. 
\end{lemma}

\def\oldhole+{8.1}

The following lemma takes as input
a bounded integral decomposition 
as produced by Lemma \ref{bddint}
and produces a signed decomposition
analogous to that in \cite[Lemma \oldhole+]{Kexist}.
In the next subsection we will reduce
Lemma \ref{hole+/bddint} to Lemma \ref{hole/bddint}, 
which will be proved in section \ref{sec:cea}.
Let $c'_2 = \oO^{-h^2} c_2$.

\begin{lemma} \label{hole+/bddint}
Suppose $\Psi \in \mb{Z}^{\mc{A}(\Phi)}$ is $c'_2$-bounded 
with $0 \le_\gG \pl^\gG \Psi \le_\gG \pl^\gG M^*$.
Let $S = (\pl^\gG \Psi)^\circ$ and
suppose $M^*(S)$ is a set.
Then there is $M^o \sub M^*$ and $M^i \sub \mc{A}(\Phi)$
such that $\pl^\gG M^o = \pl^\gG M^i + \pl^\gG \Psi$.
\end{lemma}

The proof of Theorem \ref{main} is now
quite short given these lemmas.

\medskip

\nib{Proof of Theorem \ref{main}.}
Fix a template $M^*$
as in Definition \ref{def:template} that 
satisfies all of the whp statements in the paper.
Let $M^n$ be obtained from Lemma \ref{nibble}
and $M^c$ and $S$ from Lemma \ref{cover}.
Let $J = \pl^\gG (M^* + M^n + M^c) - G$
and note that $J \in \bgen{\gG(\Phi)}$,
$0 \le_\gG J \le_\gG \pl^\gG M^*$ and $J^\circ=S$.
Then $J$ is $c_2$-bounded, so by Lemma \ref{bddint}
there is some $c'_2$-bounded
$\Psi \in \mb{Z}^{\mc{A}(\Phi)}$ with $\pl^\gG \Psi = J$. 
Then we can apply Lemma \ref{hole+/bddint} 
to obtain $M^o \sub M^*$ and $M^i \sub \mc{A}(\Phi)$
such that $\pl^\gG M^o = \pl^\gG M^i + J$.
Now $M = M^n \cup M^c \cup (M^* \sm M^o) \cup M^i$
is a $\gG(\Phi)$-decomposition of $G$. \qed

\subsection{Absorption}

\def\oldab{6}

In this subsection, 
we establish the algebraic absorbing
properties of the template, and so reduce
Lemma \ref{hole+/bddint} to Lemma \ref{hole/bddint}.
Following \cite[section \oldab]{Kexist},
with appropriate modifications 
for the more general setting here,
we will define absorbers, cascades
and cascading cliques, then estimate 
the number of cascades for any cascading cliques.
We will always be concerned with absorbing maps
that are compatible with the template
(possibly with one `bad edge'),
as in the following definition.

\begin{defn} \label{def:M*compat}
Let $\phi \in A(\Phi)$ for some $A \in \mc{A}$.
We say that $\phi$ is $M^*$-compatible 
if $\phi$ is $\pi$-compatible and
$\phi^e \in A(\Phi)$ for all $e \in \phi(Q) \cap G^*$.
Also, for $e' \in \phi(Q)$ and 
$Q^* = \phi(Q) \cap G^* \sm \{e'\}$,
we say that $\phi$ is $M^*$-compatible bar $e'$
if $\pi_e \phi = id$ and $\phi^e \in A(\Phi)$
for all $e \in Q^*$.
\end{defn}

Note that if $\phi$ is $M^*$-compatible 
and $\phi(Q) \sub G^*$ then $\gG(e)=\gG(\psi)$ 
whenever $e=Im(\psi) \in G^*$ with $\psi \sub \phi$.
To define absorbers we require some notation:
let $Ker := \{ a \in \mb{F}_p^q: aM=0 \}$ and 
\[ v^w_a = M M_{[r]}^{-1} (e_{[r]} + a) w,
\quad \ v'{}^w_a = w+Maw \quad
\text{ for } a \in Ker^r
\text{ and } w \in  \mb{F}_{p^a}^q.\]

\begin{defn} \label{def:absorb} (absorbers)
Let $\phi \in A(\Phi)$ be $M^*$-compatible
with $\phi(Q) \sub G^*_j$ and $\dim(w)=q$,
where $w:=f_j\phi \in \mb{F}_{p^a}^q$.
Suppose $\phi^w: [q] \times Ker \to [n]$ such that
\begin{enumerate}
\item $f_j\phi^w((i,a)) = w_i + a \cdot w$ 
for each $i \in [q]$, $a \in Ker$,
\item if $\phi' \in [q](Ker)$ with
$f_j \phi^w \phi' = v^w_a$ for some $a \in Ker^r$
then $A_{\phi^w \phi'}=A$ and $\phi^w \phi' \in M^*_j$. 
\end{enumerate}

We say that $\phi$ is absorbable
and call $\phi^w$ the absorber for $\phi$.
We also refer to the subgraph\footnote{
Here we use `$\TT$' rather than
the natural `A' used in \cite{Kexist} 
to avoid clashes with other uses of `A'
in the paper.} 
$\TT^{\phi(Q)}=\TT^w=\phi^w(K^r_q(Ker))$ of $G^*$
as the absorber for $\phi(Q)$.
\end{defn}
 
\def\olddim{6.3}

\def\oldshuffle{6.4}

\def\oldexchange{6.5}

We denote the edges of $\TT^w$ by $e^w_a$
where $f_j(e^w_a) = (e_I+a)w$
for some $a \in Ker^I$, $I \in Q$.
As in \cite[Lemma \olddim]{Kexist},
each edge has full dimension 
under the relevant embedding:
we have $\dim(f_j(e^w_a))=r$.
As in \cite{Kexist}, we also view
$[q](Ker)$ as a subset of $\mb{F}_p^{q \times q}$,
and then we can write Definition \ref{def:absorb}.i
as $f_j \phi^w \phi' = w + \phi' w$.
We define
\[\phi^a:=M M_{[r]}^{-1} (e_{[r]} + a)-I, 
\quad \text{ and } \quad
\phi'{}^a = Ma,\]
so that
$f_j \phi^w \phi^a = w + \phi^a w = v^w_a$ and
$f_j \phi^w \phi'{}^a = w + \phi'{}^a w = v'{}^w_a$.
We write
\[ \Psi(\phi^w) = \{ \phi^w\phi^a: a \in Ker^r\}
\quad \text{ and } \quad
\Psi'(\phi^w) = \{ \phi^w\phi'{}^a: a \in Ker^r\}, \]
noting that $\Psi(\phi^w) \sub M^*_j$
and $\phi \in \Psi'(\phi^w)$.
By \cite[Lemma \oldshuffle]{Kexist},
$\Psi(\phi^w)$ and $\Psi'(\phi^w)$ 
both give $K^r_q$-decompositions of $\TT^{\phi(Q)}$.
Furthermore, $\pl^\gG \Psi(\phi^w) 
= \pl^\gG \Psi'(\phi^w) \le_\gG \pl^\gG M^*$.

Next we recall \cite[Lemma \oldexchange]{Kexist}.

\begin{lemma} \label{exchange}
There are $K^r_q$-decompositions  $\Ups$ and $\Ups'$ 
of $\OO=K^r_q(p)$ such that
\begin{enumerate}
\item $|V(f) \cap V(f')| \le r$
for all $f \in \Ups$ and $f' \in \Ups'$,
\item if $f \in \Ups$ and $\{f',f''\} \sub \Ups'$
with $|V(f) \cap V(f')|=|V(f) \cap V(f'')|=r$

then $(V(f') \sm V(f)) \cap (V(f'') \sm V(f))=\es$.
\end{enumerate}
\end{lemma}  

We identify $\Ups$ and $\Ups'$ with subsets of $[q](p)$,
i.e.\ the set of partite maps from $[q]$ to $[q] \times [p]$.
We identify $[q]$ with $\{ (i,1): i \in [q] \} \sub V(\OO)$
and with the corresponding map $id_{[q]}$;
by relabelling we can assume $[q] \in \Ups$.
For $U' \sub U \sub [n]$
we say that $U$ is $j$-generic for $U'$ 
if $\dim(f_j(U))=\dim(f_j(U'))+|U|-|U'|$.
Now we can define cascades.

\begin{defn} \label{def:cascade} (cascades)
Let $\phi \in A(\Phi)$ be $M^*$-compatible.
Suppose $\phi^c$ is an embedding of $K^r_q(p)$ in $G^*_j$
where $\phi^c\ id_{[q]} = \phi$ and $Im(\phi^c)$ 
is $j$-generic for $Im(\phi)$, such that each 
$\phi^c\phi'$ with $\phi' \in \Ups'$
has $A_{\phi^c\phi'}=A$ and is absorbable,
with absorber $\TT^{\phi^c\phi'(Q)}=\phi^{w_{\phi'}}(K^r_q(Ker))$, 
and $C_{\phi^c} = \sum \{\TT^{\phi^c\phi'(Q)}: \phi' \in \Ups'\}$
is a set (without multiple elements). 
We call $C_{\phi^c}$ a cascade for $\phi$.
\end{defn}

To flip a cascade $C_{\phi^c}$ we replace 
\[\Psi(C_{\phi^c}) := \bigcup \{ \Psi(w_{\phi'}): \phi' \in \Ups' \}
\qquad \text{ by } \]
\[ \Psi'(C_{\phi^c}) := \{ \phi^c\phi': \phi' \in \Ups \} \cup \bigcup 
\{  \Psi'(w_{\phi'}) \sm \{\phi^c\phi'\}: \phi' \in \Ups' \}.\]
This modifies $M^*$ so as to include $\phi$.
Next we define the class of cliques for
which we will show that there are many cascades.

\begin{defn} \label{def:cascading} (cascading cliques)
Let $\mc{Q}^* = \cup_{j \in [z]} \mc{Q}_j$,
where each $\mc{Q}_j$ is the set of all
$M^*$-compatible $\phi \in \mc{A}(\Phi)$
with $\dim(f_j\phi)=q$ and $\phi(Q) \sub G^*_j$
(we call $\phi$ cascading).
\end{defn}

\def\oldrole{6.10}

The following is \cite[Lemma \oldrole]{Kexist}.

\begin{lemma} \label{cascade:role}
Suppose $\phi \in \mc{Q}^*$, $Q'=\phi(Q)$ 
and $e \in G^*$ with $|e \sm V(M^*(Q'))|=r'$.
Let $\phi' \in \Ups'$, $I \in Q$, $a \in Ker^I$.
Then there are at most $p^{2q}n^{q(p-1)-r'}$ 
cascades $C_{\phi^c}$ for $\phi$ such that the absorber
$\TT^{\phi^c\phi'(Q)}=\phi^{w_{\phi'}}(K^r_q(Ker))$
for $\phi^c\phi'$ satisfies $e^{w_{\phi'}}_a=e$.
\end{lemma}

\def\oldcascade{6.11}

The following is analogous to \cite[Lemma \oldcascade]{Kexist}.

\begin{lemma} \label{cascade}
whp for any cascading $\phi \in \mc{Q}^*$ there are
at least $\oO^{p^{q^2}} n^{q(p-1)}$ cascades for $\phi$.
\end{lemma}

\nib{Proof.}
We follow the proof in \cite{Kexist},
indicating the necessary modifications.
Suppose $\phi \in A(\Phi)$, let $Q'=\phi(Q)$, 
and condition on local events 
$\mc{E} = \cap_{e \in Q'} \mc{E}^e$
such that $\phi \in \mc{Q}_j$.
Let $U$ be the set of vertices touched by $\mc{E}$.
As $\phi$ is $M^*$-compatible, each $e \in \phi(Q)$ 
has $\phi^e \in A(\Phi) \cap M^*$
with $e \in \phi^e(Q)$ and
$\pi_e \phi^e = \pi_e \phi = id$.

Next we specify the combinatorial structure 
of a potential cascade for $\phi$.
For the base of the cascade,
we fix a $\Phi$-embedding $\phi^c$ of $[q](p)$
with $\phi^c\ id_{[q]} = \phi$ and
$Im(\phi^c) \sm Im(\phi)$ disjoint from $U$,
such that for all $\phi^c\psi \in A(\Phi)^\le_r$ 
where $\psi \in [q](p)_r$ we have
$\gG(\phi^c\psi) \le_\gG G$.
For the absorbers in the cascade,
for each $\phi' \in \Ups'$ we fix
any $\Phi$-embedding $\phi^{\phi'}$ of $[q](Ker)$
with $\phi^{\phi'} \phi'{}^0 = \phi^c \phi'$,
such that for all $\phi^{\phi'}\psi  \in A(\Phi)^\le_r$
where $\psi \in [q](Ker)_r$ we have
$\gG(\phi^{\phi'}\psi) \le_\gG G$.
If $\phi'$ is some $\phi'_e$ with 
$\phi(Im(\phi'_e) \cap [q]) = e \in \phi(Q)$
we have the additional constraint
$\phi^{\phi'_e} \phi^{a_e} = \phi^e$,
where $a_e \in Ker^r$ is such that 
$Im(\pi_e) \sub Im(\phi^{a_e})$.
We choose absorbers as disjointly as possible,
i.e.\ with `private vertices' $I_{\phi'}$ that
are pairwise disjoint and disjoint from $U \cup Im(\phi^c)$,
where $I_{\phi'} = Im(\phi^{\phi'}) \sm Im(\phi^c\phi')$
if $\phi'$ is not some $\phi'_e$
or $I_{\phi'_e} = Im(\phi^{\phi'_e}) \sm 
( Im(\phi^c\phi'_e) \cup Im(\phi^e) )$.

As $(\Phi,\gG[G]^A)$ is $(\oO,h)$-extendable,
the number of such choices
for $\phi^c$ and $\phi^{\phi'}$ given $\phi$ and $\mc{E}$
is at least $0.9\oO n^{q(p-1)+v_+}$, where 
$v_+ = \sum_{\phi'} |I_{\phi'}| = p^r q(p^{q-r}-1)-Q(q-r)$.

Now we consider the algebraic constraints
that must be satisfied for the cascade
described above to appear in the template.
We condition on $f_j\phi^c$
such that $Im(\phi^c)$ is $j$-generic for $Im(\phi)$
and define $w_{\phi'}=f_j\phi^c\phi'$ for $\phi' \in \Ups'$.
Then each $\dim(w_{\phi'})=q$.
Now $\phi^c$ will define a cascade
$C_{\phi^c} = \sum \{\TT^{\phi^c\phi'(Q)}: \phi' \in \Ups'\}$
with each $\TT^{\phi^c\phi'(Q)}=\phi^{w_{\phi'}}(K^r_q(Ker))
= \phi^{\phi'}(K^r_q(Ker))$ as in Definitions 
\ref{def:absorb} and \ref{def:cascade} if 
\begin{enumerate}
\item
$f_j\phi^{\phi'}((i,a)) = (w_{\phi'})_i + a \cdot w_{\phi'}$ 
for each $\phi' \in \Ups'$, $i \in [q]$, $a \in Ker$, and
\item
$\phi^{\phi'} \phi^a$ is activated,
$A_{\phi^{\phi'} \phi^a} = A$, and
$T_e=j$ and $\pi_e \phi^{\phi'} \phi^a=id$ 
for all $\phi' \in \Ups'$, $a \in Ker^r$,
$e \in \phi^{\phi'} \phi^a(Q)$.
\end{enumerate}
We have the same bound on the probability of these
events as in \cite{Kexist}, as the estimates there
were sufficiently crude to absorb the extra factor
of $|\mc{A}|^{p^{r(q-r+1)}}$ here
for the events $A_{\phi^{\phi'} \phi^a} = A$.
The proof of concentration of the number of cascades
is also the same (the effect of changing $A_\phi$
is analysed in the same way as the effect of
changing whether $\phi$ is activated). \qed

\medskip

The proof of \cite[Lemma \oldhole+]{Kexist}
(the cascade random greedy algorithm)
now applies to show that Lemma \ref{hole+/bddint} 
follows from the following lemma.

\begin{lemma} \label{hole/bddint}
Suppose $\Psi^0 \in \mb{Z}^{\mc{A}(\Phi)}$ is $c'_2$-bounded
with $0 \le_\gG \pl^\gG \Psi^0 \le_\gG \pl^\gG M^*$.
Let $S = (\pl^\gG \Psi^0)^\circ$ and
suppose $M^*(S)$ is a set.
Then there are $M^\pm \sub \mc{A}(\Phi)$
such that every $\phi \in M^+$ is cascading, 
$M^*(\sum_{\phi \in M^+} \phi(Q))$ 
is a set and $3c_4$-bounded,
and $\pl^\gG M^+ = \pl^\gG M^- + \pl^\gG \Psi^0$.
\end{lemma}

We will prove Lemma \ref{hole/bddint} in section \ref{sec:cea}.

\section{Clique Exchange Algorithm} \label{sec:cea}

In this section we define our Clique Exchange Algorithm,
which has three applications in this paper,
namely to the proofs of Lemmas
\ref{avoid/bddint} and \ref{hole/bddint} (in this section)
and Lemma \ref{bddrat:L} (in section \ref{sec:bddint}).
The following lemma will allow us to modify an 
integer decomposition so as to avoid unforced uses of bad sets.
For the statement we recall the lattice constant $C_0=C_0(\gG)$
from Lemma \ref{AAspan}, and that $C_0=1$ if $\gG$ is elementary.

\begin{lemma} \label{avoid/bddint}
Let $\mc{A}$ be a $\Ss^\le$-family 
with $\Ss \le S_q$ and $|\mc{A}| \le K$
and suppose $\gG \in (\mb{Z}^D)^{\mc{A}_r}$.
Let $\Phi$ be an $(\oO,h)$-extendable
$\Ss$-adapted $[q]$-complex on $[n]$, 
where $\oO<\oO_0(q,D,K)$ and $n>n_0(q,D,K)$.
Suppose $\Psi^0 \in \mb{Z}^{\mc{A}(\Phi)}$ 
is $\tT$-bounded with $n^{-1/2} < \tT < \oO^4$.
Let $B^k \sub \Phi^\circ_k$ be $\eta$-bounded for $r \le k \le q$,
where $\eta = (9q)^{-2q} \oO$.
Then there is some $M_2 \tT$-bounded
$\Psi \in \mb{Z}^{\mc{A}(\Phi)}$,
where $M_2 = q^{(2q)^{2q}} C_0^2 \oO^{-2}$, 
with $\pl^\gG \Psi = J := \pl^\gG \Psi^0$ such that
\begin{enumerate}
\item if $k>r$ then $U(\Psi)_\psi \le 1$ 
for all $\psi \in \Phi_k$,
and $U(\Psi)_\psi = 0$ if $Im(\psi) \in B^k$,
\item $U(\Psi)_\psi \le U(J)_\psi + C_0+1$ 
for all $\psi \in \Phi_r$,
and $U(\Psi)_\psi = U(J)_\psi$ 
if $Im(\psi) \in B^r$.
\end{enumerate}
\end{lemma}

Lemma \ref{avoid/bddint} and Lemma \ref{bddint}
immediately imply the following lemma,
which will be used (with smaller $q$) in the inductive proof 
of Lemma \ref{bddint} in section \ref{sec:bddint}.

\begin{lemma} \label{bddint:avoid}
Let $\mc{A}$ be a $\Ss^\le$-family 
with $\Ss \le S_q$ and $|\mc{A}| \le K$
and suppose $\gG \in (\mb{Z}^D)^{\mc{A}_r}$.
Let $\Phi$ be an $(\oO,h)$-extendable 
$\Ss$-adapted $[q]$-complex on $[n]$, 
where $n^{-h^{-3q}}<\oO<\oO_0(q,D,K)$ and $n>n_0(q,D,K)$.
Suppose $J \in \bgen{\gG(\Phi)}$ is $\tT$-bounded
with $n^{-(5hq)^{-r}} < \tT < \oO^4$. 
Suppose $B^k \sub \Phi^\circ_k$ is $\eta$-bounded for 
$r \le k \le q$, where $\eta = (9q)^{-2q} \oO$.
Then there is some $\oO_q^{-3h} \tT$-bounded
$\Psi \in \mb{Z}^{\mc{A}(\Phi)}$ 
with $\pl \Psi = J$ such that
\begin{enumerate}
\item if $k>r$ then $U(\Psi)_\psi \le 1$ 
for all $\psi \in \Phi_k$,
and $U(\Psi)_\psi = 0$ if $Im(\psi) \in B^k$,
\item $U(\Psi)_\psi \le U(J)_\psi + C_0+1$ 
for all $\psi \in \Phi_r$,
and $U(\Psi)_\psi = U(J)_\psi$ if $Im(\psi) \in B^r$.
\end{enumerate}
\end{lemma}

\subsection{Splitting Phase}

Now we start the proof of Lemma \ref{avoid/bddint}.
Suppose $\Psi^0 \in \mb{Z}^{\mc{A}(\Phi)}$ is $\tT$-bounded.
We will obtain the desired $\Psi$ 
by an algorithm similar to that in \cite{Kexist}
(but with several significant differences).

To define the first phase of the algorithm,
we recall the $K^r_q$-decompositions  $\Ups$ and $\Ups'$ 
of $\OO=K^r_q(p)$ given by Lemma \ref{exchange},
and write $\OO' = K^r_q(p) \sm Q$.

\begin{alg} (Splitting Phase)
Let $(\phi_i: i \in |\Psi^0|)$ be any ordering 
of the signed elements of $\Psi^0$, i.e.\
$\Psi^0 = \sum_i s_i \{ \phi_i \}$ with each $s_i \in \pm 1$
and $\phi_i \in A^i(\Phi)$ for some $A^i \in \mc{A}$.
We apply a random greedy algorithm
to choose $\phi^*_i \in X_{E_i}(\Phi)$
for each $i$, where $E_i = ([q](p),[q],\phi_i)$. 
We say $\phi^*_i$ uses $e \in \Phi^\circ$ 
if $e \sub Im(\phi^*_i \phi)$ for some 
$\phi \in \Ups' \cup \Ups$
and $e \sm Im(\phi_i) \ne \es$.
Let $F_i$ be the set of used $e \in \Phi^\circ_r$.
We choose $\phi^*_i \in X_{E_i}(\Phi)$ uniformly at random 
subject to not using $F_i$ or $B=\cup_k B^k$.
\end{alg}

\begin{lemma} \label{split:avoid}
whp Splitting Phase does not abort
and $F_{|\Psi^0|}$ is $M_1 \tT$-bounded,
where $M_1 = 2^{r+3} (pq)^q \oO^{-1}$.
\end{lemma}

\nib{Proof.}
For $i \in [|\Psi^0|]$ we let $\mc{B}_i$ be the bad event 
that $F_i$ is not $M_1 \tT$-bounded. Let $\tau$ be the 
smallest $i$ for which $\mc{B}_i$ holds or the algorithm 
aborts, or $\infty$ if there is no such $i$. 
It suffices to show whp $\tau=\infty$. 
We fix $i_0 \in [|\Psi^0|]$ and bound 
$\mb{P}(\tau=i_0)$ as follows. 

We claim that for any $i<i_0$ the restrictions on $\phi^*_i$ 
forbid at most half of the possible choices 
of $\phi^*_i \in X_{E_i}(\Phi)$. To see this, 
first note that $X_{E_i}(\Phi) > \oO n^{pq-q}$
as $\Phi$ is $(\oO,s)$-extendable.
As $F_i$ is $M_1 \tT$-bounded
and each $B^k$ is $\eta$-bounded,
at most $(pq)^q (q\eta + M_1 \tT) n^{pq-q}$
choices use $F_i \cup B$; the claim follows. 

Now for each $e \in \Phi^\circ_r$ let 
$r_e = \sum_{i<i_0} \mb{P}'(e \in \phi^*_i(\OO'))$,
where $\mb{P}'$ denotes conditional probability
given the choices made before step $i$.

For any $i<i_0$, writing $r'=|e \sm Im(\phi_i)|$, 
there are at most $(pq)^q n^{pq-q-r'}$ choices of $\phi^*_i$
such that $e \in \phi^*_i(\OO')$, so by the claim
$\mb{P}'(e \in \phi^*_i(\OO')) < 2(pq)^q \oO^{-1} n^{-r'}$.
Also, given $r' \in [r]$, as $\Psi^0$ 
is $\tT$-bounded there are at most
$2\tbinom{r}{r'} \tT n^{r'}$ choices of $i$
such that $|e \sm Im(\phi_i)|=r'$.
Therefore $r_e < 2^{r+2} (pq)^q \oO^{-1} \tT$.

Now fix any $f \in \Phi^\circ_{r-1}$ and 
write $U(F_i)_f =  \sum_{i<i_0} X_i$,
where $X_i$ is the number of $e \in \Phi^\circ_r$
with $f \sub e \in \phi^*_i(\OO')$.
Then each $|X_i|<|\OO'|$ and $\sum_{i<i_0} \mb{E}'X_i
= \sum \{r_e: f \sub e\} 
< 2^{r+2} (pq)^q \oO^{-1} \tT n$,
so by Lemma \ref{dom} whp 
$U(F_i)_f <  M_1 \tT  n$,
so $F_i$ is $M_1 \tT$-bounded, so $\tau>i_0$. 
Taking a union bound over $i_0$, 
whp $\tau=\infty$, as required. \qed

\medskip

We let\footnote{
Note that (e.g.) $A^i(\Phi[\phi^*_i \Ups])
= \{ \phi^*_i \psi \in A^i(\Phi): \psi \in \Ups \}$. }
$\Psi^1 = \Psi^0 + \sum_{i \in [|\Psi^0|]} s_i
( A^i(\Phi[\phi^*_i \Ups']) - A^i(\Phi[\phi^*_i \Ups]) )$.
Then $\pl^\gG \Psi^1 = \pl^\gG \Psi^0 = J$, and all signed elements
of $\Psi^0$ are cancelled, so $\Psi^1$ is supported 
on maps added during Splitting Phase.

We classify maps added during Splitting Phase as near or far, 
where the near maps are those of the form $\phi^*_i \phi$ 
for $\phi \in \Ups'$ with $|Im(\phi) \cap [q]|=r$.
Also, for each pair $(O,\phi')$ where $\phi'$
is added during Splitting Phase,
$O \in \Phi_r/\Ss$ and $O \sub \phi'\Ss$,
we call $(O,\phi')$ near if $\phi'=\phi^*_i \phi$ is near 
and $Im(O) = \phi^*_i(Im(\phi) \cap [q])$, 
otherwise we call $(O,\phi')$ far.
We fix orbit representatives $\psi^O \in O$
and say that $(O,\phi')$ has type $\tT$ where
$\gG(\phi')^O = \gG[\psi^O]^\tT = \gG(\psi^O \tT^{-1})$
(we fix any such $\tT$ for each $\gG$-atom at $O$).
We also classify maps and near pairs as positive 
or negative according to their sign in $\Psi^1$.

Note that for each orbit $O$ such that there is some
far pair on $O$ there are exactly two such far pairs 
$(O,\phi^\pm)$ and $\gG(\phi^-)^O = -\gG(\phi^+)^O$.
For each $O \in \Phi_r/\Ss$ we let $\Psi^O$ 
be the sum of all $\pm \{\phi\}$ where
$\pm (O,\phi)$ is a signed near pair in $\Psi^1$.
Then $(\pl^\gG \Psi^O)^O = (\pl^\gG \Psi^1)^O = J^O$.

\subsection{Grouping Phase}

Now we will organise the near pairs 
on each $O$ into some cancelling groups
and $U(J)_O$ ungrouped near pairs.
To do so we will introduce some additional near pairs
in which we add and subtract some given element of $\mc{A}(\Phi)$
(which has no net effect on $\Psi^1$).

Consider any orbit $O$ with $\psi^O \in \Phi_B$, $B \in Q$. 
As $J = \pl^\gG \Psi^0 \in \bgen{\gG(\Phi)} 
= \mc{L}_\gG(\Phi) \sub \mc{L}^-_\gG(\Phi)$, we have 
$f_B(J)_{\psi^O} \in \gG^B = \bgen{ \gG^\tT: \tT \in \mc{A}_B }$.
By definition of $U(J)_O$ we can express $J^O$
as a sum of $U(J)_O$ signed $\gG$-atoms,
i.e.\ $f_B(J)_{\psi^O} = \sum_\tT n^O_\tT \gG^\tT$,
so $J^O = \sum_\tT n^O_\tT \gG(\psi^O\tT^{-1})$,
where $n^O \in \mb{Z}^{\mc{A}_B}$ with $|n^O|=U(J)_O$.
Let $m^{O\pm}_\tT$ be the number of near pairs
on $O$ of type $\tT$ of each sign.
As $J^O = (\pl^\gG \Psi^O)^O$ we have
$f_B(J)_{\psi^O}= \sum_\tT (m^{O+}_\tT-m^{O-}_\tT)\gG^\tT$,
so $n^O - m^{O+} + m^{O-} \in Z_B(\gG)$.
By Lemma \ref{AAspan} we have
$n^{Oj} \in Z_B(\gG)$ for $j \in [t^O]$ 
for some $t^O \le C_0 (|n^O|+|m^{O+}|+|m^{O-}|)$ 
with each $|n^{Oj}| \le C_0$
and $n^O - m^{O+} + m^{O-} = \sum_{j \in [t^O]} n^{Oj}$.

We will assign the near pairs to
cancelling groups and ungrouped near pairs
so that for each such $O$ and $\tT$, 
there are $|n^O_\tT|$ (correctly signed) 
ungrouped near pairs on $O$ of type $\tT$,
and the $j$th group has $|n^{Oj}_\tT|$ such near pairs. 

Let $d^O_\tT = (\sum_j |n^{Oj}_\tT|) - m^{O+}_\tT
= ((|n^O_\tT| + \sum_j |n^{Oj}_\tT|) 
- (m^{O+}_\tT+m^{O-}_\tT))/2$. If $d^O_\tT>0$ then 
we need to introduce $d^O_\tT$
new near pairs of type $\tT$ on $O$ 
in the Grouping Phase below.
If $d^O_\tT \le 0$ then we do not need to introduce
any new near pairs of type $\tT$ on $O$.
If $d^O_\tT<0$ then we have $2|d^O_\tT|$
unassigned near pairs of type $\tT$ on $O$,
with which we form $|d^O_\tT|$
additional cancelling groups each containing
one positive and one negative near pair.

\begin{alg} (Grouping Phase)
Let $S^J = \{ (O^i,\tT^i): i \in [|S^J|] \}$
be such that each $(O,\tT)$ with $d^O_\tT>0$ 
appears $d^O_\tT$ times.
We apply a random greedy algorithm 
to choose $\phi_i \in X_{E_i}(\Phi)$ with 
$E_i = (\ova{q},Im(\tT^i),\psi_{O^i} (\tT^i)^{-1})$.
We say $\phi_i$ uses $e \in \Phi^\circ_r$ 
if $Im(O^i) \ne e \sub Im(\phi_i)$.
Let $F'_i$ be the set of used $e \in \Phi^\circ_r$.
We choose $\phi_i$ uniformly at random subject to 
not using $F'_i \cup F_{|\Psi^0|} \cup B$.
\end{alg}

Similarly to Lemma \ref{split:avoid},
whp Grouping Phase does not abort
and $F'_{|S^J|}$ is $M_1 \tT$-bounded.
We create new near pairs by adding and subtracting 
each $\phi_i$, and then organise the near pairs
into cancelling groups and ungrouped near pairs
as described above. 

\subsection{Elimination Phase}

\def\oldOO*{6.15}

In the Elimination Phase we replace $\Psi^1$
by $\Psi^2$ so as to remove all cancelling groups
while preserving $\pl \Psi^2 = \pl \Psi^1 = J$.
We start by recalling \cite[Definition \oldOO*]{Kexist}.

\begin{defn} \label{def:OO*}
Let $\OO_1$ and $\OO_2$ be two copies of $\OO$.
Fix $f \in \Ups$ and $f' \in \Ups'$ 
with $|V(f) \cap V(f')|=r$.
For $j=1,2$ we denote the copies of
$\Ups$, $\Ups'$, $f$, $f'$ in $\OO_j$
by $\Ups_j$, $\Ups'_j$, $f_j$, $f'_j$.
Let $\OO^*$ be obtained by identifying
$\OO_1$ and $\OO_2$ so that $f'_1=f'_2$.
Let $\Ups^+ = \Ups_1 \cup (\Ups'_2\sm\{f'_2\})$
and $\Ups^- = \Ups_2 \cup (\Ups'_1\sm\{f'_1\})$.
Then $\Ups^+$ is a $K^r_q$-decomposition of $\OO^*$ containing $f_1$
and $\Ups^-$ is a $K^r_q$-decomposition of $\OO^*$ containing $f_2$.
\end{defn}

Next we introduce some notation for octahedra 
and their associated signed characteristic vectors.

\begin{defn} \label{def:oct} (octahedra)
Let $\Phi$ be an $R$-complex and $B \sub R$.
The $B$-octahedron is $O^B = B(2)$.
For $x \in [2]^B$ we define the sign
of $x$ by $s(x)=(-1)^{\sum_i (x_i-1)}$.
For $\psi \in O^B_B$ such that
$\psi(i)=(i,x_i)$ for all $i \in B$
we also write $s(\psi)=s(x)$.
Let $\mc{O}^B(\Phi)$ be the set of $\Phi$-embeddings of $O^B$.
For $\phi \in \mc{O}^B(\Phi)$ we let $\chi(\phi)$ denote the 
`signed characteristic vector' in $\mb{Z}^{\Phi_B}$, where 
$\chi(\phi)_{\phi \circ \psi} = s(\psi)$ for $\psi \in O^B_B$,
and all other entries of $\chi(\phi)$ are zero.
\end{defn}

The following definition and lemma 
implement octahedra as signed combinations of cliques.

\begin{defn}
For $x=(x_i) \in [s]^q$ we identify $x$ 
with the partite map $x:[q] \to [q] \times [s]$ 
where each $x(i)=(i,x_i)$,
and also with the image of this map.
We write $1$ for the map with all $1(i)=1$
and identify $[q]$ with $1([q])$.

For $e \in [q](s)_r$ let $X_e=\{x:e \sub x\}$.
Suppose $w \in \{-1,0,1\}^{[s]^q}$.

We say $e \in [q](s)_r$ is bad for $w$
if $|\{x \in X_e: w_x=1\}|>1$
or $|\{x \in X_e: w_x=-1\}|>1$. 

We say $w$ is simple if no $e$ is bad for $w$.
We define $\pl w \in \mb{Z}^{[q](s)_r}$
by $\pl w_e = \sum_{x \in X_e} w_x$.
\end{defn}

\begin{lemma} \label{octmove}
Let $s=(2q)^r$. Then for any $B \in Q$ there is a simple
$w^B \in \{-1,0,1\}^{[s]^q}$ with $\pl w^B = \chi(O^B)$.
Let $w^B_r$ denote the set of $e \in [q](s)_r$
such that $w^B_x \ne 0$ for some $x \in X_e$.
We can choose $w^B$ with $w^B_1=1$ so that
$w^B_r[V(O^B) \cup [q]] = O^B_B \cup Q$.
\end{lemma}

\nib{Proof.}
We start by setting $w^B_x = (-1)^{\sum_{i=1}^r (x_i-1)}$
if $x_i \in [2]$ for $i \in B$ 
and $x_i=1$ for $i \in [q] \sm B$,
otherwise $w^B_x=0$. Then $\pl w^B_e = \chi(O^B)$. We will 
repeatedly apply transformations to $w^B$ that
preserve $\pl w^B = \chi(O^B)$ until $w^B$ becomes simple.
Suppose $w^B$ is not simple. Fix $e$ bad for $w^B$
and $x,x' \in X_e$ with $w^B_x=1$ and $w^B_{x'}=-1$.
Fix a $[q](s)$-embedding $\phi$ of $\OO^*$
as in Definition \ref{def:OO*},
where $\phi(f_1)=x$, $\phi(f_2)=x'$ and 
if $a \in \phi(V(\OO^*)) \sm (x \cup x')$ 
then $w^B_y=0$ whenever $a \in y$.
We modify $w^B$ by adding 
$-1$ to each $\phi(g)$ where $g \in \Psi^+$ and 
$1$ to each $\phi(g')$ where $g' \in \Psi^-$.
This preserves $\pl w^B_e = \chi(O^B)$ and reduces the sum 
of $|w^B_x|$ over $x \in X_e$ with $e$ bad for $w^B$.
The process terminates with $w^B$ that is simple,
and we can relabel so that the other properties hold. \qed

\medskip

Let $w^B=[q](s)[w^B_r]$ be the $[q]$-complex
obtained by restricting $[q](s)$ to $w^B_r$.

\begin{alg} \label{alg:elim} (Elimination Phase)
Let $(C^i: i \in [P])$ be any ordering of
the cancelling groups, where each
$C^i = \{ (O^i,\phi^i_j): j \in [|C^i|] \}$
for some orbit $O^i$ with representative 
$\psi_{O^i} \in \Phi_{B^i}$
and each $\phi^i_j \in A^i_j(\Phi)$.
Our random greedy algorithm
will make several choices at step $i$.
First we choose $\psi^*_i \in X_{E^i}(\Phi)$
where $E^i = (B^i(2),B^i,\psi_{O^i})$;
we say that this choice has type 1.
Then for each $j \in [|C^i|]$ we make
type 2 choices $\phi^{*i}_j \in X_{E^i_j}(\Phi)$ 
where $E^i_j = (w^{B^i_j},F^i_j,\phi'{}^i_j)$,
$F^i_j = [q] \cup (B^i_j \times [2])$,
$B^i_j = \tT^i_j(B^i)$,
$\psi_{O^i} = \phi^i_j \tT^i_j$,
$\phi'{}^i_j\mid_{[q]}=\phi^i_j$,
$\phi'{}^i_j(\tT^i_j(x),y)=\psi^*_i(x,y)$
for $x \in B^i$, $y \in [2]$.

We let $\OO'_i = B^i(2)_r \sm \{id_{B^i}\}$ and 
$\OO^i_j = w^{B^i_j}_r \sm (B^i_j(2) \cup \ova{q}_r)$.
We say that $\psi^*_i$ uses $e \in \Phi^\circ_r$ 
with type 1 if $e = Im(\psi^*_i \psi)$
for some $\psi \in \OO'_i$
(we also write $e \in \psi^*_i(\OO'_i)$).
We say that $\phi^{*i}_j$ uses $e \in \Phi^\circ_r$ 
with type 2 if $e = Im(\phi^{*i}_j\psi)$
for some $\psi \in \OO^i_j$
or if $|e|>r$ and $e \sub Im(\phi^{*i}_j x)$
for some $x \in [s]^q \sm \{1\}$ with $w^{B^i_j}_x \ne 0$.

For $\aA=1,2$ let $F^\aA_i$ be the set of 
$e \in \Phi^\circ_r$ used with type $\aA$.
We make each choice at step $i$
uniformly at random subject to not using 
$F^1_i \cup F^2_i \cup F_{|\Psi^0|} \cup F'_{|S^J|} \cup B$.
\end{alg}

We will obtain $\Psi$ from $\Psi^1$ by adding
$\sum_x  w^{B^i_j}_x A^i_j(\Phi[\phi^{*i}_j x])$
for each $i \in [P]$ and $j \in [|C^i|]$
with the opposite sign to that 
of the near pair $(O^i,\phi^i_j)$.
This cancels all cancelling groups and
preserves $\pl^\gG \Psi = \pl^\gG \Psi^1 = J$
by the following lemma, which shows that
the construction for each cancelling group
has no total effect on $\pl^\gG \Psi$,
using $\sum_{j \in [|C^i|]} \gG^{\tT^i_j} = 0$.

\begin{lemma}
With notation as in Algorithm \ref{alg:elim}, we have
\[ \pl^\gG \sum_x  w^{B^i_j}_x A^i_j(\Phi[\phi^{*i}_j x]) 
= \sum \{ \chi(\psi^*_i)_e \gG(\phi^{*i}_j e (\tT^i_j)^{-1})
 : e \in \psi^*_i O^{B^i}_{B^i} \}. \]
\end{lemma}

\nib{Proof.}
By Lemma \ref{octmove}, we have  
$( \pl^\gG \sum_x  w^{B^i_j}_x A^i_j(\Phi[\phi^{*i}_j x]) )_\psi$
equal to zero unless $\psi\Ss=e\Ss$ 
with $e \in \psi^*_i O^{B^i}_{B^i}$, in which case, 
writing $e' = e \circ (\tT^i_j)^{-1}$, it equals
$(\pl w^{B^i_j}_{e'}) \gG(\phi^{*i}_j e')
= \chi(\psi^*_i)_e \gG(\phi^{*i}_j e')$. \qed

\medskip

Recall $M_1 = 2^{r+3} (pq)^q \oO^{-1}$,
$M_2 = q^{(2q)^{2q}} C_0^2 \oO^{-2}$ and note that
$M_2 \tT < \oO^{1.5}$ for $\oO<\oO_0(q,D,K)$.

\begin{lemma} \label{elim:avoid}
whp Elimination Phase does not abort
and $F^1_P \cup F^2_P$ is $M_2 \tT/2$-bounded.
\end{lemma}

\nib{Proof.}
For $i \in [P]$ we let $\mc{B}_i$ be the bad event 
that $F^1_i$ is not $2 C_0 M_1 \tT$-bounded
or $F^2_i$ is not $M_2 \tT/4$-bounded. Let $\tau$ be the 
smallest $i$ for which $\mc{B}_i$ holds or the algorithm 
aborts, or $\infty$ if there is no such $i$. 
It suffices to show whp $\tau=\infty$. 
We fix $i_0 \in [P]$ and bound 
$\mb{P}(\tau=i_0)$ as follows. 

Consider any $i<i_0$. 
Then $F^1_i$ is $2 C_0 M_1 \tT$-bounded
and $F^2_i$ is $M_2 \tT/4$-bounded.
As $\Phi$ is $(\oO,h)$-extendable, each
$X_{E_i}(\Phi) > \oO n^{v_{E_i}}$ and
$X_{E^i_j}(\Phi) > \oO n^{v_{E^i_j}}$.
As $F_{|\Psi^0|}$ and $F'_{|S^J|}$ are $M_1 \tT$-bounded,
and each $B^k$ is $\eta$-bounded,
at most half of the choices for $\psi^*_i$
or $\phi^{*i}_j$ are forbidden due to using 
$F^1_i \cup F^2_i \cup F_{|\Psi^0|} \cup F'_{|S^J|} \cup B$.

Next we fix $e \in \Phi^\circ_r$ and estimate the 
probability of using $e$ at step $i$ with each type.
For uses of type 1 there are at most
$2^r n^{v_{E_i} - |e \sm Im(O^i)|}$ choices of $\psi^*_i$
with $e \in \psi^*_i(\OO'_i)$, 
so $\mb{P}'(e \in \psi^*_i(\OO'_i)) 
\le 2^{r+1}\oO^{-1}n^{-|e \sm Im(O^i)|}$.
Similarly, for uses of type 2 we have
$\mb{P}'(e \in \phi^{*i}_j(\OO^i_j)) 
< 2|\OO^i_j|\oO^{-1}n^{-r^i_j(e)}$,
where $r^i_j(e)$ is the minimum $|e \sm Im(\psi')|$ 
with $\psi' \sub \phi^i_j$ or $\psi'=\psi^*_i\psi$
with $\psi \in B^i_j(2)$.

For any $r' \in [r]$, as $\Phi^0$ is $\tT$-bounded,
by construction of the cancelling groups in Grouping Phase,
there are at most $\tbinom{r}{r'} C_0 \tT n^{r'}$ choices of $i$
with $|Im(O^i) \sm e|=r'$, so
$r^1_e := \sum_{i<i_0} \mb{P}'(e \in \psi^*_i(\OO'_i)) 
< 2^{2r} \oO^{-1} C_0 \tT$.
Similarly, as $F_{|\Psi^0|}$ 
and $F'_{|S^J|}$ are $M_1 \tT$-bounded
and $F^1_i$ is $2 C_0 M_1 \tT$-bounded,
there are at most 
$4C_0 M_1 \tbinom{r}{r'} \tT n^{r'}$ 
choices of $i$ with $r^i_j(e)=r'$, so
$r^2_e := \sum_{i<i_0} \sum_{j \in [|C^i|]}
\mb{P}'(e \in \phi^{*i}_j(\OO^i_j)) 
< C_0 |\OO^i_j| 2^{r+3} C_0 M_1 \tT$.
By Lemma \ref{dom} we deduce whp
$F^1_i$ is $2 C_0 M_1 \tT$-bounded and
$F^2_i$ is $M_2 \tT/4$-bounded, so $\tau>i_0$. 
Taking a union bound over $i_0$, 
whp $\tau=\infty$, as required. \qed

\medskip

For any $\psi \in \Phi_{r-1}$ we have
$U(\Psi)_\psi \le U(\Psi^0)_\psi + |(F_{|\Psi^0|}
+ F'_{|S^J|} + F^1_P + F^2_P)\mid_\psi|
< M_2 \tT n$, so $\Psi$ is $M_2 \tT$-bounded.
For $k>r$ we avoided using any $\psi \in \Phi_k$
more than once or $B^k$ at all, so
$U(\Psi)_\psi \le 1$ for all $\psi \in \Phi_k$, 
and $U(\Psi)_\psi = 0$ if $Im(\psi) \in B^k$.
For any $\psi \in \Phi_r$ we have
a contribution of $U(J)_\psi$
from ungrouped near pairs to $U(\Psi)_\psi$.
If $\psi \in B^r$ there are no other uses,
so $U(\Psi)_\psi = U(J)_\psi$,
and otherwise there are at most\footnote{
The `+1' is only needed to account
for cancelling pairs in the case $C_0=1$.} 
$C_0+1$ other uses by a cancelling group, 
so $U(\Psi)_\psi \le U(J)_\psi + C_0+1$.
This completes the proof of Lemma \ref{avoid/bddint}. 

\subsection{Proof of Lemma \ref{hole/bddint}}

The proof of Lemma \ref{hole/bddint}
is very similar to that of Lemma \ref{avoid/bddint},
so we will just show the necessary modifications.
We consider $\Psi^0 \in \mb{Z}^{\mc{A}(\Phi)}$ 
that is $c'_2$-bounded, where $c'_2 = \oO^{-h^2} c_2$.
There are no bad sets $B$. We also suppose
$0 \le_\gG \pl^\gG \Psi^0 \le_\gG \pl^\gG M^*$
and $M^*(S)$ is a set, where
$S = (\pl^\gG \Psi^0)^\circ$.
We require the following definitions for Splitting Phase.

\begin{defn} \label{def:X*}
Consider any extension $E(\phi)=([q](p),[q],\phi)$ 
where $\phi \in A(\Phi)$ with $\gG(\phi) \le_\gG G$.
Let $H(\phi) = [q](p)_r \sm \ova{q}_r$.
We let $X^*_{E(\phi),H(\phi)}$ be the set or number of extensions 
$\phi^+ \in X_{E(\phi),H(\phi)}(\Phi,\gG[\pl^\gG M^*]^A)$ such that
\begin{enumerate}
\item
$\phi^+$ is {\em rainbow}: 
$j \ne j'$ whenever $\{\psi,\psi'\} \sub [q](p)_r \sm \ova{q}_r$,
$Im(\phi^+\psi) \in G^*_j$, $Im(\phi^+\psi') \in G^*_{j'}$, and
\item each $\phi^+ \phi'$ with $\phi' \in \Ups'$ 
is $M^*$-compatible if $|Im(\phi') \cap [q]|<r$
or $M^*$-compatible bar $\phi^+(e)$
if $Im(\phi') \cap [q] = e \in Q$.
\end{enumerate}
\end{defn}

We claim whp 
\begin{equation} \label{X*}
X^*_{E(\phi),H(\phi)} > \oO n^{pq-q} (\oO^2 z\rho/2)^{Qp^r}.
\end{equation}
Indeed, the proof of Lemma \ref{ext*} already gives 
rainbow extensions $\phi^+$, and by Remark \ref{rem:e|E}.ii
we can also require 
$\pi_e \phi^+ \psi = id$ and $A_{\phi^e}=A$
for all $\psi \in [q](p)_r \sm \ova{q}$, $e = Im(\phi^+ \psi)$,
which gives $\phi^+ \in X^*_{E(\phi),H(\phi)}$.

For the modified Splitting Phase, recalling that
$F_i$ is the set of used $e \in \Phi^\circ_r$,
we let $D_i = \cup_{e \in F_i} M^*(e)$,
and choose $\phi^*_i \in X^*_{E(\phi_i),H(\phi_i)}$ 
uniformly at random subject to 
$\phi^*_i(\OO') \cap (M^*(S) \cup D_i) = \es$.
Note that each $\phi^*_i(\OO') \sub G^*$ is rainbow,
so $M^*(\phi^*_i(\OO'))$ is a set. 

The modified form of Lemma \ref{split:avoid}
is to show whp $D_{|\Psi^0|}$ is $c_3$-bounded.
Accordingly, the bad event $\mc{B}_i$
is that $D_i$ is not $c_3$-bounded.
To see that at most half of the choices of
$\phi^*_i \in X^*_{E(\phi_i),H(\phi_i)}$ 
are forbidden we use (\ref{X*}), which gives
$X^*_{E(\phi_i),H(\phi_i)} 
> \oO n^{pq-q} (\oO^2 z\rho/2)^{Qp^r}
> 4(pq)^q c_3 n^{pq-q}$.

For $e \in G^*$ we define
$r_e = \sum_{i<i_0} \mb{P}'(e \in M^*(\phi^*_i(\OO')))
= \sum_{i<i_0} \sum_{e' \in M^*(e)}
 \mb{P}'(e' \in \phi^*_i(\OO'))$.
As $\Psi^0$ is $c'_2$-bounded,
there are at most $c'_2 \tbinom{r}{r'} n^{r'}$
choices of $i$ with $|e' \sm Im(\phi_i)|=r'$,
each $\mb{P}'(e' \in \phi^*_i(\OO')) 
< 2r!|\OO'|\oO^{-1} (\oO^2 z\rho/2)^{-Qp^r} n^{-r'}$,
so $r_e < (pq)^q \oO^{-1} 
(\oO^2 z\rho/2)^{-Qp^r} 2^{r+1} c'_2$.
We conclude that whp no $\mc{B}_i$ occurs,
so whp $D_{|\Psi^0|}$ is $c_3$-bounded.

Defining $\Psi^1$ as before, we have
$\pl^\gG \Psi^1 = \pl^\gG \Psi^0 = J$
and $\Psi^1$ is supported 
on maps added during Splitting Phase,
which are now rainbow in $G^*$ and
$M^*$-compatible bar at most one edge.

As before, we classify maps $\phi'$ and pairs $(O,\phi')$ 
added in Splitting Phase as near/far and positive/negative,
and assign types to pairs. As $\gG$ is now elementary,
the next part of the algorithm becomes simpler.
Indeed, for each $O \in \Phi_r/\Ss$ we can group 
the near pairs on $O$ into cancelling groups
of size one (near pairs of type zero)
or size two (of the same type and opposite sign),
and at most one additional positive near pair $(O,\phi^O)$,
which we call `solo', where if $Im(O) = e \in S$
with $\gG(e) \ne 0$ then $\phi^O(Q) \sub G^*$,
$\phi^O$ is $M^*$-compatible bar $e$,
and $\gG(e)=\gG(\psi)$ with $\psi \sub \phi^O$.

In Grouping Phase we only need to consider 
the solo near pairs, which we denote by
$\{ (O^i,\phi^{O^i}): i \in [s'] \}$.
We let $e_i = Im(O^i) \in G^*$, $A^i = A_{\phi^{e_i}}$,
and $\tT^i \in A^i_{B^i}$ be such that
$\psi'_i := \psi_{O^i} (\tT^i)^{-1} 
= \pi_{e_i}^{-1} \in A^i$, so $\gG(\psi'_i)=\gG(e_i)$.
Writing $D'_i = \cup_{i'<i} M^*(\phi_i(Q)))$,
we choose $\phi_i \in X_{E_i}(\Phi)$ with 
$E_i = (\ova{q},Im(\tT^i),\psi_{O^i} (\tT^i)^{-1})$
uniformly at random subject to $(\phi^*_i(Q) \sm \{e_i\})
 \cap  (M^*(S) \cup D_{|\Psi^0|} \cup D'_i) = \es$ 
and $\phi_i \in \mc{Q}^*$ being cascading.
Similarly to Lemma \ref{elim:cascade} (see below),
there are at least $0.9(\oO/z)^{q-r}$
choices for each $\phi_i$,
and similarly to Lemma \ref{split:avoid}
whp $D'_{s'}$ is $c_3$-bounded.

The next definition and accompanying lemma
set up the notation for the Elimination Phase
and show that there are whp many choices for each step.

\begin{defn} \label{def:elim2}
Given $A \in \mc{A}$, $B \in Q$,
$\psi \in A(\Phi)^\le_B  = \Phi_B$ 
we let $E(\psi) = (B(2),B,\psi)$ and 
$H(\psi) = B(2)_B \sm \{id_B\}$.

Suppose $\psi \sub \phi \in A(\Phi)$ 
with $\phi(Q) \sm Im(\psi)$ rainbow in $G^*$
and $\phi$ is $M^*$-compatible bar $Im(\psi)$.
Let $X^\phi_{E(\psi),H(\psi)}$ be the set or number of
$\psi^* \in X_{E(\psi),H(\psi)}(\Phi,\gG[\pl^\gG M^*]^A)$ 
such that 
\begin{enumerate}
\item $Im(\psi^*) \cap Im(\phi) = Im(\psi)$,
\item $Q^* := (\phi(Q) \cup \{ Im(\psi^*\psi'): \psi' \in B(2)_B \} )
\sm \{Im(\psi)\}$ is rainbow in $G^*$, and 
\item for all $e=Im(\psi^*\psi')$ with $\psi' \in B(2)_B \sm \{id_B\}$ 
we have $A_{\phi^e}=A$ and $\pi_e \psi^*\psi' = id$.
\end{enumerate}

For $\psi^* \in X^\phi_{E(\psi),H(\psi)}$
we let $E^\phi_{\psi^*} = (w^B,F,\psi^* \cup \phi)$,
where $F = [q] \cup (B \times [2])$.

We let $H^\phi_{\psi^*} = w^B \sm w^B[F]$
and $v_c := v_{E^\phi_{\psi^*}}$.

Let $X^c(E^\phi_{\psi^*} )^\pm$ be the set or number of 
$\phi^+ \in X_{E^\phi_{\psi^*},H^\phi_{\psi^*}}(\Phi,\gG[\pl^\gG M^*]^A)$ 
that are `rainbow $\Ups^\pm$ cascading', 
i.e.\ $\phi^+ x \in \mc{Q}^*$ is cascading for all 
$x \in \Ups^\pm := \{ x \in [s]^q \sm \{1\}: w^B_x = \pm 1\}$,
and $j \ne j'$ whenever $\{x,x'\} \sub \Ups^\pm$ with
$\phi^+ x \in \mc{Q}_j$, $\phi^+ x' \in \mc{Q}_{j'}$.
\end{defn}

\begin{lemma} \label{elim:cascade}
For $\psi$, $\phi$, $\psi^*$ as in Definition \ref{def:elim2}
whp $X^c(E^\phi_{\psi^*} )^\pm > (\oO/z)^{3Q^2 s^r} n^{v_c}$.
\end{lemma}

The proof of Lemma \ref{elim:cascade}
requires the following analogue of Lemma \ref{e|E}.

\begin{lemma} \label{Q*|eE}
Let $S \sub \Phi^\circ_r$ with $|S|<h=z$
and $\mc{E} = \cap_{f \in S} \mc{E}^f$.
Suppose $\phi \in A(\Phi)$ with $\gG(\phi) \le_\gG G$
such that $|\phi(Q) \cap S| \le 1$ and each
$e' \in \phi(Q) \sm S$ is not touched by $\mc{E}$.
Let $j \in [z]$ be such that
$T_{e'} \ne j$ for all $e' \in S \sm \phi(Q)$.
If $\phi(Q) \cap S = \{e\}$ suppose also that 
$\pi_e\phi=id$, $e \in G^*_j$, $\phi^e \in A(\Phi)$.
Then $\mb{P}( \phi \in \mc{Q}_j \mid \mc{E} ) > (\oO/z)^{3Q^2}$.
\end{lemma}

\nib{Proof.}
Let $1_e$ be $1$ if $\phi(Q) \cap S = \{e\}$ or $0$ otherwise.
For $e' \in \phi(Q) \sm S$ let $\pi^0_{e'}:e' \to [q]$
be such that $\pi^0_{e'} \phi = id$.
For each $e' \in \phi(Q) \sm S$ we fix 
$\phi_0^{e'} \in A(\Phi)$ with $\pi^0_{e'}\phi_0^{e'}=id$
and estimate the probability that all such
$e' \in G^*_j$ with $\phi^{e'} = \phi_0^{e'}$. 
Let $U$ be the set of vertices touched by $\mc{E}$.
As $(\Phi,\gG[G]^A)$ is $(\oO,h)$-extendable, there are at least
$(1-O(n^{-1})) \oO n^{(Q-1_e)(q-r)}$ choices for all $\phi_0^{e'}$
such that the sets $Im(\phi_0^{e'}) \sm e'$ are pairwise
disjoint and disjoint from $Im(\phi) \cup U$, and for each
$e' \in \phi(Q) \sm S$ and $\psi \sub \phi_0^{e'}$ with 
$Im(\psi) \in \phi_0^{e'}(Q) \sm \{e'\}$ we have $\gG(\psi) \le_\gG G$.
The probability that $\phi_0^{e'}$ is activated, $A_{\phi_0^{e'}}=A$,
$T_f=j$ and $\pi_f\phi_0^{e'}=id$ for all such $e'$ and $f \in \phi_0^{e'}$ 
is at least $((z(q)_r)^{-Q} |\mc{A}|^{-1} \oO^2 )^{Q-1_e}$.

We condition on $f_j\!\mid_{Im(\phi)}$ such that
$\dim(f_j\phi)=q$; this occurs with probability $1-O(n^{-1})$.
For each $e' \in \phi(Q) \sm S$ there is a unique 
$y^{e'} \in \mb{F}_{p^a}^r$ such that 
$(My^{e'})_i = f_j\pi_{e'}^{-1}(i)$ for all $i \in Im(\pi_{e'})$.
With probability $(1+O(n^{-1}))(p^{-a})^{(q-r)(1_e-Q)}$
we have $f_j(\phi^{e'}(i))=(My^{e'})_i$ 
for all such $e'$ and $i \in [q] \sm Im(\pi_{e'})$.
Therefore $\mb{P}(\cap_{e'} \{ \phi^{e'} = \phi_0^{e'} \} \mid \mc{E})
> (1+O(n^{-1})) [ (z(q)_r |\mc{A}|^{-1} \oO^{-2} )^Q p^{a(q-r)} ]^{1_e-Q}$.
Summing over all choices for $\phi_0^{e'}$ gives 
$\mb{P}( \phi \in \mc{Q}_j \mid \mc{E} ) > (\oO/z)^{3Q^2}$. \qed

\medskip

\nib{Proof of Lemma \ref{elim:cascade}.}
As $(\Phi,\gG[G]^A)$ is $(\oO,h)$-extendable, 
there are at least $\oO n^{v_c}$ choices of
$\phi^+ \in X_{E^\phi_{\psi^*},H^\phi_{\psi^*}}(\Phi,\gG[G]^A)$.
We fix any such $\phi^+$ and estimate 
$\mb{P}(\phi^+ \in X^c(E^\phi_{\psi^*} )^+)$
by repeated application of Lemma \ref{Q*|eE}. 
(The same estimates will apply
to $X^c(E^\phi_{\psi^*} )^-$.)
We consider sequentially each $\phi' \in \Ups^+$, 
and fix $j \in [z]$ distinct from all previous choices 
such that (recall $Q^*$ from Definition \ref{def:elim2})
if $Im(\phi') \cap Q^*= e$ then $e \in G^*_j$.

We let $\mc{E}$ be the intersection of all local events $\mc{E}^e$ 
where $e \in Q^*$ or $e \sub Im(\phi^+\phi'')$ 
for some previously considered $\phi'' \in \Ups^+$.
We discard $O(n^{v_c-1})$ choices of $\phi^+$
such that any $e' \in \phi^+\phi'(Q)$ is touched by $\mc{E}$.
Then $\mb{P}( \phi^+\phi' \in \mc{Q}_j \mid \mc{E} ) 
> (\oO/z)^{3Q^2}$ by Lemma \ref{Q*|eE}.
Multiplying all conditional probabilities and summing 
over $\phi^+$ gives $\mb{E}X^c(E^\phi_{\psi^*} )^+ 
> (1-O(n^{-1}))((\oO/z)^{3Q^2})^{s^r - 1} n^{v_c}$;
concentration follows from Lemma \ref{lip3}. \qed

\medskip

In the modified Elimination Phase, we recall that
the cancelling groups have size one (zero near pairs)
or two (cancelling pairs), say $C^i = \{ (O^i,\phi^i) \}$ 
or $C^i = \{ (O^i,\phi^i_+), (O^i,\phi^i_-) \}$.
We adopt the notation of Definition \ref{def:elim2}
and fix representatives $\psi_{O^i} \in O^i$
as in Algorithm \ref{alg:elim}.
 
If $C^i = \{ (O^i,\phi^i) \}$ we write\footnote{
The sign of a zero pair is irrelevant;
we fix + for convenient notation.}
$\phi^i_+ = \phi^i$, $\psi_{O^i} = \phi^i_+ \tT^i_+$
and $\psi^i_+ = \psi_{O^i} (\tT^i_+)^{-1} 
  \sub \phi^i_+ \in A^i(\Phi)$,
where $\gG(\psi^i_+)=0$, we choose
$\psi^+_i \in X^{\phi^i_+}_{E(\psi^i_+),H(\psi^i_+)}$
and then $\phi^+_i \in X^c(E^{\phi^i_+}_{\psi^+_i} )^+$.

If $C^i = \{ (O^i,\phi^i_+), (O^i,\phi^i_-) \}$,
we write $\psi_{O^i} = \phi^i_\pm \tT^i_\pm$
and $\psi^i_\pm = \psi_{O^i} (\tT^i_\pm)^{-1} 
  \sub \phi^i_\pm \in A^i_\pm(\Phi)$,
where $\gG(\psi^i_+)=\gG(\psi^i_-)$.
We note that
$\gG(\phi^i_\pm)-\gG(\psi^i_\pm) \le_\gG \pl^\gG M^*$,
let  $B^i_\pm = \tT^i_\pm(B^i)$ and choose 
$\psi^+_i \in X^{\phi^i_+}_{E(\psi^i_+),H(\psi^i_+)}$ such that
\[Q^\pm_i = (\phi^i_\pm(Q) \cup \{ Im(\psi^+_i\psi') :
 \psi' \in B^i_\pm(2)_r \} )  \sm Im(O^i)\] is rainbow 
in $G^*$ (this holds by definition for $Q^+_i$ but
is an extra requirement for $Q^-_i$).
We define $\psi^-_i \in X^{\phi^i_-}_{E(\psi^i_-),H(\psi^i_-)}$
by $\psi^-_i(\tT^i_-(x),y)=\psi^+_i(\tT^i_+(x),y)$ and
choose $\phi^\pm_i \in X^c(E^{\phi^i_\pm}_{\psi^\pm_i} )^\pm$.
 
We define type 1 and 2 uses similarly to before
and let $D^\aA_i = \cup_{\psi \in F^\aA_i} M^*(Im(\psi))$.
We make the above choices uniformly at random 
such that $M^*(Q^\pm_i)$ both avoid 
$M^*(S) \cup D_{|\Psi^0|} \cup D'_{|S^J|} \cup D^1_i \cup D^2_i$.

To modify Lemma \ref{elim:avoid}, we let $\mc{B}_i$
be the event that $D^1_i \cup D^2_i$ is not $c_4$-bounded.
To see that at most half of the choices for any
$\psi^+_i$ or $\phi^\pm_i$ are forbidden, we use
$X^{\phi^i_+}_{E(\psi^i_+),H(\psi^i_+)} >
 \oO n^{q-r} (\oO^2 z\rho/2)^Q > q^q c_4 n^{q-r}$ 
(this bound is similar to that in (\ref{X*}))
and $X^c(E^{\phi^i_\pm}_{\psi^\pm_i} )^\pm > 
   (\oO/z)^{3Q^2 s^r} n^{v_c}  > (qs)^q c_4 n^{v_c}$
(by Lemma \ref{elim:cascade}).

Then for $e \in G^*$ we have
\[r^1_e := \sum_{i<i_0} \mb{P}'(e \in M^*(\psi^+_i(\OO'_i)))
= \sum_{i<i_0} \sum_{e' \in M^*(e)}
 \mb{P}'(e' \in \psi^+_i(\OO'_i))
< 2r 2^{2r} \oO^{-1}(\oO^2 z\rho/2)^{-Q} c_3,\] 
recalling $\OO'_i = B^i(2)_r \sm \{id_{B^i}\}$, and
\[ r^2_e := \sum_{i<i_0} \sum_{\psi \in \GG_i}
\mb{P}'(Im(\psi) \in M^*(e))
< (qs)^q (\oO/z)^{-3Q^2 s^r} c_3,\]
where, writing
$\OO^i_\pm = w^{B^i_\pm}_r \sm (B^i_\pm(2) \cup \ova{q}_r)$,
we let $\GG_i = \phi^+_i \OO^i_+$ if $|C^i|=1$ or 
$\GG_i = \phi^+_i \OO^i_+ \cup \phi^i_- \OO^i_-$ if $|C^i|=2$.
As before, we deduce whp no $\mc{B}_i$ occurs,
so $D^1_P \cup D^2_P$ is $c_4$-bounded.

To conclude, we obtain $\Psi$ from $\Psi^1$ 
where for each $i \in [P]$ we add 
$\sum_x w^{B^i_+}_x A^i(\Phi[\phi^+_i x])$
with the opposite sign to $(O^i,\phi^i)$ if $|C^i|=1$,
or $\sum_x w^{B^i_+}_x A^i_+(\Phi[\phi^+_i x])
- \sum_x w^{B^i_-}_x A^i_-(\Phi[\phi^-_+ x])$ if $|C^2|=1$;
this cancels all cancelling pairs and
preserves $\pl^\gG \Psi = \pl^\gG \Psi^0$.
Also, for all positive maps $\phi$ added in Grouping Phase
and not cancelled, or added during Elimination Phase,
$\phi$ is cascading, $M^*(\phi(Q))$ is a set,
all such $M^*(\phi(Q))$ are disjoint,
their union is contained in $D^1_P \cup D^2_P$,
which is $c_4$-bounded and disjoint from $M^*(S)$.
This completes the proof of Lemma \ref{hole/bddint}. \qed

\section{Integral decomposition} \label{sec:int}

In this section we give a characterisation
of the decomposition lattice $\bgen{\gG(\Phi)}$,
which generalises the degree-type conditions 
for $K^r_q$-divisibility to the labelled setting.
The characterisation is given in the second subsection,
using a characterisation of the simpler auxiliary problem of 
octahedral decomposition, which is given
in the first subsection.

\subsection{Octahedral decomposition}

A key ingredient in the results of Graver and Jurkat \cite{GJ}
and Wilson \cite{W4} (generalised in \cite{Kexist})
is the decomposition of null vectors by octahedra.
In this subsection we establish an analogous
result for adapted complexes. 
We start by defining null vectors.
Throughout, $\Phi$ is a $\Ss$-adapted $R$-complex
and $\GG$ is a finite abelian group.

\begin{defn} (null)
For $J \in \GG^\Phi$ and $\psi \in \Phi$ 
we write $\pl J_\psi = \sum J\mid_{\psi}
= \sum \{ J_\phi: \psi \sub \phi \in \Phi \}$.

We define $\pl_i J \in \GG^{\Phi_i}$
by $(\pl_i J)_\psi = \pl J_\psi$ for $\psi \in \Phi_i$. 
We say $J$ is $i$-null if $\pl_i J = 0$.

For $J \in \GG^{\Phi_j}$ we write $\pl J = \pl_{j-1} J$;
we say $J$ is null if $\pl J=0$,
i.e.\ $J$ is $(j-1)$-null.
\end{defn}

Next we introduce the symmetric analogues of
octahedra and their associated signed characteristic vectors
(recall Definitions \ref{def:symm} and \ref{def:oct}).

\begin{defn}
Given $\psi^* \in \mc{O}^B(\Phi)$ and $v \in \GG^{\Ss^B}$
let $\chi(v,\psi^*)$ denote the 
`symmetric characteristic vector' in $(\GG^{\Ss^B})^{\Phi_B}$
where $\chi(v,\psi^*)_{\psi^*\psi\tau} = s(\psi) v \tau$ 
whenever $\psi \in O^B_B$, $\tau \in \Ss^B_B$,
and all other entries of $\chi(v,\psi^*)$ are zero.
For $\Psi \in (\GG^{\Ss^B})^{\mc{O}^B(\Phi)}$ we write 
$\pl \Psi = \sum_{\psi^*} \chi(\Psi_{\psi^*},\psi^*)$.
\end{defn}

We note the linearity 
$\chi(v+v',\psi^*) = \chi(v,\psi^*) + \chi(v',\psi^*)$,
which follows from $(v+v') \tau = v \tau + v' \tau$.
 
\begin{lemma}
If $\psi^* \in \mc{O}^B(\Phi)$ and $v \in \GG^{\Ss^B}$
then $\chi(v,\psi^*)$ is symmetric and null.
\end{lemma}

\nib{Proof.}
For $\psi \in O^B_B$ and $\tau, \tau' \in \Ss^B_B$
we have $\chi(v,\psi^*)_{\psi^*\psi\tau'} \tau 
= s(\psi) v \tau' \tau
= \chi(v,\psi^*)_{\psi^*\psi\tau'\tau}$,
so $\chi(v,\psi^*)$ is symmetric.
Also, for any $\psi^\pm \in O^B_B$
that agree on $\psi' \in O^{B'}_{B'}$ 
with $|B'|=|B|-1$ we have
$\pl \chi(v,\psi^*)_{\psi^*\psi'\tau}
= \chi(v,\psi^*)_{\psi^*\psi^+\tau}
+ \chi(v,\psi^*)_{\psi^*\psi^-\tau}
= s(\psi^+) v \tau + s(\psi^-) v \tau = 0$,
so $\chi(v,\psi^*)$ is null. \qed

\medskip

The following main lemma of this subsection
shows that groups of symmetric null vectors are
generated by symmetric characteristic vectors
of octahedra when $\Phi$ is extendable.

\begin{lemma} \label{oct:lattice}
Let $\Phi$ be an $(\oO,s)$-extendable $\Ss$-adapted $R$-complex
and $B \sub R$ with $|B|=r$, where $s=3r^2$,
$n = |V(\Phi)| > n_0(r,\GG)$ is large and $\oO>n^{-1/2}$.
Suppose $H$ is a symmetric subgroup of $\GG^{\Ss^B}$
and $J \in H^{\Phi_B}$ is symmetric and null.
Then $J = \pl \Psi$ for some 
$\Psi \in H^{\mc{O}^B(\Phi)}$.
\end{lemma}

It is convenient to first reduce the proof 
of Lemma \ref{oct:lattice} to the case $B=R$.

\begin{lemma} \label{B=R}
It suffices to prove Lemma \ref{oct:lattice} when $B=R$.
\end{lemma}

\nib{Proof.}
We reduce the general case of Lemma \ref{oct:lattice}
to the case $B=R$ as follows.
Let $\Phi'$ be the $B$-complex with
$\Phi'_{B'}=\Phi_{B'}$ for all $B' \sub B$.
Then $\Phi'$ is $\Ss^B_B$-adapted
and $(\oO',s)$-extendable, with 
$\oO' = \oO n/n' > n^{1/2}/n' > n'{}^{-1/2}$, 
where $n' = |V(\Phi')|$.

Suppose $B \in C \in \mc{P}^\Ss$. 
Let $X = \{x^{B'}: B' \in C\}$ 
where for each $B' \in C$ we fix 
any representative $x^{B'} \in \Ss^B_{B'}$.
Note that any $\sS \in \Ss^B$ has a unique expression
$\sS = \tau x$ with $\tau \in \Ss^B_B$, $x \in X$.

We define $\pi: (\GG^X)^{\Ss^B_B} \to \GG^{\Ss^B}$
by $\pi(v)_{\tau x} = (v_\tau)_x$.
For any set $Y$ we define
$\pi: ((\GG^X)^{\Ss^B_B})^Y \to (\GG^{\Ss^B})^Y$
by $\pi(w)_y = \pi(w_y)$ for all $y \in Y$.
Note that for any $v \in (\GG^X)^{\Ss^B_B}$
and $\tau' \in \Ss^B_B$ we have
$\pi(v \tau') = \pi(v) \tau'$;
indeed, for any $\tau \in \Ss^B_B$ and $x \in X$
we have $\pi(v \tau')_{\tau x} = ((v \tau')_\tau)_x
= (v_{\tau' \tau})_x = \pi(v)_{\tau' \tau x}
= (\pi(v) \tau')_{\tau x}$.
We let $H' = \{h': \pi(h') \in H\}$
and note that $H'$ is a symmetric subgroup
of $(\GG^X)^{\Ss^B_B}$.

Suppose $J \in H^{\Phi_B}$ is symmetric and null.
Define $J' \in (H')^{\Phi'_B}$
by $((J'_\psi)_\tau)_x = (J_\psi)_{\tau x}$.
Note that $\pi(J')=J$.

We claim that $J'$ is symmetric and null.
To see this, note that for any $\tau, \tau' \in \Ss^B_B$
and $x \in X$ we have $((J'_\psi \tau')_\tau)_x
  = ((J'_\psi)_{\tau'\tau})_x = (J_\psi)_{\tau'\tau x}
  = (J_\psi \tau')_{\tau x} = (J_{\psi\tau'})_{\tau x}
  = ((J'_{\psi\tau'})_\tau)_x$, 
i.e.\ $J'_\psi \tau' = J'_{\psi\tau'}$,
i.e.\ $J'$ is symmetric.
Also, for any $\psi' \in \Phi'_{B'}$
with $B' \sub B$, $|B'|=r-1$ 
and $\tau \in \Ss^B_B$, $x \in X$ we have
$((\pl J'_{\psi'})_\tau)_x
= \sum \{ ((J'_\psi)_\tau)_x : \psi' \sub \psi \}
= \sum \{ (J_\psi)_{\tau x} : \psi' \sub \psi \}
= (\pl J_\psi)_{\tau x} = 0$, so $J'$ is null, as claimed.

Now by the case $B=R$ of Lemma \ref{oct:lattice}
we have $J' = \pl \Psi'$ for some
$\Psi' \in (H')^{\mc{O}^B(\Phi)}$.
Let $\Psi=\pi(\Psi') \in H^{\mc{O}^B(\Phi)}$.
It remains to show that $\pl \Psi = J$,
i.e.\ for any $\psi \in \Phi_B$ that 
$\pl \Psi_\psi = \pi(\pl \Psi'_\psi)$.
It suffices to show for any $\psi^* \in \mc{O}^B(\Phi)$
that $\chi(\Psi_{\psi^*},\psi^*)_\psi
= \pi \chi(\Psi'_{\psi^*},\psi^*)_\psi$.
Let $\tau_0 \in \Ss^B_B$ be such that 
$\psi \tau_0^{-1} \in \psi^* O^B_B$. Then
$\chi(\Psi'_{\psi^*},\psi^*)_\psi = \Psi'_{\psi^*} \tau_0 $ and 
$\chi(\pi(\Psi'_{\psi^*}),\psi^*)_\psi = \pi(\Psi'_{\psi^*}) \tau_0
= \pi(\Psi'_{\psi^*} \tau_0)$, as required. \qed

\medskip

We will prove Lemma \ref{oct:lattice} (with $B=R$)
by induction on $r$; the proof of the 
following lemma uses the induction hypothesis.

\begin{lemma} \label{oct:focus}
Let $\Phi$ be an $(\oO,s')$-extendable $\Ss$-adapted 
$B$-complex and $B' \sub B$, with $|B|=r$, $|B'|=r'$, 
$s'=3(r-r')^2$ and $n = |V(\Phi)| > n_0(r,\GG)$ large
and $\oO>n^{-1/2}$. Suppose $B' \in C' \in \mc{P}^\Ss$
and $\Phi' \sub \Phi_{C'}$ is 
such that $(\Phi,\Phi')$ is $(\oO,2)$-extendable
and $\Phi[\Phi']$ is $\Ss$-adapted.
Suppose $H$ is a symmetric subgroup of $\GG^\Ss$ 
and $J \in H^{\Phi_B}$ is null and symmetric.
Then there is $\Psi \in H^{\mc{O}^B(\Phi)}$
with $J-\pl\Psi \in H^{\Phi[\Phi']_B}$.
\end{lemma}

\nib{Proof.}
We can assume $B' \ne \es$, otherwise 
the lemma holds trivially with $\Psi=0$.
Let $B^* = B \sm B'$,
$\Ss' = \{ \sS \in \Ss: \sS\mid_{B'}=id_{B'} \}$
and $\Ss^*=\Ss/B'=\{ \sS\mid_{B^*}: \sS \in \Ss' \}$.
Let $X = \{ x^C: C \in \Ss' \sm \Ss\}$
be a set of representatives of
the right cosets of $\Ss'$ in $\Ss$.
Then any $\sS \in \Ss$ has a unique
representation $\sS = \sS'x$
with $\sS' \in \Ss'$, $x \in X$.

We define $\pi: (\GG^X)^{\Ss^*} \to \GG^{\Ss}$
by $\pi(v)_{\sS' x} = (v_{\sS^*})_x$ whenever 
$x \in X$, $\sS' \in \Ss'$, $\sS^*=\sS'\mid_{B^*}$.
Note that for any $v \in (\GG^X)^{\Ss^*}$,
$\tau' \in \Ss'$, $\tau^*=\tau'\mid_{B^*}$
we have $\pi(v \tau^*) = \pi(v) \tau'$;
indeed, $\pi(v \tau^*)_{\sS' x} = ((v \tau^*)_{\sS^*})_x
= (v_{\tau^* \sS^*})_x = \pi(v)_{\tau' \sS' x}
= (\pi(v) \tau')_{\sS' x}$.
We let $H^* = \{h^*: \pi(h^*) \in H\}$
and note that $H^*$ is a symmetric subgroup
of $(\GG^X)^{\Ss^*}$.

Consider any $\psi^* \in \Phi_{B'} \sm \Phi'$.
Recall (Lemma \ref{nhood:adapt}) that 
$\Phi^*=\Phi/\psi^*$ is $\Ss^*$-adapted.
Define $J^* \in (H^*)^{\Phi^*_{B^*}}$ by
$((J^*_{\psi/\psi^*})_{\sS^*})_x = (J_\psi)_{\sS' x}$
whenever $\psi^* \sub \psi \in \Phi_B$, $x \in X$
and $\sS^* = \sS'\mid_{B^*}$ with $\sS' \in \Ss'$.
Note that $\pi(J^*_{\psi/\psi^*})=J_\psi$.

We claim that $J^*$ is symmetric and null.
To see this, consider any
$\psi^* \sub \psi \in \Phi_B$, $x \in X$ and
$\sS^* = \sS'\mid_{B^*}$, $\tau^* = \tau'\mid_{B^*}$ 
with $\sS', \tau' \in \Ss'$.
Then $((J^*_{\psi/\psi^*} \tau^*)_{\sS^*})_x 
  = ((J^*_{\psi/\psi^*})_{\tau^*\sS^*})_x 
  = (J_\psi)_{\tau'\sS' x} 
  = (J_\psi \tau')_{\sS' x} 
  = (J_{\psi\tau'})_{\sS' x} 
  = ((J^*_{\psi\tau'/\psi^*})_{\sS^*})_x 
  = ((J^*_{(\psi/\psi^*)\tau^*})_{\sS^*})_x$,
so $J^*$ is symmetric.
Also, for any $\psi'/\psi^* \in \Phi^*_{r-r'-1}$
we have $((\pl J^*_{\psi'/\psi^*})_{\sS^*})_x
  = \sum \{ ((J^*_{\psi/\psi^*})_{\sS^*})_x:
   \psi/\psi^* \in \Phi^*_{B^*}\mid_{\psi'/\psi^*} \}
  =  \sum \{ (J_\psi)_{\sS' x}:
   \psi \in \Phi\mid_{\psi'} \} = 0$,
so $J^*$ is null, as claimed. 

By the inductive hypothesis of Lemma \ref{oct:lattice}
we have $J^* = \pl \Psi^{\psi^*}$ for some
$\Psi^{\psi^*} \in (H^*)^{\mc{O}^{B^*}(\Phi^*)}$.
For each $\phi^* \in \mc{O}^{B^*}(\Phi^*)$
we consider the $\Phi$-extension 
$E^{\psi^*}_{\phi^*} = (O^B,F,\psi^*\cup\phi^*)$,
where $F = B \cup V(O^{B^*})$.
We construct $\Psi \in H^{\mc{O}^B(\Phi)}$
by letting $\psi^*$ range over a 
set of orbit representatives for
$(\Phi_{B'} \sm \Phi')/\Ss$, letting $\phi^*$ 
range over $\mc{O}^{B^*}(\Phi^*)$, and adding 
$\pi(\Psi^{\psi^*}_{\phi^*})  \{\phi\}$ to $\Psi$ for some
$\phi \in X_{E^{\psi^*}_{\phi^*}}(\Phi,\Phi')$.

To complete the proof, it suffices to show
that $(\pl \Psi)_\psi = J_\psi$ for any 
$\psi \in \Phi_B$ with $\psi^*\Ss \sub \psi\Ss$ for 
some representative $\psi^*$ used in the construction.
As $\pl \Psi$ and $J$ are both symmetric,
it suffices to prove this when $\psi^* \sub \psi$.
As $J_\psi = \pi(J^*_{\psi/\psi^*})$ and
$J^*_{\psi/\psi^*} = \pl \Psi^*_{\psi/\psi^*}$
it suffices to show that
$\pi( \chi(\Psi^{\psi^*}_{\phi^*},\phi^*)_{\psi/\psi^*} )
= \chi(\pi(\Psi^{\psi^*}_{\phi^*}),\phi)_\psi$
for any $\phi^*$ and $\phi$ as above.
Both sides are zero unless
$\psi (\tau')^{-1} \in \phi O^B_B$
for some $\tau' \in \Ss$.
Then $\tau' \in \Ss'$ as $\psi^* \sub \psi$.
Let $\tau^* = \tau'\mid_{B^*}$.
Then $(\psi/\psi^*)(\tau^*)^{-1} \in \phi^* O^{B^*}_{B^*}$.
We have $\chi(\Psi^{\psi^*}_{\phi^*},\phi^*)_{\psi/\psi^*}
 = \Psi^{\psi^*}_{\phi^*} \tau^*$
and $\chi(\pi(\Psi^{\psi^*}_{\phi^*}),\phi)_\psi
 = \pi(\Psi^{\psi^*}_{\phi^*}) \tau'
 = \pi( \Psi^{\psi^*}_{\phi^*} \tau^* )$, as required. \qed

\medskip

\nib{Proof of Lemma \ref{oct:lattice}.}
We can assume $B=R$ by Lemma \ref{B=R}.
For the purposes of induction we prove a slightly
stronger statement, in which we replace the assumption
$\oO>n^{-1/2}$ by the weaker assumption
$\oO>(2r)^{r^5}n^{-0.6}$.
Fix any $\Phi$-embedding $\psi_0$ of $O^B$. 
Let $V_0=Im(\psi_0)$ and 
$\tau:V(\Phi)\sm V_0 \to B$ be uniformly random.
For $(i,x_i) \in V(O^B)$ let $\tau(\psi_0(i,x_i))=i$.
Let $\Phi^\tau$ be the set of 
$\phi \in \Phi$ with $\tau \phi = id$.
Consider the $\Phi[\Phi^\tau]$-extension $E_0=(B(3),V(O^B),\psi_0)$. 
For $j \in [r]$ we let 
\[L^j = \bigcup \{ \psi^* e \Ss: \psi^* \in X_{E_0}(\Phi[\Phi^\tau]),
e \in B(3)_{B'}, B' \in \tbinom{B}{j} \}.\]
The main part of the proof lies in showing that
we can reduce the support of $J$ to $\Phi[L^r]_B$.

Before doing so, we start by making the support 
disjoint from $V_0$. We identify $B$ with $[r]$ 
and for each $j \in B$ we let $L'_j = \bigcup \{ \psi\Ss:
\psi \in \Phi_j, \psi(j) \notin V_0 \}$. We define 
$\Phi'{}^j$ for $0 \le j \le r$ by $\Phi'{}^0=\Phi$ 
and $\Phi'{}^j = \Phi'{}^{j-1}[L'_j]$ for $j \in [r]$. 
Then $\Phi' := \Phi'{}^r = \Phi[ V(\Phi) \sm V_0 ]$.

We claim that each $(\Phi'{}^{j-1},L'_j)$
is $(\oO - O(n^{-1}),s)$-extendable. To see this,
consider any $(\Phi'{}^{j-1},L'_j)$-extension
$(E,H')$ where $E=(H,F,\phi)$ is a $\Phi'{}^{j-1}$-extension 
and $H' \sub H \sm H[F]$. Note that
if $\phi^* \in X_E(\Phi'{}^{j-1})$ with
$Im(\phi^*) \cap V_0 = Im(\phi) \cap V_0$
then $\phi^* \in X_E(\Phi'{}^{j-1},L'_j)$.
Thus $X_E(\Phi'{}^{j-1},L'_j) > X_E(\Phi'{}^{j-1}) - O(n^{v_E-1})$.
As $\Phi$ is $(\oO,s)$-extendable, the claim follows.
Now by Lemma \ref{oct:focus} applied to each 
$(\Phi'{}^{j-1},\Phi'{}^j)$ successively, 
there is $\Psi' \in H^{\mc{O}^B(\Phi)}$ 
with $J' = J - \pl \Psi' \in H^{\Phi'_B}$.

Next we will reduce the support of $J'$ to $\Phi[L^r]_B$.
We define $\Phi^0=\Phi'$ and $\Phi^j=\Phi^{j-1}[L^j]$ 
for $j \in [r]$, and show that whp each
$(\Phi^{j-1},L^j)$ is $((2r)^{-jrs}\oO,s-3j)$-extendable. 
We show by induction on $j \in [r]$ that any
$\Phi^{j-1}$-extension $E=(H',F,\phi)$ 
of rank $s-3j$ is $(2r)^{-jrs}\oO$-dense in $(\Phi^{j-1},L^j)$.
Note that $\Phi^0=\Phi'$ is $(\oO-O(n^{-1}),s)$-extendable,
and the induction statement for $j$
implies that $\Phi^j$ is $((2r)^{-jrs}\oO,s-3j)$-extendable.
Thus for the induction step we can assume that $\Phi^{j-1}$ 
is $((2r)^{-(j-1)rs}\oO-O(n^{-1}),s-3j+3)$-extendable.

We can assume that $H' \sub B([s-3j+3] \sm [3])$.
For each $e \in H'_{B'}$ with $B' \in \tbinom{B}{j}$
and $e \sm F \ne \es$
we fix a $B(s-3j+3)$-embedding $\psi^e$ of $B(3)$ 
such that $\psi^e$ is the identity on $O^B$,
$e=\psi^e e'$ where $e' \in B(3)_{B'}$
with $e'(x)=(x,3)$ for all $x \in B'$, 
and $\psi^e$ is otherwise disjoint from $H'$, 
i.e.\ $\psi^e(V(B(3))) \cap V(H') = Im(e)$.

Let $E^+=(H^+,F^+,\phi^+)$, 
where $H^+ = H' \cup \bigcup \{ \psi^e e':
e \in H'_{B'}, B' \in \tbinom{B}{j}, e' \in B(3) \}$,
$F^+ = F \cup V(O^B)$ and $\phi^+$ restricts to
$\phi$ on $F$ and $\psi_0$ on $V(O^B)$. 
We claim that $\phi^+$ is a $\Phi^{j-1}$-embedding of $H^+[F^+]$.
To see this, consider any $f \in H^+[F^+]$ with $|Im(f)|=i<j$,
and write $f = f^1 \cup f^2$, 
where $Im(f^1) \sub F$ and $f^2 \in O^B$.
As $\phi$ is a $\Phi^{j-1}$-embedding of $H'[F]$,
we have $\phi f^1 = \psi^1 e^1\sS$ for some
$\psi^1 \in X_{E_0}(\Phi[\Phi^\tau])$, $\sS \in \Ss$, $e^1 \in B(3)$.
Then $\phi^+ f = \psi^1 (e^1 \cup f^2 \sS^{-1}) \sS$, 
so $\phi^+ f \in L^i$, which proves the claim. 
Thus $E^+$ is a $\Phi^{j-1}$-extension.

Let $X$ be the number of $\phi^* \in X_{E^+}(\Phi^{j-1})$ 
such that $\tau(\phi^*(i,x_i))=i$ for all 
$(i,x_i) \in V(H^+) \sm F^+$. Then 
$X_E(\Phi^{j-1},L^j) \ge X n^{v_E-v_{E^+}}$ by construction of $E^+$.
As $\Phi^{j-1}$ is $((2r)^{-(j-1)rs}\oO-O(n^{-1}),s-3j+3)$-extendable,
we have $X_{E^+}(\Phi^{j-1}) \ge (2r)^{-(j-1)rs}\oO n^{v_{E^+}}- O(n^{v_E-1})$, 
so $\mb{E}X_E(\Phi^{j-1},L^j) \ge r^{-v_{E^+}} (2r)^{-(j-1)rs}\oO n^{v_E}
- O(n^{v_E-1})$. Changing any value of $\tau$ affects 
$X_E(\Phi^{j-1},L^j)$ by $O(n^{v_E-1})$, so by Lemma \ref{lip3} 
whp $X_E(\Phi^{j-1},L^j) \ge (2r)^{-jrs} \oO n^{v_E}$.
This completes the induction step, so each
$(\Phi^{j-1},L^j)$ is $((2r)^{-jrs}\oO,s-3j)$-extendable. 

Now by Lemma \ref{oct:focus} repeatedly applied to 
each $(\Phi^{j-1},\Phi^j)$ successively, 
always with $s' \ge s-3r \ge 3(r-1)^2$ 
and extendability parameter 
$\oO' > (2r)^{-3r^4} (2r)^{r^5}n^{-0.6}
> (2(r-1))^{(r-1)^5}n^{-0.6}$,
there is $\Psi^0 \in H^{\mc{O}^B(\Phi)}$ 
with $J^0 = J' - \pl \Psi^0 \in H^{\Phi[L^r]_B}$.
We will define null $J^1,\dots,J^r \in H^{\Phi[L^r]_B}$
such that $J^j_\psi=0$ whenever $|Im(\psi) \cap V_0|<j$,
via the following construction of `reducing octahedra'.

Let $L^*$ be the set of $e \in L^r$
with $Im(e) \sm V_0 \ne \es$ such that $e=\psi^e_0 f^e_0$ 
for some $\psi^e_0\in X_{E_0}(\Phi[\Phi^\tau])$ 
and $f^e_0 \in B(3)_B$ (i.e.\ we can take $\sS=id$ 
in the definition of $L^r)$. For each $e \in L^*$
we fix some $\Phi$-embedding $\psi^e$ of $O^B$ of the form 
$\psi^e = \psi^e_0 \mid_{F^e} \pi^e$,
where $f^e_1 \in O^B$ and $\pi^e O^B$ is the copy
of $O^B$ in $B(3)$ spanned by $F^e := Im(f^e_0) \cup Im(f^e_1)$,
identified so that $f^e_0$ has sign $1$.
Note that for any $e' \in \psi^e O^B_B$ with $e' \ne e$
we have $|Im(e') \cap V_0| > |Im(e) \cap V_0|$ and $e' \in L^r$.
We define $J^1 = J^0 - \pl \Psi^1$, where
$\Psi^1 = \sum_{e \in L^*} J^0_e \{\psi^e\}$.
As $J^0$ and each $\chi(J^0_e,\psi^e)$ are symmetric,
we have $J^1_\psi=0$ whenever $Im(\psi) \cap V_0 = \es$.

Given $J^j$ with $0<j<r$, we define $J^{j+1}$ as follows.
Consider any $f \in \Phi^\circ_{r-j}$ disjoint from $V_0$,
write $B'=\tau(f) \sub B$, and suppose $f=Im(\psi^*)$,
where $\psi^* \in \Phi_{B'}$ with $\tau\psi^*=id$.
For any $\psi \in \Phi_B$ with $J^j_\psi \ne 0$
and $\psi^*\Ss \sub \psi\Ss$,
by definition of $L^r$ we can pick a representative
$\psi$ of $\psi\Ss$ with $\psi^* \sub \psi$ and $\tau\psi=id$.
Furthermore, for any $x \in B \sm B'$ and $\psi' \in \Phi_{B-x}$
with $\psi^* \sub \psi'$ and $\tau\psi'=id$,
if there is $\psi' \sub \psi$ with $J^j_\psi \ne 0$
then there are exactly two such $\psi$,
say $\psi^\pm$, obtained from each other by 
interchanging $\psi_0((x,0))$ and $\psi_0((x,1))$, where
as $J^j$ is null we have $J^j_{\psi^-}=-J^j_{\psi^+}$.
Thus there is $a_f \in H$ such that 
$J^j_\psi = \pm a_f$ whenever
$\psi^* \sub \psi$ and $\tau\psi=id$, where the sign 
is that of $\psi_0^{-1}\psi\mid_{B \sm B'}$ in $O^{B \sm B'}$. 
Fix $e$ with $\psi^* \sub e$, $\tau e = id$, $J^j_e=a_f$.
By symmetry, we have $J^j_\psi = \chi(a_f,\psi^e)_\psi$
whenever $J^j_\psi \ne 0$ with $\psi^*\Ss \sub \psi\Ss$.
We add $a_f \{ \psi^e \}$ to $\Psi^{j+1}$ for each such 
$f$, $e$ and let $J^{j+1} = J^j - \pl \Psi^{j+1}$.

We conclude with $J^r$ such that $J^r_e$ is zero
unless $e \in \psi_0 O^B_B \Ss$. As $J^r$ is symmetric and null,
we have $J^r = \chi(a,\psi_0)$ for some $a \in H$. 
Then $\Psi := \Psi' + a \{\psi_0\} + \sum_{j=0}^r \Psi^j$ 
has $\pl \Psi = J$. \qed

\medskip

Next we give two quantitative versions of Lemma \ref{oct:lattice}.
These will be used in the next subsection 
to prove two quantitative versions of
the main lemma of this section, which will in turn
both be used in the proof of Lemma \ref{bddint}
in the next section. We make the following definitions.

\begin{defn} \label{Guse} ($G$-use)
Suppose $H$ is a symmetric subgroup of $\GG^\Ss$
and $G$ is a symmetric generating set of $H$.
For $v \in \GG$ we write $|v|_G$
for the minimum possible $\sum_{g \in G} |c_g|$
where $v = \sum_{g \in G} c_g g$ with all $c_g \in \mb{Z}$.
For $\Psi \in H^{\mc{O}^B(\Phi)}$
we write $|\Psi|_G = \sum |\Psi_\phi|_G$.
If $J \in H^{\Phi_B}$ is symmetric
we write $|J|_G = \sum'_\psi |J_\psi|_G$,
where the sum is over any choice of
orbit representatives for $\Phi_B/\Ss$.
\end{defn}

The following lemma quantifies 
the total `$G$-use' of 
the octahedral decomposition $\Psi$
in terms of that of $J$.
We define $C(i) = 2^{(9i+2)^{i+5}}$.

\begin{lemma} \label{oct:lattice:total}
Let $\Phi$ be an $(\oO,s)$-extendable
$\Ss$-adapted $B$-complex with $|B|=r$, $s=3r^2$,
$n = |V(\Phi)| > n_0(r,\GG)$ large and $\oO>n^{-1/2}$.
Suppose $H$ is a symmetric subgroup of $\GG^\Ss$
and $G$ is a symmetric generating set of $H$.
Suppose $J \in H^{\Phi_B}$ is symmetric and null.
Then there is $\Psi \in H^{\mc{O}^B(\Phi)}$
with $\pl \Psi = J$ such that $|\Psi|_G \le C(r) |J|_G$.
\end{lemma} 

Following the proof of Lemma \ref{oct:lattice},
we quantify the total $G$-use in Lemma \ref{oct:focus}.

\begin{lemma} \label{oct:focus:total}
In Lemma \ref{oct:focus}, we can choose $\Psi$
with $|\Psi|_G \le C(r-r') |J|_G$ and 
$|J-\pl\Psi|_G \le 2^r C(r-r') |J|_G$.
\end{lemma}

The proof of Lemma \ref{oct:focus:total}
is the same as that of Lemma \ref{oct:focus},
noting also that when we apply the inductive hypothesis
each $|\Psi^{\psi^*}|_G \le C(r-r')|J^*|_G$, 
so $|\Psi|_G \le C(r-r')|J|_G$, and this
also gives $|J-\pl\Psi|_G \le 2^r C(r-r') |J|_G$.

\medskip

\nib{Proof of Lemma \ref{oct:lattice:total}.}
We will estimate the total $G$-use
during the proof of Lemma \ref{oct:lattice}.
We write $\Psi' = \sum_{j=1}^r \Psi'{}^j$, $J'{}^0 = J$ 
and $J'{}^j=J'{}^{j-1} - \pl \Psi'{}^j$ for $j>0$,
so $J'{}^r=J'$. By Lemma \ref{oct:focus:total}
each $|\Psi'{}^j|_G \le C(r-1) |J'{}^{j-1}|_G$
and $|J'{}^j|_G \le 2^r C(r-1) |J'{}^{j-1}|_G$,
so $|J'|_G \le 2^{r^2} C(r-1)^r |J|_G$.

Similarly, we write $\Psi^0 = \sum_{j=1}^r \Psi^{0,j}$, $J^{00}=J'$
and $J^{0,j}=J^{0,j-1} - \pl\Psi^{0,j}$ for $j>0$. Then $J^{0r}=J^0$, 
and each $\Phi^j=\Phi^{j-1}[L^j]$ is obtained by repeated
restriction to each $L^j_{B'}$ with $B' \in \tbinom{B}{j}$,
so writing $r_j = \tbinom{r}{j}$ we have
$|\Psi^{0,j}|_G \le C(r-j)^{r_j} |J^{0,j-1}|_G$,
and $|J^{0,j}|_G \le (2^r C(r-j))^{r_j} |J^{0,j-1}|_G$, so
$|J^0|_G \le 2^{r2^r} |J'|_G \prod_{i=0}^{r-1} C(i)^{r_i}$.

Next we have $|\Psi^1|_G \le |J^0|_G$
and $|J^1|_G \le 2^r |J^0|_G$.
For $j>0$ we have $|\Psi^{j+1}|_G \le |J^j|_G$
and $|J^{j+1}|_G \le 2^r |J^j|_G$, 
so $|\Psi^r|_G \le 2^{r^2} |J^0|_G 
\le 2^{r2^r+2r^2} |J|_G  C(r-1)^{2r} \prod_{i=0}^{r-2} C(i)^{r_i}$.
Recalling $C(i) = 2^{(9i+2)^{i+5}}$, we see that
$\Psi := \Psi' + a \{\psi_0\} + \sum_{j=0}^r \Psi^j$ 
has $|\Psi|_G < C(r) |J|_G$. \qed

\medskip

In our second quantitative version,
we suppose $\GG=\mb{Z}^D$ is free,
and consider rational decompositions,
where we now bound $G$-uses on every
function in $\Phi_r$ (as opposed to
the total bound in the previous version).

\begin{defn}
Suppose $H$ is a symmetric subgroup of $(\mb{Z}^D)^\Ss$
and $G$ is a symmetric generating set of $H$.
For $v \in \mb{Q} H$ we write $|v|_G$
for the minimum possible $\sum_{g \in G} |c_g|$
where $v = \sum_{g \in G} c_g g$ with all $c_g \in \mb{Q}$.
For $\Psi \in (\mb{Q} H)^{\mc{O}^B(\Phi)}$ 
and $\psi \in \Phi_B$ we write
$U_G(\Psi)_\psi = \sum \{ |\Psi_\phi|_G : 
 \psi=\phi\psi', \psi' \in O^B_B \}$.
\end{defn}

\begin{lemma} \label{oct:lattice:Q}
Let $\Phi$ be an $(\oO,s)$-extendable
$\Ss$-adapted $B$-complex 
where $n = |V(\Phi)| > n_0(r,\GG)$ is large, $|B|=r$,
$n^{-1/2} < \oO < \oO_0(r)$ and $s = 3r^2$.
Suppose $H$ is a symmetric subgroup of $(\mb{Z}^D)^\Ss$
and $G$ is a symmetric generating set of $H$.
Suppose $J \in (\mb{Q}H)^{\Phi_B}$ is symmetric and null 
with $|J_\psi|_G \le \tT$ for all $\psi \in \Phi_B$.
Then there is $\Psi \in (\mb{Q}H)^{\mc{O}^B(\Phi)}$
with $\pl \Psi = J$ such that
$U_G(\Psi)_\psi \le C(r,\oO) \tT$
for all $\psi \in \Phi_B$, where
$C(i,\oO) = 2^{C(i)} \oO^{-(9i)^{i+4}}$.
\end{lemma} 

Again we require the corresponding quantitative
version of Lemma \ref{oct:focus}.

\begin{lemma} \label{oct:focus:Q}
Let $\Phi$ be an $(\oO,s')$-extendable $\Ss$-adapted $B$-complex 
and $B' \sub B$, with $|B|=r$, $|B'|=r'$, $s'=3(r')^2$,
$n = |V(\Phi)| > n_0(r,\GG)$ large
and $n^{-1/2} < \oO < \oO_0(r)$.
Suppose $B' \in X \in \mc{P}^\Ss$
and $\Phi' \sub \Phi_X$ is 
such that $(\Phi,\Phi')$ is $(\oO,2)$-extendable
and $\Phi[\Phi']$ is $\Ss$-adapted.
Suppose $H$ is a symmetric subgroup of $(\mb{Z}^D)^\Ss$
and $G$ is a symmetric generating set of $H$.
Suppose $J \in (\mb{Q}H)^{\Phi_B}$ 
is symmetric and null with all $|J_\psi|_G \le \tT$.
Then there is $\Psi \in (\mb{Q}H)^{\mc{O}^B(\Phi)}$
with $J-\pl\Psi \in (\mb{Q}H)^{\Phi[\Phi']_B}$
such that for all $\psi \in \Phi_B$ 
both $U_G(\Psi)_\psi$ and $|(J-\pl\Psi)_\psi|_G$
are at most $2^r \oO^{-1} C' \tT$, where $C'=C(r-r',\oO)$.
\end{lemma}

\nib{Proof.}
We follow the proof of Lemma \ref{oct:focus}.
For each $\psi^* \in \Phi_{B'} \sm \Phi'$, 
we define $\Phi^*$, $J^*$, $\Psi^{\psi^*}$ as before,
where by the inductive hypothesis of Lemma \ref{oct:lattice:Q}
there is $\Psi^{\psi^*} \in (\mb{Q} H^*)^{\mc{O}^{B^*}(\Phi^*)}$
with $\pl \Psi^{\psi^*} = J^*$ and all
$U_{G^*}(\Psi^{\psi^*})_{\psi/\psi^*} \le C'\tT$, where 
$G^* = \{ h^* \in H^*: \pi(h^*) \in G \}$.
For each orbit representative $\psi^*$ 
and $\phi^* \in \mc{O}^{B^*}(\Phi^*)$
we add $\pi(\Psi^{\psi^*}_{\phi^*}) \mb{E}\{\phi\}$ to $\Psi$, where 
$\phi$ is uniformly random in $X_{E^{\psi^*}_{\phi^*}}(\Phi,\Phi')$.
Then $J-\pl\Psi \in (\mb{Q}H)^{\Phi[\Phi']_B}$
as in the proof of Lemma \ref{oct:focus}.

For any $\psi \in \Phi_B$ we have 
$U_G(\Psi)_\psi \le \sum_{\psi^* \in \Phi_{B'} \sm \Phi'}
\sum_{\phi^* \in \mc{O}^{B^*}(\Phi^*)} |\Psi^{\psi^*}_{\phi^*}|_{G^*}
\sum_{\psi^B \in O^B_B} \mb{P}(\psi=\phi\psi^B)$,
where $\phi$ is as above.
Note that any given $\psi^B$ 
can only contribute to $U_G(\Psi)_\psi$
if $\psi(x)=\phi^*\psi^B(x)$ for all $x \in B^*$
and $\psi(x)=\psi^*(x)$ for all $x \in B' \sm D$, 
where $D = D(\psi^B) = \{ x \in B': \psi^B(x)=(x,2) \}$.

For each $\psi' \in S_D := \{ \psi' \in \Phi_B
: \psi'\mid_{B \sm D} = \psi\mid_{B \sm D} \}$ 
and $\psi^*=\psi'\mid_{B'}$ we have 
$\mb{P}(\psi=\phi\psi^B) < \oO^{-1} n^{-|D|}$,
as $(\Phi,\Phi')$ is $(\oO,2)$-extendable.
We have $|S_D| < n^{|D|}$ and all
$|\Psi^{\psi^*}_{\phi^*}|_{G^*} \le C'\tT$,
so summing over $\psi^B$ and $\psi' \in S_{D(\psi^B)}$
gives $U_G(\Psi)_\psi \le 2^r \oO^{-1} C' \tT$.
The same bound applies to $|(J-\pl\Psi)_\psi|_G$. \qed

\medskip

\nib{Proof of Lemma \ref{oct:lattice:Q}.}
We follow the proof of Lemma \ref{oct:lattice},
estimating uses of any fixed $\psi \in \Phi_B$, and then average
over all choices of the initial $\Phi$-embedding $\psi_0$ of $O^B$. 
We use the same notation as in the proofs of
Lemmas \ref{oct:lattice} and \ref{oct:lattice:total}.

Letting $\tT'_0=\tT$, and $\tT'_j$ for $j>0$
be such that all $|J'{}^j_\psi|_G \le \tT'_j$,
by Lemma \ref{oct:focus:Q}, for each $j>0$
both $U_G(\Psi'{}^j)_\psi$ and
$|J'{}^j_\psi|_G$ are at most $2^{r-1} (\oO-O(n^{-1}))^{-1}
 C(r-1,\oO-O(n^{-1})) \tT'_{j-1}$,
so we can define 
$\tT'_j = 2^r \oO^{-1} C(r-1,\oO) \tT'_{j-1}$.
Then both $U_G(\Psi')_\psi$ and $|J'_\psi|_G$ are at most $2\tT'_r$.

Similarly, letting $\tT_0=2\tT'_r$ and $\tT_j$ for $j>0$
be such that all $|J^{0,j}_\psi|_G \le \tT_j$,
by Lemma \ref{oct:focus:Q}, for each $j>0$ both 
$U_G(\Psi^{0,j})_\psi$ and $|J^{0,j}_\psi|_G$ 
are at most $2\tT_{j-1}
[ 2^{r-j} ((2r)^{-jrs}\oO)^{-1} C(r-j,(2r)^{-jrs}\oO) ]^{r_j}$,
so we can define $\tT_j = \tT_{j-1}
[ 2^r ((2r)^{-jrs}\oO)^{-1} C(r-j,(2r)^{-jrs}\oO) ]^{r_j}$.
Then $U_G(\Psi^0)_\psi$ and $|J^0_\psi|_G$ 
are at most $\tT^* := 2\tT'_r (2^r \oO^{-1} (2r)^{sr^2})^{2^r}
\prod_{i=0}^{r-1} C(i, (2r)^{-sr^2} \oO)^{r_i}$.

Now write $B'=\{x: \psi(x) \notin V_0\}$ 
and $\psi'=\psi\mid_{B'}$. For $0 \le j < r-|B'|$ we have
$U_G(\Psi^{j+1})_{\psi'} \le 2^r |(J^j\mid_{\psi'})|_G$,
so $\sum_{j=1}^{r-|B'|} U_G(\Psi^j)_\psi
\le 2^{r^2} |(J^0\mid_{\psi'})|_G
\le 2^{r^2} \tT^* n^{r-|B'|}$.

Letting $\ov{\Psi}$ be the average of $\Psi$
(from the proof of Lemma \ref{oct:lattice})
over all $\Phi$-embeddings $\psi_0$ of $O^B$,
as $\Phi$ is $(\oO,2)$-extendable,
$U_G(\ov{\Psi})_\psi \le 2^{r^2 +1} \tT^* \sum_{B' \sub B} 
n^{2r-(r-|B'|)} (\oO n^{2r})^{-1} n^{r-|B'|}
< 2^{(r+1)^2} \oO^{-1} \tT^* < C(r,\oO) \tT$,
for $\oO < \oO_0(r)$, using $s=3r^2$,
$C(i,\oO) = 2^{C(i)} \oO^{-(9i)^{i+4}}$
and $C(i) = 2^{(9i+2)^{i+5}}$. \qed

\subsection{Lattices}

In this subsection we use the octahedral decompositions
from the previous subsection to characterise
the decomposition lattice $\bgen{\gG(\Phi)}$.
For the following lemma, we recall (Lemma \ref{LsubL-}) 
that $\bgen{\gG(\Phi)} \sub \mc{L}^-_\gG(\Phi)$.
We will show that if $\Phi$
is extendable then this inclusion becomes 
an equality when we restrict to null vectors.

\begin{lemma} \label{lattice:null}
Let $\Ss \le S_q$, $\mc{A}$ be a $\Ss^\le$-family
and $\gG \in (\mb{Z}^D)^{\mc{A}_r}$. 
Let $\Phi$ be a $\Ss$-adapted 
$(\oO,s)$-extendable $[q]$-complex
with $s=3r^2$, $n = |V(\Phi)| > n_0(q,D)$ large
and $\oO > n^{-1/2}$. Suppose $B \in C \in \mc{P}^\Ss_r$ and
$J \in \mc{L}^-_\gG(\Phi)  \cap (\mb{Z}^D)^{\Phi_C}$ is null.
Then $J \in \bgen{\gG(\Phi)}$.
\end{lemma}

\nib{Proof.}
Recall (Lemma \ref{gB}) that 
$f_B(J) \in (\gG^B)^{\Phi_B}$ is symmetric,
and (Lemma \ref{taut}) that $\gG^B$ is symmetric.
We also note that $f_B(J)$ is null,
as for any $\psi' \in \Phi_{B'}$
with $B' \sub B$, $|B'|=r-1$ and $\sS \in \Ss^B$
we have $(\pl f_B(J)_{\psi'})_\sS
= \sum \{ (f_B(J)_\psi)_\sS: \psi' \sub \psi \}
= \sum \{ J_{\psi\sS}: \psi' \sub \psi \}
= \pl J_{\psi'\sS} = 0$.
By Lemma \ref{oct:lattice} we have
$\Psi \in (\gG^B)^{\mc{O}^B(\Phi)}$ with
$f_B(J) = \pl \Psi = \sum_{\psi^*} \chi(\Psi_{\psi^*},\psi^*)$.
By definition of $\gG^B$ we can fix 
integers $m^\tT_{\psi^*}$ so that each 
$\Psi_{\psi^*} = \sum_{\tT \in \mc{A}_B} m^\tT_{\psi^*} \gG^\tT$.

We will show for any $\psi^* \in \mc{O}^B(\Phi)$ and $\tT \in A_B$
that there is $\Psi^{\psi^* \tT} \in \mb{Z}^{A(\Phi)}$ with
$f_B(\pl^\gG \Psi^{\psi^* \tT}) = \chi(\gG^\tT,\psi^*)$.
This will suffice to prove the lemma. Indeed, letting
$\Psi^* = \sum_{\psi^*,\tT} m^\tT_{\psi^*} \Psi^{\psi^* \tT}$,
we will have $\pl^\gG \Psi^* 
 = \sum_{\psi^*,\tT} m^\tT_{\psi^*} \chi(\gG^\tT,\psi^*)
 = \sum_{\psi^*} \chi(\Psi_{\psi^*},\psi^*) =  \pl \Psi = J$.

Consider $\psi^* \in \mc{O}^B(\Phi)$, $A \in \mc{A}$, $\tT \in A_B$,
let $F = \tT(B) \times [2]$ and define $\psi_0:F \to V(\Phi)$ by
$\psi_0(\tT(i),x)=\psi^*(i,x)$ for $i \in B$, $x \in [2]$.
We claim that $\psi_0 \in \mc{O}^{\tT(B)}(\Phi)$,
and so $E=([q](2),F,\psi_0)$ is a $\Phi$-extension.
To see this, note that for any $\ups \in O^{\tT(B)}_{\tT(B)}$,
defining $\ups' \in O^B_B$ such that
$\ups'(i)=(i,x)$ when $\ups(\tT(i))=(\tT(i),x)$,
as $\psi^* \in \mc{O}^B(\Phi)$ we have 
$\psi_0 \ups \tT = \psi^* \ups' \in \Phi_B$.
As $\Phi$ is $\Ss$-adapted, we deduce 
$\psi_0 \ups \in \Phi_{\tT(B)}$, which proves the claim. 

As $\Phi$ is $(\oO,s)$-extendable,
we can choose $\psi^+ \in X_E(\Phi)$.
We let $\Psi^{\psi^* \tT}$ 
be the sum of all $\pm \{ \psi^+ \circ \phi \}$
over all $A(2)$-embeddings $\phi$ of $A$ 
such that $\phi(i)=(i,1)$ when $i \notin \tT(B)$,
where the sign is that of $\phi\mid_{\tT(B)}$ in $O^{\tT(B)}$,
i.e.\ that of $\psi' \in O^B$ defined
by $\psi' = (\psi^*)^{-1} \psi$,
where $\psi = \psi_0\phi\tT$ with $\psi_0$ as above.
Then all entries in $\pl^\gG \Psi^{\psi^* \tT}$
cancel in $\pm$ pairs, except that for each
$\pm \{ \psi^+ \circ \phi \}$ in $\Psi^{\psi^* \tT}$
there is a non-cancelling term
$\gG(\psi^+ \phi)^{\psi\Ss} = \gG[\psi]^\tT$
(by Lemma \ref{atom=}.i)
as $\psi = \psi_0\phi\tT = \psi^+\phi\tT$.
Now $f_B(\gG[\psi]^\tT)_\psi = \pm \gG^\tT 
= \chi(\gG^\tT,\psi^*)_\psi$,
where the sign is that of
$\psi' = (\psi^*)^{-1} \psi \in O^B$, so by symmetry
$f_B(\pl^\gG \Psi^{\psi^* \tT}) = \chi(\gG^\tT,\psi^*)$,
as required to prove the lemma. \qed

\medskip

To motivate the characterisation of decomposition lattices
in general, it may be helpful to consider 
the decomposition lattice of triangles in a complete tripartite
graph (a very special case of Theorem \ref{Hdecomp:partite1}).
Say $T$ is a complete tripartite graph with parts $(A,B,C)$.
It is not hard to show that $J \in \mb{Z}^T$ is in the decomposition 
lattice iff $J$ is `balanced', in that each $a \in A$
has $\sum_{b \in B} J_{ab} = \sum_{c \in C} J_{ac}$
and similarly for each $b \in B$ and $c \in C$.
At first sight this seems a rather different condition to the 
tridivisibility condition that arises for nonpartite
triangle decomposition (even degrees and $3 \mid \sum J$).
However, we can unify the conditions by lifting $J$ to a vector
$J^+ \in (\mb{Z}^3)^T$ in which we assign different basis vectors to
the three bipartite subgraphs, say $(1,0,0)$ to $(B,C)$,
$(0,1,0)$ to $(A,C)$ and $(0,0,1)$ to $(A,B)$. We want to characterise
when $J^+$ is in the lattice generated by all vectors 
$v(abc) \in (\mb{Z}^3)^T$, where for $a \in A$, $b \in B$, $c \in C$ 
we let $v(abc)_{bc}=(1,0,0)$, $v(abc)_{ac}=(0,1,0)$, $v(abc)_{ab}=(0,0,1)$
and $v(abc)_{xy}=0$ otherwise. 
The lifted vertex degree condition is that for any $a \in A$ 
we have $\sum_{x \in V(T)} J^+_{ax}$ in the lattice generated 
by $(0,1,0)$ and $(0,0,1)$, and similarly 
for any $b \in B$ and $c \in C$; 
this is equivalent to $J$ being balanced.
This example suggests the form of the degree
conditions in the following definition.

\begin{defn} \label{def:L}
Let $\Phi$ be a $[q]$-complex,
$\mc{A}$ be a $[q]$-complex family
and $\gG \in (\mb{Z}^D)^{\mc{A}_r}$. 
For $J \in (\mb{Z}^D)^{\Phi_r}$ we define
$J^\sharp \in ((\mb{Z}^D)^Q)^\Phi$ by 
$(J^\sharp_{\psi'})_B = 
\sum \{ J_\psi: \psi' \sub \psi \in \Phi_B \}$
for $B \in Q$, $\psi' \in \Phi$.
Similarly, we define
$\gG^\sharp \in ((\mb{Z}^D)^Q)^{\cup \mc{A}}$
by $(\gG^\sharp_{\tT'})_B = 
\sum \{ \gG_\tT: \tT' \sub \tT \in A_B \}$
for $B \in Q$, $\tT' \in A \in \mc{A}$.
We let $\mc{L}_\gG(\Phi)$ be the set of all
$J \in (\mb{Z}^D)^{\Phi_r}$ such that
$(J^\sharp)^O \in \sgen{\gG^\sharp[O]}$
for any $O \in \Phi/\Ss$. 
\end{defn}

Note that $\gG^\sharp = (\gG^\sharp_i: 0 \le i \le r)$
where each $\gG^\sharp_i$ is a vector system for $\mc{A}_i$.
For convenient use in Lemma \ref{lattice} 
we will reformulate Definition \ref{def:L}
in terms of iterated independent shadows,
in the sense of the following definition.

\begin{defn} \label{def:ishadow}
(independent shadows) 

For each $B \in [q]_i$ and any set $S$
we define $\pi_B: \mb{Z}^S \to \mb{Z}^{[q]_i \times S}$
by $\pi_B(v) = e_B \otimes v$,
i.e.\ $\pi_B(v)$ is $v$ in the block of $S$
coordinates belonging to $B$ and $0$ otherwise.

Let $\Phi$ be a $[q]$-complex.
We define $\pi_i: (\mb{Z}^S)^{\Phi_i} 
\to (\mb{Z}^{[q]_i \times S})^{\Phi_i}$
by $\pi_i(J)_\psi = \pi_B(J_\psi)$ for 
$\psi \in \Phi_B$, $B \in [q]_i$.
We also write $J^+ = \pi_i(J)$ for $J \in (\mb{Z}^S)^{\Phi_i}$.

For $i=r,\dots,0$ we define $D_i$ by $D_r=[q]_r \times [D]$ 
and $D_i = [q]_i \times D_{i+1}$ for $0 \le i < r$.
For $J \in (\mb{Z}^{D_{i+1}})^{\Phi_{i+1}}$
we define $\pl^* J \in (\mb{Z}^{D_i})^{\Phi_i}$ 
by $\pl^* J = \pi_i(\pl_i J) = (\pl J)^+$.

For $J \in (\mb{Z}^D)^{\Phi_r}$ and $0 \le i \le r$
we define $\pl^*_i J = (\pl^*)^{r-i} (J^+)$.

Similarly, given a $[q]$-complex family $\mc{A}$
and $\gG \in (\mb{Z}^D)^{\mc{A}_r}$
we define $\pl^*_i \gG = (\pl^*)^{r-i} (\pi_r(\gG))$,
where for $\gG' \in (\mb{Z}^{D_{i+1}})^{\mc{A}_{i+1}}$
and $\tT' \in A \in \mc{A}$
we define $(\pl \gG')_{\tT'} = 
 \sum \{ \gG'_\tT: \tT \in A\mid_{\tT'} \}$
and $(\pl^* \gG')_{\tT'} = \pi_i((\pl \gG')_{\tT'})$.

We let $\mc{L}^i_\gG(\Phi)$ be the set of all
$J \in (\mb{Z}^D)^{\Phi_r}$ such that
$\pl^*_i J \in \mc{L}^-_{\pl^*_i \gG}(\Phi)$,
i.e.\ $(\pl^*_i J)^O \in \sgen{(\pl^*_i \gG)[O]}$
for any $O \in \Phi_i/\Ss$. 
\end{defn}

\begin{rem} \label{rem:unravel}
Unravelling the iterative definitions
in Definition \ref{def:ishadow}, we see that
each $D_i = [D] \times \prod_{j=i}^r [q]_j$, and 
for each $B^i \in [q]_i$, $\psi^i \in \Phi_{B^i}$,
we have $(\pl^*_i J)_{\psi^i}$ supported 
on `full chains from $B^i$', i.e.\
if $((\pl^*_i J)_{\psi^i})_C \ne 0$ 
then $C=(B^i,\dots,B^r)$ with each $B^j \in [q]_j$ 
and $B^j \sub B^{j+1}$ for $i \le j < r$.

We obtain $((\pl^*_i J)_{\psi^i})_C \in \mb{Z}^D$
by summing $J_{\psi^r}$
over all choices of $(\psi^i,\dots,\psi^r)$
with each $\psi^j \in \Phi_{B^j}$ and 
$\psi^j \sub \psi^{j+1}$ for $i \le j < r$.

Similarly, for $A \in \mc{A}$, $\tT^i \in A_{B^i}$
we obtain $((\pl^*_i \gG)_{\tT^i})_C \in \mb{Z}^D$
by summing $\gG_{\tT^r}$
over all choices of $(\tT^i,\dots,\tT^r)$
with each $\tT^j \in A_{B^j}$ and 
$\tT^j \sub \tT^{j+1}$ for $i \le j < r$.
\end{rem}

Now we reformulate Definition \ref{def:L}
using Definition \ref{def:ishadow}.

\begin{lemma} \label{lem:L} 
$\mc{L}_\gG(\Phi) = \cap_{i=0}^r \mc{L}^i_\gG(\Phi)$.
\end{lemma}

\nib{Proof.}
Fix $\psi^O \in O \in \Phi_i/\Ss$, 
say with $\psi^O \in \Phi_{B'}$.
We need to show $(J^\sharp)^O \in \sgen{\gG^\sharp[O]}$
iff $(\pl^*_i J)^O \in \sgen{(\pl^*_i \gG)[O]}$.

We have $(J^\sharp)^O \in \sgen{\gG^\sharp[O]}$ iff
there is $n \in \mb{Z}^{\mc{A}_{B'}}$ with 
$(J^\sharp)^O = \sum_\sS \{ n_\sS \gG^\sharp(\psi^O \sS^{-1}) \}$,
i.e.\ all $J^\sharp_{\psi^O \sS'} = 
\sum_\sS \{ n_\sS \gG^\sharp_{\sS\sS'} \}$, 
i.e.\ for each $A \in \mc{A}$, $\sS' \in A[B']$, $B \in Q$ we have
$\sum \{ J_\psi: \psi^O \sS' \sub \psi \in \Phi_B \} =
\sum_{\sS,\tT} \{ n_\sS \gG_\tT: \sS\sS' \sub \tT \in A_B \}$.

We have $(\pl^*_i J)^O \in \sgen{(\pl^*_i \gG)[O]}$
iff there is $n' \in \mb{Z}^{\mc{A}_{B'}}$
such that for each $A \in \mc{A}$, $\sS' \in A[B']$ and chain 
$C = (B^i,\dots,B^r)$ from $B^i=B'$ to $B^r=B \in Q$,
we have $\sum_{\psi^i,\dots,\psi^r} J_{\psi^r} 
= \sum_{\sS,\tT^i,\dots,\tT^r} n'_\sS \gG_{\tT^r}$,
where we sum over $(\tT^i,\dots,\tT^r)$
and $(\psi^i,\dots,\psi^r)$
with $\tT^i=\sS\sS'$, $\psi^i=\psi^O \sS'$,
each $\tT^j \in A_{B^j}$, $\psi^j \in \Phi_{B^j}$,
and each $\tT^j \sub \tT^{j+1}$, $\psi^j \sub \psi^{j+1}$. 

Setting $n' = n$ we see that the two conditions
are identical, as any such chain $C$ and $\sS',\sS$ 
uniquely specifies all $\psi^j = \psi^r\mid_{B^j}$
and $\tT^j = \tT^r\mid_{B^j}$. \qed
 
\medskip

Our next lemma shows that the two definitions
of $\pl^*_i$ (for vectors and for vector systems)
are compatible with each other.

\begin{lemma} \label{plpl}
If $\Psi \in \mb{Z}^{\mc{A}(\Phi)}$ then
$\pl^{\pl^*_i \gG} \Psi = \pl^*_i (\pl^\gG \Psi)$.  
\end{lemma}

\nib{Proof.}
By linearity we can assume $\Psi=\{\phi\}$
for some $\phi \in A(\Phi)$. We prove 
the identity by induction on $i=r,\dots,0$.
In the base case $i=r$ we have
$\pl^{\pl^*_r \gG} \Psi
= (\pl^*_r \gG)(\phi)
= (\pi_r \gG)(\phi)
= \pi_r(\gG(\phi))
= \pl^*_r (\pl^\gG \Psi)$, where
in the third equality we used
$(\pi_r \gG)(\phi)_{\phi\tT}
= (\pi_r \gG)_\tT
= \pi_r(\gG_\tT)
= \pi_r(\gG(\phi)_{\phi\tT})$.
For the induction step with $i<r$ we have
$\pl^{\pl^*_i \gG} \Psi
= (\pl^*_i \gG)(\phi)
= (\pl^* \pl^*_{i+1} \gG)(\phi)
= \pl^* (\pl^*_{i+1} \gG)(\phi)
= \pl^* \pl^{\pl^*_{i+1} \gG} \Psi
= \pl^* \pl^*_{i+1} (\pl^\gG \Psi)
=  \pl^*_i (\pl^\gG \Psi)$,
where in the third equality,
writing $\gG' = \pl^*_{i+1} \gG$,
for $\tT' \in A \in \mc{A}$ we used 
$(\pl^* \gG')(\phi)_{\phi\tT'} = (\pl^* \gG')_\tT 
= \pi_i \sum \{ \gG'_\tT: \tT \in A\mid_{\tT'} \}
= \pi_i \sum \{ \gG'_\tT: \phi\tT \in \Phi\mid_{\phi\tT'} \}
= \pl^* ( \gG'(\phi)_{\phi\tT'} )$. \qed

\medskip

Now we come to the main lemma of this section,
which characterises the decomposition lattice.

\begin{lemma} \label{lattice}
Let $\Ss \le S_q$, $\mc{A}$ be a $\Ss^\le$-family
and $\gG \in (\mb{Z}^D)^{\mc{A}_r}$. 
Let $\Phi$ be a $\Ss$-adapted 
$(\oO,s)$-extendable $[q]$-complex with $s=3r^2$, 
$n = |V(\Phi)| > n_0(q,D)$ large and $\oO > n^{-1/2}$.
Then $\bgen{\gG(\Phi)} = \mc{L}_\gG(\Phi)$.
\end{lemma}

\nib{Proof.} 
Note for any molecule $\gG(\phi) \in \gG(\Phi)$ that 
$(\pl^*_i \gG)(\phi) \in (\pl^*_i \gG)(\Phi)
\sub \mc{L}^-_{\pl^*_i \gG}(\Phi)$,
so $\bgen{\gG(\Phi)} \sub \mc{L}^i_\gG(\Phi)$.
We now show that the reverse inclusions hold.
Suppose $J \in \mc{L}_\gG(\Phi)$.
We will define $\Psi^0,\dots,\Psi^r \in \mb{Z}^{\mc{A}(\Phi)}$ so that,
letting $J^0=J-\pl^\gG\Psi^0$ and $J^i=J^{i-1}-\pl^\gG\Psi^i$ for $i \in [r]$,
each $\pl^*_i J^i = 0$. We will then have $J^r=0$, so 
$J = \sum_{i=0}^r \pl^\gG\Psi^i \in \bgen{\gG(\Phi)}$, as required.

We start by noting that\footnote{ 
Recall that $\es$ may denote the empty set
or the function with empty domain.
We write $\es^A$ for the copy of $\es$ in $A$.} 
$\pl^*_0 J_\es \in (\pl^*_0 \gG)_\es$ as $J \in \mc{L}_0(\Phi)$,
so we have integers $k_A$ with 
$\pl^*_0 J_\es = \sum_{A \in \mc{A}} k_A (\pl^*_0 \gG)_{\es^A}$.
We can take $\Psi^0 = \sum_{A \in \mc{A}} k_A \{\phi^A\}$ 
for any choices of $\phi^A \in A(\Phi)$. 
Then $J^0=J-\pl^\gG\Psi^0$ has $\pl^*_0 J^0 = 0$.
It remains to define $\Psi^i$ given $J^{i-1}$ for some $i \in [r]$.

Note that $J^{i-1} \in \mc{L}^i_\gG(\Phi)$ as $J \in \mc{L}^i_\gG(\Phi)$
and $\gG(\Phi) \sub \mc{L}^i_\gG(\Phi)$, so
$\pl^*_i J^{i-1} \in \mc{L}^-_{\pl^*_i \gG}(\Phi)$.
We write $\pl^*_i J^{i-1} = \sum_{C \in \mc{P}^\Ss_i} J^C$,
where each $J^C$ is uniquely defined
by $J^C_\psi = \pl^*_i J^{i-1}_\psi$ for $\psi \in \Phi_C$.

We claim that each $J^C$ is null. 
Indeed, for any $\psi \in \Phi_{i-1}$  
we have $0 = \pl^*_{i-1} J^{i-1}_\psi
= ( \pl \pl^*_i J^{i-1} )^+_\psi
= \sum_C (\pl J^C)^+_\psi$,
so by linear independence 
each $(\pl J^C)_\psi=0$, as claimed.

By Lemma \ref{lattice:null}
each $J^C \in \bgen{(\pl^*_i \gG)(\Phi)}$,
so we have $\Psi^i \in \mb{Z}^{\mc{A}(\Phi)}$
with $\pl^{\pl^*_i \gG} \Psi^i = \pl^*_i J^{i-1}$.
Letting $J^i = J^{i-1}-\pl^\gG\Psi^i$,
by Lemma \ref{plpl} we have
$\pl^*_i J^i = 0$, as required. \qed

\medskip

To illustrate the use of Lemma \ref{lattice},
we give the following characterisation of the
decomposition lattice for nonpartite hypergraph decomposition,
thus giving an independent proof of (a generalisation of)
a result of Wilson \cite{W5} (a similar generalisation 
is implicit in \cite{GKLO2}). 

\begin{lemma} \label{Hlattice:nonpartite}
Let $H$ be an $r$-graph on $[q]$ and $\Phi$ be 
an $(\oO,s)$-extendable $S_q$-adapted $[q]$-complex
where $n=|V(\Phi)|>n_0(q)$ is large, $s=3r^2$ and $\oO>n^{-1/2}$.
Let $H(\Phi) = \{ \phi(H): \phi \in \Phi_q\}$.
Suppose $G \in \mb{N}^{\Phi^\circ_r}$.
Then $G \in \bgen{H(\Phi)}$ iff $G$ is $H$-divisible.
\end{lemma}

\nib{Proof.}
As in example $i$ in subsection \ref{sec:vvd},
we have $G \in \bgen{H(\Phi)}$ iff $G^* \in \bgen{\gG(\Phi)}$,
with $G^* \in \mb{N}^{\Phi_r}$ defined by
$G^*_\psi = G_{Im(\psi)}$ for $\psi \in \Phi_r$,
and $\gG \in \{0,1\}^{A_r}$ with $A=S_q^\le$
and $\gG_\tT = 1_{Im(\tT) \in H}$.
By Lemma \ref{lattice} we have
$\bgen{\gG(\Phi)} = \mc{L}_\gG(\Phi)$.
By Definition \ref{def:L} we need to show
that  $G$ is $H$-divisible iff
$((G^*)^\sharp)^O \in \sgen{\gG^\sharp[O]}$
for any $O \in \Phi/S_q$. 

Fix any $O \in \Phi/S_q$,
write $e=Im(O) \in \Phi^\circ$ and $i=|e|$.
Then $((G^*)^\sharp)^O \in (\mb{Z}^Q)^O = \mb{Z}^{Q \times O}$
is a vector supported on the coordinates
$(B,\psi')$ with $B' \sub B \in Q$ and $\psi' \in O \cap \Phi_{B'}$
in which every nonzero coordinate is equal:
we have $((G^*)^\sharp_{\psi'})_B) = 
\sum \{ G^*_\psi: \psi' \sub \psi \in \Phi_B \} = (r-i)! |G(e)|$.
Similarly, $\sgen{\gG^\sharp[O]}$ is generated by vectors
with the same support that are constant on the support, 
where the constant can be $(r-i)! |H(f)|$ for any $f \in [q]_i$,
and so any multiple of $(r-i)! gcd_i(H)$.
Therefore $((G^*)^\sharp)^O \in \sgen{\gG^\sharp[O]}$
iff $gcd_i(H)$ divides $|G(e)|$, as required. \qed

\medskip

For the proof of Lemma \ref{bddint} in the next section,
we also require two quantitative versions of Lemma \ref{lattice},
analogous to those given above for Lemma \ref{oct:lattice}.
We recall the notation for $G$-use from Definition \ref{Guse},
and for $B \in Q$ fix the symmetric generating set
$G = G^B = \{ \gG^\tT: \tT \in \mc{A}_B \}$ for $\gG^B$.
Then we have the following relationship between $G$-use 
and use in the sense of Definition \ref{def:use+bdd}.

\begin{lemma} \label{U=fB}
Suppose $B \in C \in \mc{P}^\Ss_r$ and
$J \in \mc{L}^-_\gG(\Phi)  \cap (\mb{Z}^D)^{\Phi_C}$.
Then $U(J^O) = |f_B(J^O)|_G$ for each $O \in \Phi_r/\Ss$,
so $U(J)_\es = |f_B(J)|_G$.
\end{lemma}

\nib{Proof.}
We fix an orbit representative $\psi \in O$ 
and write $U(J^O)$ as the minimum possible value of
$\sum \{ |x^\tT_\psi| \}$ over all expressions of
$J^O = \sum \{ x^\tT_\psi \gG[\psi]^\tT \}$ 
as a $\mb{Z}$-linear combination of $\gG$-atoms at $O$.
We can write any such expression using
some fixed representative $\psi \in \Phi_B$;
then $f_B(J^O)_\psi = \sum_\tT x^\tT_\psi \gG^\tT$.
As $|f_B(J^O)|_G$ is the minimum value of
$\sum \{ |x^\tT_\psi| \}$ over all such expressions 
the lemma follows. \qed

\medskip

Our first quantitative version of Lemma \ref{lattice} 
will bound the total use $U(\Psi)$ of atoms by $\Psi$
in terms of that by $J$, where
$U(\Psi) = \sum_\phi |\Psi_\phi| U(\gG(\phi))
\le q^{2r} |\Psi|$, where $|\Psi| = \sum_\phi |\Psi_\phi|$.
We start by stating
the analogous statement for Lemma \ref{lattice:null}.

\begin{lemma} \label{lattice:null:total}
Let $\Ss \le S_q$, $\mc{A}$ be a $\Ss^\le$-family
and $\gG \in (\mb{Z}^D)^{\mc{A}_r}$. 
Let $\Phi$ be a $\Ss$-adapted 
$(\oO,s)$-extendable $[q]$-complex with $s=3r^2$,
$n = |V(\Phi)| > n_0(q,D)$ large and $\oO > n^{-1/2}$.
Suppose $B \in C \in \mc{P}^\Ss_r$ and
$J \in \mc{L}^-_\gG(\Phi)  \cap (\mb{Z}^D)^{\Phi_C}$ is null.
Then there is $\Psi^* \in \mb{Z}^{\mc{A}(\Phi)}$ with
$\pl^\gG \Psi^* = J$ and $U(\Psi^*) \le 2^r q^{2r} C(r) U(J)$.
\end{lemma}

The proof of Lemma \ref{lattice:null:total}
is the same as that of Lemma \ref{lattice:null}.
When we apply Lemma \ref{oct:lattice:total}
to $f_B(J) \in (\gG^B)^{\Phi_B}$ we obtain
$\Psi \in (\gG^B)^{\mc{O}^B(\Phi)}$ with $f_B(J) = \pl \Psi$ 
and $|\Psi|_G \le C(r) |f_B(J)|_G$. We write each 
$\Psi_{\psi^*} = \sum_{\tT \in \mc{A}_B} m^\tT_{\psi^*} \gG^\tT $
with $\sum_{\tT \in \mc{A}_B} |m^\tT_{\psi^*}| = |\Psi_{\psi^*}|_G$.
Defining $\Psi^*$ as in the proof of Lemma \ref{lattice:null},
we have $|\Psi^*| \le 2^r |\Psi|_G \le 2^r C(r) |f_B(J)|_G$,
so $U(\Psi^*) \le 2^r q^{2r} C(r) U(J)$ by Lemma \ref{U=fB}.

\begin{lemma} \label{lattice:total}
Let $\Ss \le S_q$, $\mc{A}$ be a $\Ss^\le$-family
with $\gG \in (\mb{Z}^D)^{\mc{A}_r}$. 
Let $\Phi$ be a $\Ss$-adapted 
$(\oO,s)$-extendable $[q]$-complex with $s=3r^2$,
$n = |V(\Phi)| > n_0(q,D)$ large and $\oO > n^{-1/2}$.
Suppose $J \in \mc{L}_\gG(\Phi)$.
Then there is $\Psi \in \mb{Z}^{\mc{A}(\Phi)}$
with $\pl^\gG\Psi = J$ and
$U(\Psi) \le 2^{(9q)^{q+2}} U(J)$.
\end{lemma}

\nib{Proof.} 
We follow the proof of Lemma \ref{lattice}.
We can choose the $k_A$ so that 
$U(\Psi^0) = \sum |k_A| \le U(J)$ 
and $U(J^0) \le q^{2r} U(J)$.
For each $i \in [r]$, by Lemma \ref{lattice:null:total}
we can take $U(\Psi^i) \le 2^i q^{2i} C(i) U(\pl^*_i J^{i-1})
  \le (2rq^2)^i C(i) U(J^{i-1})$, and so 
$U(J^i) \le q^{2r} (2r)^i C(i) U(J^{i-1})$,
where $U(\pl^*_i J^{i-1})$ is defined
with respect to $(\pl^*_i \gG)$-atoms.
Then $U(\Psi) \le 2U(J^0) \prod_i q^{2r} 2^{r+i} C(i)
\le 2^{(9q)^{q+2}} U(J)$, as $C(i) = 2^{(9i+2)^{i+5}}$. \qed

\medskip

Our second quantitative version of Lemma \ref{lattice} 
is analogous to Lemma \ref{oct:lattice:Q}:
we seek a rational decomposition $\Psi$ for which we bound
the usage of every orbit in terms of that in $J$. We start with
the analogous statement for Lemma \ref{lattice:null}.

\begin{lemma} \label{lattice:null:Q}
Let $\Ss \le S_q$, $\mc{A}$ be a $\Ss^\le$-family
and $\gG \in (\mb{Z}^D)^{\mc{A}_r}$. 
Let $\Phi$ be a $\Ss$-adapted 
$(\oO,s)$-extendable $[q]$-complex with $s=3r^2$,
$n = |V(\Phi)| > n_0(q,D)$ large and $\oO > n^{-1/2}$.
Suppose $C \in \mc{P}^\Ss_r$ and
$J \in \mb{Q}\mc{L}^-_\gG(\Phi) 
\cap (\mb{Q}^D)^{\Phi_C}$ is null with
$U(J)_\psi \le \eps$ for all $\psi \in \Phi_B$.
Then there is $\Psi^* \in \mb{Q}^{\mc{A}(\Phi)}$ 
with $\pl^\gG \Psi^* = J$ and 
$U(\Psi^*)_O \le q^{2r} C(r,\oO) \oO^{-1} \eps$ 
for each orbit $O \in \Phi_r/\Ss$.
\end{lemma}

\nib{Proof.}
We follow the proof of Lemma \ref{lattice:null}.
Applying Lemma \ref{oct:lattice:Q}
to $f_B(J) \in (\mb{Q}\gG^B)^{\Phi_B}$ gives
$\Psi \in (\mb{Q}\gG^B)^{\mc{O}^B(\Phi)}$ with $f_B(J) = \pl \Psi$
and $U_G(\Psi)_\psi \le C(r,\oO) \eps$
for all $\psi \in \Phi_B$. We write each 
$\Psi_{\psi^*} = \sum_{\tT \in \mc{A}_B} m^\tT_{\psi^*} \gG^\tT$
with $\sum_{\tT \in \mc{A}_B} |m^\tT_{\psi^*}| = |\Psi_{\psi^*}|_G$.

We let $\Psi^* = \sum_{\psi^*,\tT} m^\tT_{\psi^*} \Psi^{\psi^* \tT}$,
where we modify the definition of each $\Psi^{\psi^* \tT}$
by averaging over the choice of $\psi^+ \in X_E(\Phi)$.
To estimate $U(\Psi^*)_O$ for some $O \in \Phi_r/\Ss$
we fix any representative $\psi' \in O$,
say with $\psi' \in \Phi_{B'}$. For each 
$A \in \mc{A}$, $\tT \in A_B$, writing
$r' = |\tT(B) \sm B'|$, there are at most $n^{r'}$
choices of $\psi \in \Phi_B$ such that $\psi\tT^{-1}$ 
agrees with $\psi'$ on $\tT(B) \cap B'$
and each has $U_G(\Psi)_\psi \le C(r,\oO) \eps$.
For each $\psi^* \in \mc{O}^B(\Phi)$ with 
$\psi = \psi^* \ups$ for some $\ups \in O^B_B$,
letting $\phi$ be the $A(2)$-embedding of $A$
such that $\phi(i)=(i,1)$ for $i \notin \tT(B)$
and $\phi(i)=(i,x)$ when $i=\tT(j)$
with $j \in B$ and $\ups(j)=(j,x)$,
for uniformly random $\psi^+ \in X_E(\Phi)$ we have
$\mb{P}(\psi' \sub \psi^+ \phi) < \oO^{-1} n^{-r'}$,
as $\Phi$ is $(\oO,s)$-extendable.
Summing over $\psi'$ and $\tT$ gives
$U(\Psi^*)_O \le q^{2r} C(r,\oO) \oO^{-1} \eps$. \qed

\medskip

We conclude this section by proving the second quantitative
version of Lemma \ref{lattice}, which can be viewed
as a rational version of Lemma \ref{bddint}, and will
form the basis of the `randomised rounding' aspect
of the proof referred to in the introduction.

\begin{lemma} \label{lattice:Q}
Let $\Ss \le S_q$, $\mc{A}$ be a $\Ss^\le$-family
with $\gG \in (\mb{Z}^D)^{\mc{A}_r}$. 
Let $\Phi$ be a $\Ss$-adapted 
$(\oO,s)$-extendable $[q]$-complex with $s=3r^2$,
$n = |V(\Phi)| > n_0(q,D)$ large 
and $n^{-1/2} < \oO < \oO_0(r)$.
Suppose $J \in \mb{Q}\mc{L}_\gG(\Phi)$
with $U(J)_O \le \eps$ for all $O \in \Phi_r/\Ss$.
Then there is $\Psi \in \mb{Q}^{\mc{A}(\Phi)}$
with $\pl^\gG\Psi = J$ and all
$U(\Psi)_O \le C(q,\oO) \eps$.
\end{lemma}

\nib{Proof.}
We follow the proof of Lemma \ref{lattice}.
We can choose the $k_A$ so that 
$\sum_A |k_A| \le U(J) < \eps n^r$.
We define $\Psi^0$ by averaging over each
choice of $\phi^A \in A(\Phi)$. 
Then for any orbit $O$ we have 
$\mb{P}(O \sub \phi^A\Ss) < q^r \oO^{-1} n^{-r}$,
as $\Phi$ is $(\oO,s)$-extendable,
so $U(\Psi^0)_O < \sum_A |k_A| q^r n^{-r}
< q^r \oO^{-1} \eps$ and similarly
$U(J^0)_O \le q^{2r} \oO^{-1} \eps$.

By Lemma \ref{lattice:null:Q} we can construct 
$\Psi^i$ and $J^i = J^{i-1}-\pl^\gG\Psi^i$
as in the proof of Lemma \ref{lattice}
so that all $U(\Psi^i)_O \le \eps_i$
and $U(J^i)_O \le q^r \eps_i$,
where $\eps_0 = q^r \oO^{-1} \eps$ and 
$\eps_i = q^{r+2i} C(i,\oO) \oO^{-1} \eps_{i-1}$.
Then all $U(\Psi)_O \le 2 q^{2r} \oO^{-1} \eps
\prod_{i \in [r]} (q^{r+2i} C(i,\oO) \oO^{-1})
< C(q,\oO) \eps$, recalling that 
$C(i,\oO) = 2^{C(i)} \oO^{-(9i)^{i+4}}$
and $C(i) = 2^{(9i+2)^{i+5}}$. \qed

\section{Bounded integral decomposition} \label{sec:bddint}

To complete the proof of Theorem \ref{main},
it remains to prove Lemma \ref{bddint}.
The high-level strategy is similar
to the randomised rounding and focussing argument
from \cite{Kexist} (version 1), although there are
some additional complications in the general setting.
The proof is by induction on $q$.
In the inductive step we can assume
Lemma \ref{bddint:avoid} for smaller values of $q$
(this will be used in the proof of Lemma \ref{bddint:reduce}).
Note that we do not assume that $\gG$ is elementary,
as this property is not preserved by the inductive step.

\subsection{Proof modulo lemmas}

We start by stating two key lemmas
and using them to deduce Lemma \ref{bddint};
the remainder of the section will then
be devoted to proving the key lemmas.
The first lemma is an approximate version 
of Lemma \ref{bddint}; the second will allow
us to focus the support in a smaller set of vertices.
The proof of Lemma \ref{bddint} is then to alternate
applications of these lemmas until the support 
is sufficiently small that it suffices to use 
the total use quantitative version of the 
decomposition lattice lemma.
In the statements of the lemmas we denote
the labelled complex by $\Phi$, but note
that they will be applied to restrictions
of $\Phi$ as in the statement of Lemma \ref{bddint},
so we allow for weaker lower bounds on
the extendability and number of vertices.
Throughout we fix a $\Ss^\le$-family $\mc{A}$
with $\Ss \le S_q$ and $|\mc{A}| \le K$,
suppose $\gG \in (\mb{Z}^D)^{\mc{A}_r}$ and let 
$\oO_0:=\oO_0(q,D,K)$, $n_0:=n_0(q,D,K)$ be as in Lemma \ref{bddint}.

\begin{lemma} \label{bddint:approx}
Let $\Phi$ be an $(\oO',h)$-extendable $\Ss$-adapted $[q]$-complex
on $[n]$, where $n^{-h^{-2q}}<\oO'<\oO_0$ and $n>n_0^{1/2r}$.
Suppose $J \in \bgen{\gG(\Phi)}$ is $\tT$-bounded,
with $n^{-(4hq)^{-r}}<\tT<1$. 
Then there is some $(\oO')_q^{3qh} \tT$-bounded $J' \in (\mb{Z}^D)^{\Phi_r}$ 
and $(\oO')_q^{-1} \tT$-bounded $\Psi \in \mb{Z}^{\mc{A}(\Phi)}$ 
with $\pl^\gG \Psi = J-J'$.
\end{lemma}

\begin{lemma} \label{bddint:reduce}
Let $\Phi$ be an $(\oO',h)$-extendable $\Ss$-adapted $[q]$-complex
on $[n]$, where $n^{-h^{-2q}}<\oO'<\oO_0$ and $n>n_0^{1/2r}$.
Let $V' \sub V(\Phi)$ with $|V'|=n/2$
be such that $(\Phi,V')$ is $(\oO',h)$-extendable wrt $V'$.
Suppose $J \in \bgen{\gG(\Phi)}$ is $\tT$-bounded,
with $n^{-(3hq)^{-r}} < \tT < (\oO')_q^{3qh}$. Then there is some 
$(\oO')_q^{-2qh} \tT$-bounded $J' \in (\mb{Z}^D)^{\Phi[V']_r}$
and $(\oO')_q^{-2qh} \tT$-bounded $\Psi \in \mb{Z}^{\mc{A}(\Phi)}$
with $\pl^\gG \Psi = J-J'$.
\end{lemma}

\nib{Proof of Lemma \ref{bddint}.}
Let $\mc{A}$ be a $\Ss^\le$-family 
with $\Ss \le S_q$ and $|\mc{A}| \le K$
and suppose $\gG \in (\mb{Z}^D)^{\mc{A}_r}$.
Let $\Phi$ be an $(\oO,h)$-extendable $\Ss$-adapted $[q]$-complex
on $[n]$, where $n^{-h^{-3q}}<\oO<\oO_0$ and $n>n_0$.
Suppose $J \in \bgen{\gG(\Phi)}$ is $\tT$-bounded, with 
$n^{-(5hq)^{-r}} < \tT < 1$.
We need to show that there is some $\oO_q^{-2h} \tT$-bounded
$\Psi \in \mb{Z}^{\mc{A}(\Phi)}$ with $\pl^\gG \Psi = J$. 

Let $t$ be such that $n^{1/2r}/2 < 2^{-t}n \le n^{1/2r}$.
Choose $V_t \sub \dots \sub V_1 \sub V_0 = V(\Phi)$
with $|V_i|=2^{-i}n$ uniformly at random.
By Lemma \ref{extrandom2} (a simple concentration argument 
given in the next subsection) whp all $(\Phi[V_i],V_{i+1})$ 
are $(\oO',h)$-extendable wrt $V_{i+1}$, 
where $\oO'=(\oO/2)^h > n^{-h^{-2q}}$.
We define $\tT$-bounded 
$J_i \in \bgen{\gG(\Phi[V_i])}$ as follows.

Let $J_0=J$. Given $J_i$ with $0 \le i < t$, 
we apply Lemma \ref{bddint:approx} to obtain some 
$(\oO')_q^{3qh} \tT$-bounded $J'_i \in (\mb{Z}^D)^{\Phi[V_i]_r}$ 
and $(\oO')_q^{-1} \tT$-bounded 
$\Psi_i \in \mb{Z}^{\mc{A}(\Phi[V_i])}$
with $\pl^\gG \Psi_i = J_i-J'_i$.
Note that $J'_i \in \bgen{\gG(\Phi[V_i])}$.

Next we apply Lemma \ref{bddint:reduce} to $J'_i$ 
(with $(\oO')_q^{3qh}\tT$ in place of $\tT$) 
to obtain some $\tT$-bounded 
$J_{i+1} \in (\mb{Z}^D)^{\Phi[V_{i+1}]_r}$ 
and $\tT$-bounded $\Psi'_i \in \mb{Z}^{\mc{A}(\Phi[V_i])}$ 
with $\pl^\gG \Psi'_i = J'_i-J_{i+1}$.

To continue the process, we need to show
$J_{i+1} \in \bgen{\gG(\Phi[V_{i+1}])}$.
To see this, first note $J_{i+1} \in \bgen{\gG(\Phi)}
= \mc{L}_\gG(\Phi)$ by Lemma \ref{lattice}.
Now for any $O \in \Phi[V_{i+1}]/\Ss$,
by definition of $\mc{L}_\gG(\Phi)$
we have $(J_{i+1}^\sharp)^O \in \sgen{\gG^\sharp[O]}$,
so $J_{i+1} \in \mc{L}_\gG(\Phi[V_{i+1}])
= \bgen{\gG(\Phi[V_{i+1}])}$, 
again by Lemma \ref{lattice}.

We conclude with some $\tT$-bounded 
$J_t \in \bgen{\gG(\Phi[V_t])}$, where $|V_t| \le n^{1/2r}$.
By Lemma \ref{lattice:total} there is 
$\Psi_t \in \mb{Z}^{\mc{A}(\Phi[V_t])}$ such that 
$\pl^\gG \Psi_t = J_t$ and $U(\Psi_t) \le 2^{(9q)^{q+2}} U(J_t)$.
Let $\Psi = \Psi_t + \sum_{i=0}^{t-1} (\Psi_i+\Psi'_i)$.
Then $\pl^\gG \Psi = J_t + \sum_{i=0}^{t-1} (J_i-J'_i + J'_i-J_{i+1}) = J$.

Also, for any $\psi \in \Phi_{r-1}$ we have
$U(\Psi)_\psi \le U(\Psi_t)_\psi  + 
\sum_{i=0}^{t-1} ( U(\Psi_i)_\psi + U(\Psi'_i)_\psi )
< 2^{(9q)^{q+2}} \tT (n^{1/2r})^r + \sum_{i=0}^{t-1}
( (\oO')_q^{-1} \tT 2^{-i}n + \tT 2^{-i}n ) < 2 (\oO')_q^{-1}\tT n$, 
so $\Psi$ is (say) $\oO_q^{-2h} \tT$-bounded. \qed

\subsection{Random subgraphs}

In the next subsection we will extend the
rational decomposition lemma (Lemma \ref{lattice:Q})
to a version relative to a sparse random subgraph $L$.
We establish some preliminary properties of $L$
in this subsection. First we show that whp $L$ is
`typical' in $\Phi$, in that specifying that certain
edges of an extension should belong to $L$ 
scales the number of extensions in the expected way.

\begin{defn}
Let $\Phi$ be an $[q]$-complex and $L \sub \Phi^\circ_r$.
Let $d_\Phi(L) = |L| |\Phi^\circ_r|^{-1}$.
We say $L$ is $(c,s)$-typical in $\Phi$
if for any $\Phi$-extension $E=(H,F,\phi)$ of rank $s$ 
and $H' \sub H^\circ_r \sm H^\circ_r[F]$ we have
$X_{E,H'}(\Phi,L) = (1 \pm c) d_\Phi(L)^{|H'|} X_E(\Phi)$.
\end{defn}

\begin{lemma} \label{extrandom}
Let $\Phi$ be an $(\oO,s)$-extendable $[q]$-complex.
Suppose $L$ is $\nu$-random in $\Phi^\circ_r$,
where $\nu > n^{-(3sq)^{-r}}$, $n=|V(\Phi)|$.
Then whp $L$ is $(n^{-1/3},s)$-typical in $\Phi$.
In particular, $\Phi[L]$ is $(\oO',s)$-extendable,
where $\oO' = 0.9\nu^{Qs^r} \oO$.
\end{lemma}

\nib{Proof.}
First note by the Chernoff bound that
whp $d_\Phi(L) = (1 \pm n^{-0.4}) \nu$.
Let $E=(H,F,\phi)$ be any $\Phi$-extension of rank $s$,
$H' \sub H^\circ_r \sm H^\circ_r[F]$ and $X = X_{E,H'}(\Phi,L)$.
Note that $\mb{E}X = \nu^{|H'|} X_E(\Phi)$,
where $X_E(\Phi) > \oO n^{v_E}$.
Also, for any $k \in [r]$ there are $O(n^k)$ choices of 
$f \in \Phi^\circ_r$ with $f \sm \phi(F)|=k$, 
and for each such $f$, changing whether 
$f \in L$ affects $X$ by $O(n^{v_E-k})$.
Thus $X$ is $O(n^{2v_E-1})$-varying, so by Lemma \ref{lip3} 
whp $X = (1 \pm n^{-1/3}) \nu^{|H'|} X_E(\Phi)$. 
In particular, $X > \oO' n^{v_E}$. \qed

\medskip

Similarly, we obtain the following variant
form of the previous lemma that was
used in the previous subsection.

\begin{lemma} \label{extrandom2}
Let $\Phi$ be an $(\oO,s)$-extendable $[q]$-complex on $[n]$.
Suppose $S$ is uniformly random in $\tbinom{[n]}{m}$, 
where $m > \oO^{-1}\log n$ and $n$ is large.
Then $(\Phi,S)$ is $((\oO/2)^h,s)$-extendable wrt $S$
with probability at least $1-e^{-(\oO m)^2/20}$.
\end{lemma}

\nib{Proof.}
It suffices to estimate the probability that any simple 
$\Phi$-extension of rank $s$ is $\oO/2$-dense in $(\Phi,S)$.
Let $E=(H,F,\phi)$ be any $\Phi$-extension of rank $s$ 
with $F = V(H) \sm \{x\}$ for some $x \in V(H)$.
Note that $X_E(\Phi,S) 
= \sum_{\phi^+ \in X_E(\Phi)} 1_{\phi^+(x) \in S}$ and
$X_E(\Phi) > \oO n$ as $\Phi$ is $(\oO,s)$-extendable.
Then $X_E(\Phi,S)$ is hypergeometric
with $\mb{E}X_E(\Phi,S) > \oO m$, so 
$\mb{P}(X_E(\Phi,S) < \oO m/2) < e^{-(\oO m)^2/12}$.
The lemma follows by taking a union bound over at most
$qsn^{qs} < e^{(\oO m)^2/48}$ choices of $E$. \qed

\medskip

\def\oldJrandom{2.21}

We also require the following refined notion of boundedness 
that operates with respect to all small extensions in $L$.
The following lemma is analogous to 
\cite[Lemma \oldJrandom]{Kexist}.

\begin{defn} \label{def:bddwrt}
Let $\Phi$ be an $[q]$-complex, $L \sub \Phi^\circ_r$,
and $J \in (\mb{Z}^D)^{\Phi_r}$.
Let $E=(H,F,\phi)$ with $H \sub [q](s)$ be a $\Phi$-extension,
$G \sub H_r^\circ \sm H_r^\circ[F]$, $\psi \in [q](s)_r$ and $e=Im(\psi)$.
We write $X^{e,J}_{E,G}(\Phi,L) 
= \sum_{\phi^* \in X_{E,G}(\Phi,L)} U(J)_{\phi^*\psi}$.
We say that $J$ is $(\tT,s)$-bounded wrt $(\Phi,L)$
if $X^{e,J}_{E,G}(\Phi,L) < \tT d_\Phi(L)^{|G|} |V(\Phi)|^{v_E}$
for any such $E$, $G$ and $e$ 
with $e \notin G$ and $e \sm F \ne \es$.
\end{defn}

\begin{lemma} \label{Jrandom}
Let $\Phi$ be an $[q]$-complex with $|V(\Phi)|=n$.
Suppose $J \in (\mb{Z}^D)^{\Phi_r}$ is $\tT$-bounded,
with $\tT>n^{-0.01}$, and all $U(J)_\psi<n^{0.1}$. 
Let $L$ be $\nu$-random in $\Phi^\circ_r$, 
where $\nu > n^{-(3sq)^{-r}}$.
Then whp $J$ is $(1.1\tT, s)$-bounded wrt $(\Phi,L)$.
\end{lemma}

\nib{Proof.}
Let $E=(H,F,\phi)$ with $H \sub [q](s)$ be a $\Phi$-extension,
$G \sub H_r^\circ \sm H_r^\circ[F]$ and $\psi \in [q](s)_r$ 
with $e:=Im(\psi) \notin G$ and $e \sm F \ne \es$.
Write $X = X^{e,J}_{E,G}(\Phi,L) = \sum_{\phi^* \in X_E(\Phi)} 
1_{\phi^* \in X_{E,G}(\Phi,L)} U(J)_{\phi^*\psi}$.
As $J$ is $\tT$-bounded we have
$\sum_{\phi^* \in X_E(\Phi)} U(J)_{\phi^*\psi} < \tT n^{v_E}$.
For each $\phi^* \in X_E(\Phi)$ we have
$\mb{P}(\phi^* \in X_{E,G}(\Phi,L)) = \nu^{|G|}$,
so $\mb{E}X < \tT \nu^{|G|} n^{v_E}$.
For any $k \in [r]$ there are $O(n^k)$ choices of 
$f \in \Phi^\circ_r$ with $|f \sm \phi(F)|=k$, 
and for each such $f$, changing whether 
$f \in L$ affects $X$ by $O(n^{v_E-k+0.1})$.
Thus $X$ is $O(n^{2v_E-0.8})$-varying, so by Lemma \ref{lip3} 
whp $X < 1.1\tT d_\Phi(L)^{|G|} n^{v_E}$. \qed

\subsection{Rational decompositions}

In this subsection we prove the following result,
which is a version of Lemma \ref{lattice:Q}
relative to a sparse random subgraph $L$;
note the key point that we incur a loss in
boundedness that depends only on $q$,
not on the density of $L$.
Throughout, as in the hypotheses of Lemma \ref{bddint:approx},
we let $\mc{A}$ be a $\Ss^\le$-family 
with $\Ss \le S_q$ and $|\mc{A}| \le K$,
suppose $\gG \in (\mb{Z}^D)^{\mc{A}_r}$,
and let $\Phi$ be an $(\oO,h)$-extendable $\Ss$-adapted $[q]$-complex
on $[n]$, where $n^{-h^{-2q}}<\oO<\oO_0$ and $n>n_0^{-1/2}$.
(For convenient notation we rename $\oO'$ as $\oO$.)

\begin{lemma} \label{bddrat:L}
Suppose $L \sub \Phi^\circ_r$ is $(c,h)$-typical in $\Phi$
with $d_\Phi(L) = \nu \ge n^{-(3hq)^{-r}}$.
Let $J \in \bgen{\gG(\Phi)}_{\mb{Q}} \cap \mb{Q}^{\Phi[L]_r}$.
Suppose $J$ is $\tT$-bounded with all $|J_\psi| < n^{0.1}$
and $(\tT, h)$-bounded wrt $(\Phi,L)$.
Then there is some $\oO_q^{-0.9} \tT$-bounded
$\Psi \in \mb{Q}^{\mc{A}(\Phi[L])}$ with $\pl^\gG \Psi = J$.
\end{lemma}

The proof of Lemma \ref{bddrat:L} uses the
following result, which reduces to the case
when we bound the use of every orbit.

\begin{lemma} \label{flat:Q}
For any $\tT$-bounded $J \in (\mb{Q}^D)^{\Phi_r}$
there is some $J' \in (\mb{Q}^D)^{\Phi_r}$
and $\Psi \in \mb{Q}^{\mc{A}(\Phi)}$ 
such that $\pl \Psi = J-J'$,
and for any $O \in \Phi_r/\Ss$ both 
$U(\Psi)_O - U(J)_O$ and $U(J')_O$ 
are at most $q^q \oO^{-1} \tT$.
\end{lemma}

\nib{Proof.}
For each $O \in \Phi_r/\Ss$
we fix a representative $\psi^O \in \Phi_{B^O}$
and $n^O \in \mb{Q}^{\mc{A}_{B^O}}$ 
with $|n^O| = U(J)_O$ and
$f_{B^O}(J)_{\psi^O} = \sum_\tT n^O_\tT \gG^\tT$,
so $J^O = \sum_\tT n^O_\tT \gG(\psi^O\tT^{-1})$.
For each such $(O,\tT)$ we let
$E^O_\tT = (\ova{q},Im(\tT),\psi^O \tT^{-1})$.
We let $\Psi = \sum_{O,\tT} n^O_\tT 
|X_{E^O_\tT}(\Phi)|^{-1} X_{E^O_\tT}(\Phi)$,
where if $\tT \in A \in \mc{A}$ we choose
the copy of $\phi \in A(\Phi)$ 
for each $\phi \in X_{E^O_\tT}(\Phi)$.
Then $J$ is exactly cancelled by the $\gG$-atoms 
of $\pl^\gG \Psi$ corresponding to
$\gG(\phi)^O$ for $\phi \in A(\Phi)$.
To estimate the remaining contributions
of $\pl^\gG \Psi$, note that for any 
$O \in \Phi_r/\Ss$ and $r' \in [r]$ we have 
$\sum \{ |n^{O'}_\tT|: |Im(O') \sm Im(O)|=r' \} 
\le \tbinom{r}{r'} \tT n^{r'}$. For each such $(O',\tT)$,
there are at least $\oO n^{q-r}$ choices of $\phi$,
of which at most $(q-r)!n^{q-r-r'}$ contain $Im(O)$,
so for random $\phi \in X_{E^{O'}_\tT}(\Phi)$ we have
$\mb{P}(O \sub \phi\Ss) \le (q-r)! \oO^{-1} n^{-r'}$.
Summing over $r'$ gives the stated bounds on
$U(\Psi)_O - U(J)_O$ and $U(J')_O$. \qed

\medskip

\nib{Proof of Lemma \ref{bddrat:L}.}
We start by defining $\Psi^0 \in \mb{Q}^{\mc{A}(\Phi)}$ 
such that $\pl^\gG \Psi^0 = J$ and 
$U(\Psi^0)_O < U(J)_O + (C^*-1) \tT$
for all $O \in \Phi_r/\Ss$,
where $C^* = 2C(q,\oO) q^q \oO^{-1}$
(so $\Psi^0$ is $C^*\tT$-bounded).
Then we will modify $\Psi^0$ to obtain $\Psi$
using a version of the Clique Exchange Algorithm. 

First we apply Lemma \ref{flat:Q}
to obtain $\Psi' \in \mb{Q}^{\mc{A}(\Phi)}$
and $J' = J - \pl^\gG \Psi$ so that
all $U(\Psi)_O - U(J)_O$ and $U(J')_O$ 
are at most $q^q \oO^{-1} \tT$.
Then by Lemma \ref{lattice:Q} there is 
$\Psi'' \in \mb{Q}^{\mc{A}(\Phi)}$ 
such that $\pl^\gG \Psi'' = J'$ and all
$U(\Psi'')_O < C(q,\oO) q^q \oO^{-1} \tT$.
We take $\Psi^0 = \Psi' + \Psi''$.

We apply two Splitting Phases,
the first in $\Phi$ and the second in $\Phi[L]$.
For the first, we fix $N_0 \in \mb{N}$ such that
$N_0\Psi^0 \in \mb{Z}^{\mc{A}(\Phi)}$, and list the signed 
elements of $N_0\Psi^0$ as $(s_i \phi_i: i \in [|N_0\Psi^0|])$,
where each $s_i \in \pm 1$.
For each $i$, say with $\phi_i \in A^i(\Phi)$, we consider 
the $\Phi$-extension $E_i = ([q](p),[q],\phi_i)$,
and define $\Psi^1 \in \mb{Q}^{\mc{A}(\Phi)}$ by
$\Psi^1 = \Psi^0 + N_0^{-1} \sum_{i \in [|N_0 \Psi^0|]} 
s_i \mb{E}_{\phi^*_i \in X_{E_i}(\Phi)}
( A^i(\Phi[\phi^*_i \Ups']) - A^i(\Phi[\phi^*_i \Ups]) )$.
Then $\pl^\gG \Psi^1 = \pl^\gG \Psi^0 = J$, 
and all signed elements in $\Psi^0$ are cancelled.

We claim for any $O \in \Phi_r/\Ss$ that 
$\GG_O := U(\Psi^1)_O - U(\Psi^0)_O 
\le r!\oO^{-1} (2pq)^r C^* \tT$.
To see this, we estimate
$\GG_O \le \mb{E}_{i \in [|N_0 \Psi^0|]} 
\mb{P}(Im(O) \in \phi^*_i(\OO'))$
(recall $\OO' = K^r_q(p) \sm Q$).
For any $r' \in [r]$, 
as $\Psi^0$ is $C^*\tT$-bounded
there are at most $N_0 \tbinom{r}{r'} C^* \tT n^{r'}$
choices of $i$ such that $|Im(O) \sm Im(\phi_i))|=r'$.
For each such $i$, as $\Phi$ is $(\oO,h)$-extendable
there are at least $\oO n^{pq-q}$ choices of
$\phi^*_i \in X_{E_i}(\Phi)$, of which at most
$r!|\OO'| n^{v_{E_i}-r'}$ have $Im(O) \in \phi^*_i(\OO')$,
so $\mb{P}(Im(O) \in \phi^*_i(\OO')) < r!\oO^{-1} |\OO'| n^{-r'}$.
Therefore $\GG_O \le N_0^{-1} \sum_{r' \in [r]}
N_0 \tbinom{r}{r'} C^* \tT n^{r'} \cdot r!\oO^{-1} |\OO'| n^{-r'}
< r!\oO^{-1} (2pq)^r C^* \tT$, as claimed.

Also, for any $\phi \in \mc{A}(\Phi)$ 
we claim that $|\Psi^1_\phi| < \oO^{-1} p^q n^{r-q} 
\sum_{O \sub \phi\Ss} (qC^*\tT + U(J)_O)$.
To see this, we consider separately the
contributions from $\phi_i$ according
to $r'=|Im(\phi) \cap Im(\phi_i)|$.

First we consider $0 \le r' < r$. 
As $\Psi^0$ is $C^*\tT$-bounded there are at most 
$N_0 Q C^* \tT n^{r-r'}$
such choices of $\phi_i$. For each such $i$,
there are at least $\oO n^{pq-q}$ choices of
$\phi^*_i \in X_{E_i}(\Phi)$,
of which at most $p^q n^{pq-q-(q-r'))}$
have $\phi = \phi^*_i\phi'$ for some $\phi' \in [q](p)$,
so the total such contribution to $|\Psi^1_\phi|$
is at most $N_0^{-1} \sum_{r'} N_0 Q C^* \tT n^{r-r'}
\cdot \oO^{-1} p^q n^{r'-q} 
= r p^q Q C^* \oO^{-1} \tT n^{r-q}$.

It remains to consider $r'=r$. For each $O \sub \phi\Ss$
there are at most 
$N_0 ( U(J)_O + (C^*-1) \tT )$ choices of $\phi_i$ containing $Im(O)$.
For each of these with $|Im(\phi_i) \cap Im(\phi)|=r$ 
we have $\phi = \phi^*_i\phi'$ for some $\phi' \in [q](p)$
with probability at most $\oO^{-1} p^q n^{r-q}$,
so the total such contribution to $|\Psi^1_\phi|$ is at most
$(U(J)_O + (C^*-1) \tT ) \oO^{-1} p^q n^{r-q}$.
The claim follows.

In the second Splitting Phase, we fix $N_1 \in \mb{N}$ such that
$N_1\Psi^1 \in \mb{Z}^{\mc{A}(\Phi)}$, and list the signed 
elements of $N_1\Psi^1$ as $(s_i \phi_i: i \in [|N_1\Psi^1|])$,
where each $s_i \in \pm 1$.
For each $i$, say with $\phi_i \in A^i(\Phi)$, we consider 
the $\Phi$-extension $E_i = ([q](p),[q],\phi_i)$,
and define $\Psi^2 \in \mb{Q}^{\mc{A}(\Phi)}$ by
$\Psi^2 = \Psi^1 + N_1^{-1} \sum_{i \in [|N_1 \Psi^1|]} 
s_i \mb{E}_{\phi^*_i \in X_{E_i,L}(\Phi)}
( A^i(\Phi[\phi^*_i \Ups']) - A^i(\Phi[\phi^*_i \Ups]) )$. 
Then $\pl \Psi^2 = \pl \Psi^1 = J$, 
and all signed elements in $\Phi^1$ are cancelled,
so $\Psi^2$ is supported on maps $\phi$ such that
there is at most one $e \in Q$ with $\phi(e) \notin L$.
(This was the same as the first Splitting Phase except
that we changed $X_{E_i}(\Phi)$ to $X_{E_i}(\Phi,L)$.)

We claim for any $O \in \Phi_r/\Ss$ that 
$\GG'_O := U(\Psi^2)_O - U(\Psi^1)_O 
\le (pq)^{2q} \oO^{-2} C^*\tT \nu^{-1}$,
where $\nu:=d_\Phi(L)$.
To see this, we estimate
$\GG'_O \le \mb{E}_{i \in [|N_1 \Psi^1|]} 
\mb{P}(Im(O) \in \phi^*_i(\OO'))$.
We fix $f' \in \OO'$ and consider the
contribution from $i$ with $O = \phi^*_if' \Ss$.
Consider any $\phi' \in \Phi$ with $O = \phi'f' \Ss$
and the $\Phi$-extension $E_{f'}=(H,f',\phi')$, 
where $H = [q](p)[[q] \cup f']$
is the restriction of $[q](p)$ to $[q] \cup f'$.
Let $H' = H^\circ_r \sm ([q]_r \cup \{f'\})$.

The number of $i$ with $\phi_i=\phi^*\mid_{[q]}$ 
for some $\phi^* \in X_{E_{f'},H'}(\Phi,L)$ is at most
\begin{align*}
N_1 \sum_{\phi^* \in X_{E_{f'},H'}(\Phi,L)} |\Psi^1_{\phi^*\mid_{[q]}}|
& < \sum_{\phi^* \in X_{E_{f'},H'}(\Phi,L)} N_1 \oO^{-1} p^q n^{r-q} 
\sum_{\psi \sub \phi^*\mid_{[q]}} (qC^*\tT + U(J)_\psi) \\
& = N_1 \oO^{-1} p^q n^{r-q} \brak{ qQC^*\tT X_{E_{f'},H'}(\Phi,L)
+ \sum_{e' \in Q} X^{e',J}_{E_{f'},H'}(\Phi,L) } \\
& <  2N_1 \oO^{-1} p^q n^{r-q} 
qQC^*\tT \nu^{|H'|} n^{q-|f' \cap [q]|},
\end{align*}
as $L$ is $(c,h)$-typical and 
$J$ is $(\tT,h)$-bounded wrt $L$.

For each such $i$, there are at least
$0.9 \oO n^{pq-q} \nu^{|\OO'|}$ choices
of $\phi^*_i \in X_{E_i}(\Phi,L)$,
of which the number containing $\phi'$ is at most 
$1.1 \nu^{|\OO'|-|H' \cup \{f'\}|} n^{pq-q-r+|f' \cap [q]|}$,
so $\mb{P}(\phi' \sub \phi^*_i)
< 1.3  \oO^{-1} \nu^{-|H'|-1} n^{|f' \cap [q]-r}$.
Thus 
\begin{align*} 
\GG'_O & < q^r N_1^{-1} \sum_{f' \in \OO'}
2N_1 \oO^{-1} p^q n^{r-q} qQC^*\tT \nu^{|H'|} n^{q-|f' \cap [q]|}
\cdot 1.3  \oO^{-1} \nu^{-|H'|-1} n^{|f' \cap [q]|-r} \\
& <  (pq)^{2q} \oO^{-2} C^*\tT \nu^{-1},
\ \text{ using } 
p>2^{8q},  \text{ as claimed.}
\end{align*}
We fix $N_2 \in \mb{N}$ 
such that $N_2\Psi^2 \in \mb{Z}^{\mc{A}(\Phi)}$,
and classify signed elements of $N_2\Psi^2$ as before:
recall that a pair $(O,\phi')$ is near or far,
has the same sign as that of $\phi'$ in $N_2\Psi^2$,
and has a type $\tT$ determined 
by an orbit representative $\psi^O \in O$.
For $\phi \in \mc{A}(\Phi)$ let $B_\phi$ 
be the number of pairs $(O,\phi)$ in $N_2\Psi^2$
such that $Im(O) \notin L$; note that all such pairs
are near and $Im(O)$ is uniquely determined by $\phi$.

We claim that each
$B_\phi < 3N_2  \oO^{-2} p^q  \nu^{-Q+1} n^{r-q}
\sum_{\psi \sub \phi} (qQC^*\tT + U(J)_\psi)$,
To see this, we fix
$\psi' \in \Ups' \cup \Ups \sm \{[q]\}$
and consider the contributions
from $i$ with $\phi^*_i \psi'=\phi$.
We consider the $\Phi$-extension $E_{\psi'}=(H,F',\phi')$,
where $F'=Im(\psi')$ (so $|F' \cap [q]|=r$), 
$H=[q](p)[[q] \cup F']$ and $\phi = \phi'\psi'$.
Let $H' = H^\circ_r \sm ([q]_r \cup F'_r)$.

The number of $i$ with $\phi_i=\phi^*\mid_{[q]}$ 
for some $\phi^* \in X_{E_{\psi'},H'}(\Phi,L)$ is at most
\begin{align*}
N_1 \sum_{\phi^* \in X_{E_{\psi'},H'}(\Phi,L)} |\Psi^1_{\phi^*\mid_{[q]}}|
& \le \sum_{\phi^* \in X_{E_{\psi'},H'}(\Phi,L)} N_1 \oO^{-1} p^q  n^{r-q} 
\sum_{\psi \sub \phi^*\mid_{[q]}} (qC^*\tT + U(J)_\psi)\\
& = N_1 \oO^{-1} p^q n^{r-q} \brak{ qQC^*\tT X_{E_{\psi'},H'}(\Phi,L)
+ \sum_{e' \in Q} X^{e',J}_{E_{\psi'},H'}(\Phi,L) } \\
& <  2N_1 \oO^{-1} p^q  n^{r-q} 
\nu^{|H'|} n^{q-|F' \cap [q]|} (qQC^*\tT + U(J)_O) \\
& = 2N_1 \oO^{-1} p^q \nu^{|H'|} (qQC^*\tT + U(J)_O),
\end{align*}
where for $e' = F' \cap [q]$ we let 
$O \sub \phi\Ss$ be such that $Im(O)=\phi(e')$ 
and use $X^{e',J}_{E_{\psi'},H'}(\Phi,L) 
\le X_{E_{\psi'},H'}(\Phi,L) U(J)_O$.

For each such $i$, there are at least
$0.9 \oO n^{pq-q} \nu^{|\OO'|}$ choices
of $\phi^*_i \in X_{E_i}(\Phi,L)$,
of which the number containing $\psi'$ is at most 
$1.1 \nu^{|\OO'|-|H'|-(Q-1)} n^{pq-q-q+|F' \cap [q]|}$
(as $|F'_r \sm [q]_r| = Q-1$),
so $\mb{P}(\psi' \sub \phi^*_i)
< 1.3  \oO^{-1} \nu^{-|H'|-Q+1} n^{r-q}$.
Then 
\begin{align*}
B_\phi & <  N_2 N_1^{-1} \sum_{\psi':|F' \cap [q]|=r}
2N_1 \oO^{-1} p^q \nu^{|H'|} (qQC^*\tT + U(J)_O)
\cdot  1.3  \oO^{-1} \nu^{-|H'|-Q+1} n^{r-q} \\
& < 3N_2  \oO^{-2} p^q  \nu^{-Q+1} n^{r-q}
\sum_{\psi \sub \phi} (qQC^*\tT + U(J)_\psi), \text{ as claimed.}
\end{align*}

In Elimination Phase, we consider each orbit 
$O \in \Phi_r/\Ss$ with $Im(O) \notin L$,
say with representative $\psi^O \in \Phi_B$,
and let $E_O = (B(2),B,\psi^O)$.
For each $\psi^*_O \in X_{E_O}(\Phi,L)$
and each signed near pair $\pm (O,\phi)$ in $N_2\Psi^2$, 
say of type $\tT$, i.e.\ $\phi\tT=\psi^O$,
we let $E^\phi(\psi^*_O) = (w^{B^\phi_O},F,\phi_0)$,
where $B^\phi_O = \tT(B)$,
$F = [q] \cup (B^\phi_O \times [2])$,
$\phi_0\mid_{[q]} = \phi$ and
$\phi_0(\tT(x),y)=\psi^*_O(x,y)$
for $x \in B$, $y \in [2]$.

We let $\Psi = \Psi^2 + N_2^{-1} \sum_{(O,\phi)} \pm
\mb{E}_{\psi^*_O \in X_{E_O}(\Phi,L)}
\mb{E}_{\psi^+ \in X_{E^\phi(\psi^*_O)}(\Phi,L)}
\sum_{x \in [q](s)} w^{B^\phi_O}_x \psi^* x$,
where the sign is that of $(O,\phi)$,
and each $\psi^* x \in A^\phi(\Phi)$
where $\phi \in A^\phi(\Phi)$.
Then $\pl^\gG \Psi = J$, similarly to the 
previous version of the algorithm,
treating all near pairs on $O$ 
as a single cancelling group,
which is valid as $J^O=0$ when $Im(O) \notin L$.
All pairs $(O,\phi)$ with $Im(O) \notin L$ are cancelled, 
so $\Psi \in \mb{Q}^{\mc{A}(\Phi[L])}$.

We claim for any $O' \in \Phi_r/\Ss$ that 
$\GG''_{O'} := U(\Psi)_{O'} - U(\Psi^2)_{O'} 
\le \oO^{-3} (pq)^{2q} C^*\tT \nu^{-1}$.
To see this, we estimate
\[ \GG''_{O'} \le N_2^{-1} \sum_{(O,\phi)}
\mb{E}_{\psi^*_O \in X_{E_O}(\Phi,L)}
\mb{P}( Im(O') \in \psi^+(\OO') )
=  N_2^{-1} \sum_{(O,\phi)}
\mb{P}( Im(O') \in \psi^\phi_O(\OO') ),\]
with each $\psi^\phi_O$ uniformly random
in $X_{E^\phi_O}(\Phi,L)$ where
$E^\phi_O = ( [q](s), [q], \phi)$.

We fix $f' \in [q](s)_r \sm \ova{q}_r$ and consider the
contribution from near pairs $(O,\phi)$ 
with $O' = \psi^\phi_O f' \Ss$.
Consider any $\phi' \in \Phi_r$ with $O' = \phi' f' \Ss$
and the $\Phi$-extension $E_{f'}=(H,f',\phi')$, 
where $H = [q](p)[[q] \cup f']$. For $B \in [q]_r$
let $H^B = H^\circ_r \sm \{B,f'\}$.

The number of signed near pairs
$\pm (O,\phi)$ in $N_2 \Psi^2$
with $\phi = \phi^*\mid_{[q]}$, $Im(O) \notin L$,
$O=\phi^*\mid_B \Ss$ for some $B \in [q]_r$ and
$\phi^* \in X_{E_{f'},H^B}(\Phi,L)$ is at most 
\begin{align*}
\sum_{\phi^* \in X_{E_{f'},H^B}(\Phi,L)} B_{\phi^*\mid_{[q]}} 
& < \sum_{\phi^* \in X_{E_{f'},H^B}(\Phi,L)}
3N_2 \oO^{-2} p^q  \nu^{-Q+1} n^{r-q}
\sum_{\psi \sub \phi} (qQC^*\tT + U(J)_\psi) \\
& = 3N_2 \oO^{-2} p^q  \nu^{-Q+1} n^{r-q} \sum_{e^* \in Q}
[ qQC^*\tT X_{E_{f'},H^B}(\Phi,L) + X_{E_{f'},H^B}^{e^*,J}(\Phi,L) ]\\
& < 3N_2 \oO^{-2} p^q  \nu^{-Q+1} n^{r-q} 2Q 
qQC^*\tT \nu^{|H^B|} n^{q-|f' \cap [q]|},
\end{align*}
as $L$ is $(c,h)$-typical and 
$J$ is $(\tT,h)$-bounded wrt $L$.

For each such $(O,\phi)$, there are at least 
$0.9 \oO n^{qs-q} \nu^{Q(s^r-1)}$ choices for 
$\psi^\phi_O \in X_{E^\phi_O}(\Phi,L)$,
of which the number with $O' = \psi^\phi_O f' \Ss$ 
is at most $1.1 n^{qs-q-r+|f' \cap [q]|} \nu^{Qs^r-|H|}$,
so $\mb{P}(O' = \psi^\phi_O f'\Ss) 
< 1.3 \oO^{-1} \nu^{-|H^B|+Q-2} n^{-r+|f' \cap [q]|}$,
as $Qs^r-|H|=Q(s^r-1)-|H^B|+Q-2$. Thus 
\begin{align*}
\GG''_{O'} &< N_2^{-1} \sum_{f' \in [q](s)_r \sm \ova{q}_r}
3N_2 \oO^{-2} p^q  \nu^{-Q+1} n^{r-q} 2Q 
qQC^*\tT \nu^{|H^B|} n^{q-|f' \cap [q]|} \\
& \hspace{3cm} \cdot
1.3 \oO^{-1} \nu^{-|H^B|+Q-2} n^{-r+|f' \cap [q]|} \\
& < \oO^{-3} (pq)^{2q} C^*\tT \nu^{-1}, \ \text{ as claimed.}
\end{align*}

Finally, for any $f \in \Phi_{r-1}$ we have
\begin{align*}
U(\Psi)_f & \le U(\Psi^0)_f
+ \sum \{ \GG_O: f\Ss \sub O \}
+ \sum \{ \GG'_O + \GG''_O:
 f\Ss \sub O, Im(O) \in L \} \\
& < C^*\tT n 
+ r!\oO^{-1} (2pq)^r C^* \tT n
+ 1.1 q^r \nu n ( 
(pq)^{2q} \oO^{-2} C^*\tT \nu^{-1}
+  \oO^{-3} (pq)^{2q} C^*\tT \nu^{-1} ) \\
& < 2\oO^{-3} q^r (pq)^{2q} C^*\tT n.
\end{align*}
Recalling that
$C^* = 2C(q,\oO) q^q \oO^{-1}$,
$C(i,\oO) = 2^{C(i)} \oO^{-(9i)^{i+4}}$,
$C(i) = 2^{(9i+2)^{i+5}}$,
$\oO_q:=\oO^{(9q)^{q+5}}$,
and $\oO<\oO_0$ we see that
$\Psi$ is $\oO_q^{-0.9} \tT$-bounded. \qed

\subsection{Approximation}

In this subsection we prove Lemma \ref{bddint:approx}
(approximate bounded integral decomposition)
by randomly rounding Lemma \ref{bddrat:L}
(rational decomposition with respect
to a sparse random subgraph).
Throughout, as in the hypotheses of Lemma \ref{bddint:approx},
we let $\mc{A}$ be a $\Ss^\le$-family 
with $\Ss \le S_q$ and $|\mc{A}| \le K$,
suppose $\gG \in (\mb{Z}^D)^{\mc{A}_r}$,
and let $\Phi$ be an $(\oO,h)$-extendable $\Ss$-adapted $[q]$-complex
on $[n]$, where $n^{-h^{-2q}}<\oO<\oO_0$ and $n>n_0$.
We start with some preliminary lemmas 
for flattening and focussing a vector.

\begin{lemma} \label{bddint:flat}
If $J \in (\mb{Z}^D)^{\Phi_r}$ is
$\tT$-bounded with $\tT > n^{-1/2}$ then
there is some $J' \in (\mb{Z}^D)^{\Phi_r}$
and $\Psi \in \mb{Z}^{\mc{A}(\Phi)}$ 
such that $\pl \Psi = J-J'$,
$J'$ and $\Psi$ are $q^q \oO^{-1} \tT$-bounded,
and $U(J')_O < n^{0.1}$ for all $O \in \Phi_r/\Ss$.
\end{lemma}

\nib{Proof.}
Similarly to the proof of Lemma \ref{flat:Q},
for each $O \in \Phi_r/\Ss$ with 
representative $\psi^O \in \Phi_{B^O}$, 
we have $n^O \in \mb{Z}^{\mc{A}_{B^O}}$ 
with $|n^O| = U(J)_O$ and
$f_{B^O}(J)_{\psi^O} = \sum_\tT n^O_\tT \gG^\tT$,
so $J^O = \sum_\tT n^O_\tT \gG(\psi^O\tT^{-1})$.
Let $S$ be the intset where each $(O,\tT)$
appears $|n^O_\tT|$ times with the sign of $n^O_\tT$.
For each $(O,\tT)$ in $S$ we add to $\Psi$
with the same sign as $(O,\tT)$ 
a uniformly random $\phi \in X_E(\Phi)$
with $E = (\ova{q},Im(\tT),\psi^O \tT^{-1})$;
then each $\gG(\psi^O\tT^{-1})$ is cancelled 
by the corresponding $\gG(\phi)^O$,
where $\tT \in A$, $\phi \in A(\Phi)$.

For any $\psi \in \Phi_{r-1}$ and $k \in [r]$ there are at most 
$\tbinom{r-1}{k-1} \tT n^k$ signed elements $(O,\tT)$ 
of $S$ with $|Im(O) \sm Im(\psi)|=k$. For each such $(O,\tT)$,
there are at least $\oO n^{q-r}$ choices of $\phi$,
of which at most $(q-r)!n^{q-r-(k-1)}$ contain $\phi$,
so $\mb{P}(\psi \sub \phi) \le (q-r)! \oO^{-1} n^{-k+1}$. 
Then $U(\Psi)_\psi$ is a sum of bounded independent variables
with mean at most $0.9 q^q \oO^{-1} \tT n$, so whp 
$J'$ and $\Phi$ are $q^q \oO^{-1} \tT$-bounded.
Similarly, whp $U(J')_O < n^{0.1}$ for all $O \in \Phi_r/\Ss$. \qed

\begin{lemma} \label{bddint:random}
Suppose $J \in (\mb{Z}^D)^{\Phi_r}$ 
is $\tT$-bounded with $\tT > n^{-1/2}$.
Let $L \sub \Phi^\circ_r$ be $(c,h)$-typical in $\Phi$
with $d_\Phi(L) > n^{-(3hq)^{-r}}$.
Suppose $J$ is $(\tT, h)$-bounded wrt $(\Phi,L)$.
Then there is some $J' \in (\mb{Z}^D)^{\Phi[L]_r}$
and $\Psi \in \mb{Z}^{\mc{A}(\Phi)}$ 
such that $\pl \Psi = J-J'$, 
$J'$ and $\Psi$ are $2^{2r} \oO^{-1} \tT$-bounded, 
and $J'$ is $(q! 2^{2r} \oO^{-1} \tT, h)$-bounded wrt $(\Phi,L)$.
\end{lemma}

\nib{Proof.}
We apply the same procedure as in the proof 
of Lemma \ref{bddint:flat}, 
replacing $X_E(\Phi)$ by $X_E(\Phi,L)$.
To spell this out, for each orbit $O$ with 
representative $\psi^O \in \Phi_{B^O}$, 
we have $n^O \in \mb{Z}^{\mc{A}_{B^O}}$ 
with $|n^O| = U(J)_O$ and
$f_{B^O}(J)_{\psi^O} = \sum_\tT n^O_\tT \gG^\tT$,
so $J^O = \sum_\tT n^O_\tT \gG(\psi^O\tT^{-1})$.
Let $S$ be the intset where each $(O,\tT)$
appears $|n^O_\tT|$ times with the sign of $n^O_\tT$.
For each $(O,\tT)$ in $S$
we add to $\Psi$ with the same sign as $(O,\tT)$ 
a uniformly random $\phi^O_\tT \in X_E(\Phi,L)$
with $E = (\ova{q},Im(\tT),\psi^O \tT^{-1})$;
then each $\gG(\psi^O\tT^{-1})$ 
is cancelled by $\gG(\phi^O_\tT)^O$,
where $\tT \in A$, $\phi^O_\tT \in A(\Phi)$.

We claim for any $e' \in L$ that 
$E_{e'} :=  \sum_{(O,\tT)} \sum_{B' \in [q]_r \sm \{Im(\tT)\}}
\mb{P}(\phi^O_\tT(B')=e') < 1.3q! 2^r \oO^{-1} \tT d(L)^{-1}$. 
To see this, first note that for any $k \in [r]$,
as $J$ is $(\tT, h)$-bounded wrt $(\Phi,L)$ there are 
at most $\tbinom{r}{k} d(L)^{\tbinom{k+r}{r}-2} \tT n^k$
signed elements $(O,\tT)$ of $S$ with $|Im(O) \sm Im(e')|=k$. 
As $\Phi$ is $(\oO,h)$-extendable 
and $L$ is $(c,h)$-typical in $\Phi$,
for each such $(O,\tT)$, there are at least
$0.9d(L)^{Q-1} \oO n^{q-r}$ choices of $\phi^O_\tT$,
of which at most $1.1q!d(L)^{Q-\tbinom{k+r}{r}} n^{q-r-k}$
have $\phi^O_\tT(B')=e'$ for some $B' \ne Im(\tT)$,
so $\mb{P}(\phi^O_\tT(B')=e') 
< 1.3q!\oO^{-1}d(L)^{1-\tbinom{k+r}{r}} n^{-k}$.
The claim follows.
Now for any $f \in \Phi_{r-1}$,
by typicality $|L(Im(f))| < 1.1 d(L) n$,
so by the claim $U(\Psi)_f$ is a sum of bounded independent variables
with mean at most $1.5 q! 2^r \oO^{-1} \tT n$,
so whp $J'$ and $\Psi$ are
$q! 2^{r+1} \oO^{-1}\tT$-bounded.

Finally, consider any $\Phi$-extension $E=(H,F,\phi)$ 
with $H \sub [q](h)$, any $G \sub H^\circ_r \sm H^\circ_r[F]$ 
and $e \in [q](h)^\circ_r \sm G$ with $e \sm F \ne \es$.
Write $e=Im(\psi)$ with $\psi \in [q](h)$,
$E^+=(H^+,F,\phi)$ with $H^+=H \cup \psi^\le$
and $G^+=G \cup \{e\}$. Note that
$X^{e,J'}_{E,G}(\Phi,L) = 
\sum_{\phi^+ \in X_{E^+,G^+}(\Phi,L)} U(J')_{\phi^+(e)}$.
As $L$ is $(c,h)$-typical in $\Phi$ we have
$X_{E^+,G^+}(\Phi,L) < 1.1 d(L)^{|G|+1} n^{v_E}$, 
so by the claim $\mb{E}X^{e,J'}_{E,G}(\Phi,L) 
<  0.9 q! 2^{2r} \oO^{-1}\tT d(L)^{|G|} n^{v_E}$.
Any choice of $\phi^O_\tT$ affects $X^{e,J'}_{E,G}(\Phi,L)$
by $O(n^{v_E-1})$, so by Lemma \ref{dom} whp $X^{e,J'}_{E,G}(\Phi,L) 
< q! 2^{2r} \oO^{-1} \tT d(L)^{|G|} n^{v_E}$.
Thus $J'$ is $(q! 2^{2r} \oO^{-1} \tT, h)$-bounded 
wrt $(\Phi,L)$. \qed

\medskip

We also need the following estimate for the expected deviation from the mean 
of a random variable that is a sum of independent indicator variables.

\begin{lemma} \label{diff-binom}
There is $C>0$ such that for any sum of independent indicator variables
$X$ with mean $\mu$ we have $\mb{E}|X-\mu| \le C\sqrt{\mu}$.
\end{lemma}

\nib{Proof.}
We can assume $\mu>1$, otherwise we use the bound
$\mb{E}|X-\mu| \le 2\mu \le 2\sqrt{\mu}$.
Write $\mb{E}|X-\mu| = \sum_{t \ge 0} |t-\mu| \mb{P}(X=t) = E_0+E_1$,
where $E_i$ is the sum of $|t-\mu| \mb{P}(X=t)$ 
over $|t-\mu| > \tfrac{1}{2}C\sqrt{\mu}$ for $i=0$
or $|t-\mu| \le \tfrac{1}{2}C\sqrt{\mu}$ for $i=1$.
Clearly, $E_1 \le \tfrac{1}{2}C\sqrt{\mu}$,
and by Chernoff bounds, for $C$ large,
\begin{equation}
E_0 \le \sum_{a>\tfrac{1}{2}C\sqrt{\mu}} a(e^{-a^2/2\mu}+e^{-a^2/2(\mu+a/3)}) 
\le \tfrac{1}{2}C\sqrt{\mu}. \tag*{$\Box$}
\end{equation} 

Now we give the proof of our first key lemma,
on approximate integral decompositions.

\medskip

\nib{Proof of Lemma \ref{bddint:approx}.}
Suppose $J \in (\mb{Z}^D)^{\Phi_r}$ is $\tT$-bounded.
By Lemma \ref{bddint:flat}
there is some $J^0 \in (\mb{Z}^D)^{\Phi_r}$
and $\Psi^0 \in \mb{Z}^{\mc{A}(\Phi)}$ 
such that $\pl \Psi^0 = J-J^0$,
$J^0$ and $\Psi^0$ are $q^q \oO^{-1} \tT$-bounded,
and $U(J^0)_O < n^{0.1}$ for all $O \in \Phi_r/\Ss$.

Let $L$ be $\nu$-random in $\Phi^\circ_r$,
where $\nu = n^{-(3hq)^{-r}}$. By Lemma \ref{extrandom}
whp $L$ is $(n^{-1/3},h)$-typical in $\Phi$.
and by Lemma \ref{Jrandom}
whp $J^0$ is $(2q^q \oO^{-1} \tT, h)$-bounded wrt $(\Phi,L)$.

By Lemma \ref{bddint:random}
there is some $J^1 \in (\mb{Z}^D)^{\Phi[L]_r}$
and $\Psi^1 \in \mb{Z}^{\mc{A}(\Phi)}$ 
such that $\pl \Psi^1 = J^0-J^1$,
$J^1$ and $\Psi^1$ are $(2q)^{2q} \oO^{-2} \tT$-bounded, 
and $J^1$ is $((2q)^{2q} \oO^{-2} \tT, h)$-bounded wrt $(\Phi,L)$.

By Lemma \ref{bddrat:L} there is some $\oO_q^{-1} \tT/2$-bounded
$\Psi^* \in \mb{Q}^{\mc{A}(\Phi[L])}$ with $\pl \Psi^* = J^1$.

We obtain $\Psi^2 \in \mb{Z}^{\mc{A}(\Phi[L])}$ from $\Psi^*$ 
by randomised rounding as follows.
For each $\phi$ with $\Psi^*_\phi \ne 0$, let $s_\phi \in \pm 1$ 
be the sign of $\Psi^*_\phi$, let $m_\phi = \bfl{s_\phi \Psi^*_\phi}$, 
let $X_\phi$ be independent Bernoulli random variables 
such that $\Psi^*_\phi = s_\phi(m_\phi + \mb{E}X_\phi)$,
and let $\Psi^2_\phi = s_\phi(m_\phi + X_\phi)$.
Note that each $\mb{E}\Psi^2_\phi = \Psi^*_\phi$, 
so $\mb{E} \pl \Psi^2 = \pl \Psi^* = J^1$.

Let $J'=J^1-\pl \Psi^2$. For any $O \in \Phi[L]_r/\Ss$ we have 
$(J')^O = \sum_\phi \gG(\phi)^O s_\phi (\mb{E}X_\phi-X_\phi) $.
Separating the positive and negative contributions
for each $\gG$-atom $a$ at $O$ we can write
$U(J')_O \le \sum_{a \in \pm \gG[O]} |Y_a-\mu_a|$,
where each $Y_a$ is a sum of independent indicator
variables with mean $\mu_a \le U(\Psi^*)_O$.
By Lemma \ref{diff-binom} we have 
$\mb{E}U(J')_O < C \sqrt{U(\Psi^*)_O}$,
where $C$ depends only on $q$, $D$ and $K$.
For any $f \in \Phi_{r-1}$, whp $|L(Im(f))|<1.1 \nu n$,
so writing $\sum'_O$ for the sum over $O \in \Phi_r/\Ss$
with $f \sub Im(O) \in L$, by Cauchy-Schwartz
\begin{align*}
\mb{E}U(J')_f & < C \sum'_O \sqrt{U(\Psi^*)_O}
\le C (1.1\nu n \sum'_O U(\Psi^*)_O)^{1/2} \\
& < C (\nu q^r \oO_q^{-1} \tT)^{1/2} n < \oO_q^{3qh} \tT n/2,
\end{align*}
as $\nu = n^{-(3hq)^{-r}}$, $\tT > n^{-(4hq)^{-r}}$,
$\oO_q=\oO^{(9q)^{q+5}}$, $\oO > n^{-h^{-2q}}$.

Any rounding decision affects $U(J')_f$ by at most $1$,
so by Lemma \ref{dom} whp $U(J')_f < \oO_q^{3qh} \tT n$
for all $f \in \Phi_{r-1}$. 
Similarly, whp $U(\Psi^2)_f < 0.9 \oO_q^{-1} \tT n$
for all $f \in \Phi_{r-1}$.
Thus $J'$ is $\oO_q^{3qh} \tT$-bounded and 
$\Psi = \Psi^0 + \Psi^1 + \Psi^2$ is
$\oO_q^{-1} \tT$-bounded with $\pl \Psi = J - J'$. \qed

\subsection{Lifts and neighbourhood lattices}

This subsection contains some preliminaries for
the proof of our second key lemma 
in the following subsection.
Throughout we suppose 
$\mc{A}$ is a $\Ss^\le$-family,
$\gG \in (\mb{Z}^D)^{\mc{A}_r}$ and
$\Phi$ is a $\Ss$-adapted $[q]$-complex.
The construction in the following definition
is a technical device for working with
neighbourhood lattices.

\begin{defn} \label{def:lift} (lifts)

Fix a representative $\psi^O$
for each orbit $O \in \Phi_r/\Ss$.
Given $J \in \GG^{\Phi_r}$ we define
$J^\ua \in (\GG^{\Ss^\le})^{\Phi_r}$,
where for each $O \in \Phi_r/\Ss$ we let
$(J^\ua_{\psi^O})_\sS = J_{\psi^O \sS}$ where defined,
and set all other entries to zero.
We call $J^\ua$ a lift of $J$.
For $J' \in (\GG^{\Ss^\le})^{\Phi_r}$ we write 
$\pi(J') := \sum \{ (J'_\psi)_\sS \{\psi\sS\} : 
\psi \in \Phi_r, \sS \in \Ss^\le \}$.
\end{defn}

Note that $J = \pi(J^\ua)$.

\begin{defn} (lifted vector systems)
Let $\mc{A}^\ua$ be obtained from $\mc{A}$
by including a copy $A^{\ov{\tT}}$ of $A$
for each $A \in \mc{A}$ and 
$\ov{\tT} = (\tT^B: B \in Q)$
with each $\tT^B \in \Ss^B$.

Suppose $\gG \in \GG^{\mc{A}_r}$.
Let $\gG^\ua \in (\GG^{\Ss^\le})^{\mc{A}^\ua_r}$ 
where for $\tT \in A^{\ov{\tT}}$
we have $\gG^\ua_\tT=0$ unless
$\tT=\tT^B$ for some $B \in Q$, when
$(\gG^\ua_{\tT^B})_\sS$ is $\gG_{\tT^B \sS}$
if $\tT^B \in \Ss_{B'}$, $\sS \in \Ss^{B'}$ 
or $0$ otherwise.
\end{defn}

\begin{lemma} \label{lift:props} $ $
\begin{enumerate}
\item The set of $\gG^\ua$-molecules
is the set of lifts of $\gG$-molecules.
\item The set of $\gG^\ua$-atoms
is the set of lifts of $\gG$-atoms.
\item If $J \in \bgen{\gG(\Phi)}$ then 
any lift $J^\ua$ of $J$ is in $\bgen{\gG^\ua(\Phi)}$,
with $U(J^\ua)_O = U(J)_O$ for all $O \in \Phi_r/\Ss$.
\end{enumerate}
\end{lemma}

\nib{Proof.}
Clearly (ii) and (iii) follow from (i).
For (i), consider any $\phi \in A(\Phi)$
and $\phi^\ua \in A^{\ov{\tT}}(\Phi)$ with $\phi^\ua=\phi$.
We claim $\gG^\ua(\phi^\ua)$ is a lift of $\gG(\phi)$.
Indeed, for any $O \in \Phi_r/\Ss$ with $O \sub \phi\Ss$
and $\psi \in O$ we have $\gG^\ua(\phi^\ua)_\psi=0$,
except for $\psi = \phi \tT^B \in O$ 
where $B = \phi^{-1}(Im(O))$,
when each $(\gG^\ua(\phi^\ua)_\psi)_\sS 
= (\gG^\ua(\phi^\ua)_{\phi^\ua \tT^B})_\sS
= (\gG^\ua_{\tT^B})_\sS = \gG_{\tT^B \sS}
= \gG(\phi)_{\psi \sS}$.
Conversely, any lift of $\gG(\phi)$
with respect to orbit representatives $\psi^O$
can be expressed as $\gG^\ua(\phi^\ua)$
with $\phi^\ua \in A^{\ov{\tT}}(\Phi)$
where $\tT^B = \phi^{-1} \psi^O$
for each $B \in Q$, $Im(O) = \phi(B)$. \qed

\medskip

We also require some notation and basic properties
of neighbourhood lattices (recall the notation
of Definitions \ref{def:nhood} and \ref{nhood:perm}).

\begin{defn} \label{def:quotient} (quotients)
Let $\Ss^*=\Ss/B^*$ where $B^* \sub [q]$
with $r^*=r-|B^*|>0$.
Suppose $\mc{A}^*$ is a $\Ss^*$-family
that includes a copy $A^{\tT^*}$ of $(\Ss^*)^\le$ 
for each $A \in \mc{A}$ and $\tT^* \in A_{B^*}$.
We call $\gG^* \in (\mb{Z}^D)^{\mc{A}^*_{r^*}}$
a $B^*$-quotient of $\gG$ if for each
$\tT^* \in A_{B^*}$ there is 
$\sS^*=\sS^*(\tT^*) \in \Ss$ with $\sS^* \tT^* = id_{B^*}$
such that $\gG^*_{\tT'} = \gG_\tT$ whenever 
$\tT^* \sub \tT \in A_r$,
$\tT' = (\sS^* \tT)/id_{B^*} \in A^{\tT^*}$.
\end{defn}

\begin{lemma} \label{nhood:span}
Let $J \in (\mb{Z}^D)^{\Phi_r}$, $\psi^* \in \Phi_{B^*}$, 
$\Phi^*=\Phi/\psi^*$ and $J^*=J/\psi^*$. 
Suppose $\gG^* \in (\mb{Z}^D)^{\mc{A}^*_{r^*}}$ is
a $B^*$-quotient of $\gG \in (\mb{Z}^D)^{\mc{A}_r}$. Then
\begin{enumerate}
\item $\phi' \in \mc{A}^*(\Phi^*)$ iff 
$\phi'=(\phi \circ (\sS^*)^{-1})/\psi^*$ for some 
$\phi \in A(\Phi)$, $\tT^* \in A_{B^*}$, $\sS^*=\sS^*(\tT^*)$.
\item for such $\phi'$, $\phi$ we have
$\gG(\phi)/\psi^* = \gG^*(\phi')$,
\item if $J \in \bgen{\gG(\Phi)}$ then
$J^* \in \bgen{\gG^*(\Phi^*)}$.
\end{enumerate}
\end{lemma}

\nib{Proof.}
For (i), first consider any $\phi' \in \mc{A}^*(\Phi^*)$,
say $\phi' \in A^{\tT^*}(\Phi^*)$, and define $\phi$ 
such that $\phi'=(\phi \circ (\sS^*)^{-1})/\psi^*$.
As $\phi' \in \Ss^{*\le}(\Phi^*) = \Phi^*_{q-|B^*|}$
we have $\phi \circ (\sS^*)^{-1} \in \Phi$,
so $\phi \in \Phi_q = A(\Phi)$. Conversely, consider any 
$\phi'=(\phi \circ (\sS^*)^{-1})/\psi^*$ with 
$\phi \in A(\Phi)$, $\tT^* \in A_{B^*}$, $\sS^*=\sS^*(\tT^*)$.
For any $\tT' = (\sS^* \tT)/id_{B^*} \in A^{\tT^*}$
where $\tT^* \sub \tT \in A$ we have
$\phi\tT \in \Phi$ as $\phi \in A(\Phi)$,
so $\phi'\tT' = (\phi\tT)/\psi^* \in \Phi^*$,
so $\phi' \in A^{\tT^*}(\Phi^*)$.
This proves (i). Furthermore, if $\tT \in A_r$ 
then $(\gG(\phi)/\psi^*)_{\phi'\tT'}
= \gG(\phi)_{\phi\tT} = \gG_\tT = \gG^*_{\tT'} 
= \gG^*(\phi')_{\phi'\tT'}$, so (ii) holds, 
and (iii) is immediate from (ii). \qed

\subsection{Reducing support}

In this section we prove Lemma \ref{bddint:reduce},
using the inductive hypothesis of Lemma \ref{bddint} 
if $r>1$. Our argument will also prove the case $r=1$,
which is the base of the induction.
For convenient notation we rename $\oO'$ as $\oO$.
Let $V' \sub V(\Phi)$ with $|V'|=n/2$
be such that $(\Phi,V')$ is $(\oO,h)$-extendable wrt $V'$.
Suppose $J \in \bgen{\gG(\Phi)}$ 
is $\tT$-bounded, where $\tT < \oO_q^{3qh}$.
We need to find some $\oO_q^{-2qh} \tT$-bounded
$J' \in (\mb{Z}^D)^{\Phi[V']_r}$
and $\oO_q^{-2qh} \tT$-bounded $\Psi \in \mb{Z}^{\mc{A}(\Phi)}$
with $\pl \Psi = J-J'$.

We will define $J=J^0,\dots,J^r \in \bgen{\gG(\Phi)}$ 
so that $J^j_\psi=0$ whenever $|Im(\psi) \cap V'|<j$
and $J^j$ is $\tT_j$-bounded, where
$\tT_0=\tT$, $\tT_1 = 2^q \oO^{-1} \tT$
and for $0 < j < r$ we let
$\tT_{j+1} = 2^{r^2 + 2} \eta^{-r+1} M^{IH}_j \tT_j$,
where $\eta = (9q)^{-2q} \oO$
and $M^{IH}_j=\oO_{q-r+j}^{-3h}
= \oO^{-3h(9(q-r+j))^{q-r+j+5}}$.
As $\tT < \oO_q^{3qh}$ we see that 
$n^{-(5h(q-r+j))^{-j}} < n^{-(5hq)^{-r}}
< \tT_j < \oO_{q-r+j}^{3h(q-r+j)}$,
i.e.\ $\tT_j$ satisfies the necessary
bounds to apply Lemma \ref{bddint:avoid}
with $(q-r+j,j)$ in place of $(q,r)$,
and also that $\tT_r < \oO_q < 2^{-r}\eta$.

We start with $J^0=J$. To define $J^1=J-\pl^\gG \Psi^0$, 
for each orbit $O \in \Phi_r/\Ss$ with $Im(O) \cap V' = \es$, 
we fix a representative, say $\psi^O \in \Phi_{B^O}$, 
recall $f_{B^O}(J)_{\psi^O} \in \gG^{B^O}$ by Lemma \ref{gB},  
and find $n^O \in \mb{Z}^{\mc{A}_{B^O}}$ with $|n^O|= U(J)_O$
and $f_{B^O}(J)_{\psi^O} = \sum_\tT n^O_\tT \gG^\tT$.
Then $J^O = \sum_\tT n^O_\tT \gG[\psi^O]^\tT 
= \sum_\tT n^O_\tT \gG(\psi^O\tT^{-1})$. 

Let $S$ be the intset where each $(O,\tT)$
appears $|n^O_\tT|$ times with the sign of $n^O_\tT$.
For each $(O,\tT)$ in $S$ we add to $\Psi^0$ a uniformly random 
$\phi$ with $\psi^O\tT^{-1} \sub \phi \in A(\Phi)$
and $Im(\phi) \sm Im(\psi^O) \sub V'$,
with the same sign as that of $(O,\tT)$ in $S$.
Then $\gG(\psi^O\tT^{-1})$ is cancelled by $\gG(\phi)$
in $J^1=J-\pl^\gG \Psi^0$, and all other $\psi$
with $\gG(\phi)_\psi \ne 0$ have $Im(\psi) \cap V' \ne \es$.
As $(\Phi,V')$ is $(\oO,h)$-extendable wrt $V'$
there are at least $\oO (n/2)^{q-r}$ choices
for each $\phi$, so similarly to Lemma \ref{bddint:flat},
whp $J^1$ and $\Psi^0$ are $\tT_1$-bounded.
If $r=1$ this completes the construction,
so henceforth we suppose $r>1$.

Given $J^j$ with $0<j<r$ we will let $J^{j+1}=J^j-\pl^\gG\Psi^j$,
where $|Im(\phi) \cap V'|=j$ whenever $\Psi^j_\phi \ne 0$.
To define $\Psi^j$, we fix any (see Definition \ref{def:lift})
lift $J^{j\ua}$ of $J^j$
with orbit representatives $\psi^O$ for $O \in \Phi_r/\Ss$,
such that writing $B^O = \{i: \psi^O(i) \in V \sm V'\}$,
we have $\psi^{O'}\mid_{B^O} = \psi^O\mid_{B^O}$
for any orbit $O'$ that contains some $\psi$
with $\psi\mid_{B^O} = \psi^O\mid_{B^O}$.
Then for each $\psi^* \in \Phi_{r-j}$ with $Im(\psi^*) \cap V'=\es$
such that $\psi^*$ is some $\psi^O\mid_{B^O}$
we consider $\Phi^*=\Phi/\psi^*$, 
$J^* = J^{j\ua}/\psi^*$, $\Ss^*=\Ss/B^O$
and some $B^O$-quotient (see Definition \ref{def:quotient})
$\gG^* \in ((\mb{Z}^D)^{\Ss^\le})^{\mc{A}^*_{r^*}}$ of $\gG^\ua$.
Note that $J^*$ is supported in $\Phi^*[V']_{r^*}$,
and $\Phi^*[V']$ is $(\oO,h)$-extendable by Lemma \ref{nhood:ext}.
Then $J^* \in \bgen{\gG^*(\Phi^*)}$ by 
Lemmas \ref{lift:props} and \ref{nhood:span},
so $J^* \in \bgen{\gG^*(\Phi^*[V'])}$ by Lemma \ref{lattice}
(as in the earlier proof of Lemma \ref{bddint} modulo lemmas).

We will write $J^* = \pl^{\gG^*} \Psi^*$ for some 
$\Psi^* = \Psi^{\psi^*} \in \mb{Z}^{\mc{A}^*(\Phi^*[V'])}$
and define $\Psi^j$ as the sum over all such $\psi^*$
and $\phi' \in \mc{A}^*(\Phi^*[V'])$ of
$\Psi^{\psi^*}_{\phi'}\{ \phi \}$, where
$\phi'=(\phi \circ (\sS^*)^{-1})/\psi^*$ for some 
$\phi \in A^{\ov{\tT}}(\Phi)$, $\tT^* \in A^{\ov{\tT}}_{B^O}$, 
$\sS^*=\sS^*(\tT^*)$; recall from Lemma \ref{nhood:span} that any
$\phi' \in \mc{A}^*(\Phi^*)$ can be thus expressed, 
and then $\gG^\ua(\phi)/\psi^* = \gG^*(\phi')$.
Then $J^{j+1} = J^j-\pl^\gG\Psi^j$
will be as required if $\pl^\gG \Psi^j_\psi = J^j_\psi$
for every $\psi \in \Phi_r$ with $|Im(\psi) \cap V'|=j$.

To see this, consider any $O \in \Phi_r/\Ss$ with $|B^O|=r-j$ 
and $\psi = \psi^O \sS \in O$. Let $\psi^*=\psi^O\mid_{B^O}$,
and define $J^*$ and $\gG^*$ for $\psi^*$ as above.
For each $\psi' \in O$ with $\psi'\mid_{B^O}=\psi^*$
let $\sS'=\sS'(\psi')$ be such that $\psi = \psi' \sS'$.
Then \begin{align*} \pl^\gG \Psi^j_\psi 
& = \pi( \pl^{\gG^\ua} \Psi^j )_\psi 
= \sum_{\psi'} (\pl^{\gG^\ua} \Psi^j_{\psi'})_{\sS'}
= \sum_{\psi'} \sum_\phi \Psi^j_\phi (\gG^\ua(\phi)_{\psi'})_{\sS'} \\
& = \sum_{\psi'} \sum_\phi \Psi^j_\phi 
((\gG^\ua(\phi)/\psi^*)_{\psi'/\psi^*})_{\sS'} 
= \sum_{\psi'} \sum_{\phi' \in  \mc{A}^*(\Phi^*)} 
\Psi^{\psi^*}_{\phi'}  (\gG^*(\phi')_{\psi'/\psi^*})_{\sS'} \\
& = \sum_{\psi'} (\pl^{\gG^*} \Psi^{\psi^*}_{\psi'/\psi^*})_{\sS'}
= \sum_{\psi'} (J^*_{\psi'/\psi^*})_{\sS'}
= (J^*_{\psi^O/\psi^*})_\sS
= (J^{j\ua}_{\psi^O})_\sS = J^j_\psi, \ \text{ as required.}
\end{align*}

We will construct $\Psi^*$ for each $\psi^*$
as above sequentially using Lemma \ref{bddint:avoid}
applied to $J^* \in \bgen{\gG^*(\Phi^*[V'])}$
which is valid by Lemma \ref{avoid/bddint} and
the inductive hypothesis of Lemma \ref{bddint},
as $\mc{A}^*$ is a $(\Ss^*)^\le$-family,
$\Phi^*[V']$ is $\Ss^*$-adapted (by Lemma \ref{nhood:adapt})
and $(\oO,h)$-extendable (by Lemma \ref{nhood:ext}), 
and $J^*$ is $\tT_j$-bounded, as this is true of $J^j$
and so $J^{j\ua}$ by Lemma \ref{lift:props}.
Lemma \ref{bddint:avoid} 
will give $\Psi^*$ that is $M^{IH}_j \tT_j$-bounded,
and also provides the option
to avoid using any $\eta$-bounded sets 
$B^*_p \sub \Phi^*[V']_p$ for $j \le p \le q-(r-j)$,
which we will define below so as to
maintain boundedness throughout the algorithm.
During the construction of $\Psi^j$, 
we say that $\psi \in \Phi$ is full if 
\begin{enumerate}
\item $|Im(\psi)|=r$ with $U(\Psi^j)_\psi > 0$, 
\item $|Im(\psi)|=r-1$ with $U(\Psi^j)_\psi 
  > \tT_{j+1} n/4 - C_0-1$, or
\item $|Im(\psi)|=i<r-1$ with more than $2^{-r} \eta n - C_0-1$
  full elements of $\Phi_{i+1}\mid_\psi$.
\end{enumerate}

There will be no uses of full sets by $\Psi^j$ 
apart from at most $U(J^j)_\psi$ `forced' uses of $\psi$ 
for each $\psi \in \Phi_r$, 
and at most $C_0+1$ further unforced uses.
This implies that for $\psi$ with $|Im(\psi)|=i<r$,
if $\psi$ is not full then $\Phi_{i+1}\mid_\psi$
has at most $2^{-r} \eta n$ full elements
(if $i=r-1$ we use $\tT_r<2^{-r}\eta$ here).

We claim that there 
is no full $\psi' \in \Phi$ with $|Im(\psi') \cap V'|=j-1$. 
Indeed, suppose we have such $\psi'$. 
Let $\psi^a = \psi'[V']$ and $\psi^b=\psi' \sm \psi^a$.
Then \begin{align*}
& ( \tT_{j+1} n/4 - C_0-1) (2^{-r}\eta n - C_0-1)^{r-1-|Im(\psi')|} \\
& < \sum \{ U(\Psi^j)_\psi: \psi' \sub \psi \in \Phi_{r-1} \} 
= \sum \{ U(\Psi^{\psi^*})_{\psi^a} : \psi^b \sub \psi^* \} \\
& < n^{r-j-|\psi^b|} \cdot M^{IH}_j\tT_j n 
= M^{IH}_j\tT_j n^{r-|Im(\psi')|},
\end{align*}
as all $\Psi^{\psi^*}$ are $M^{IH}_j\tT_j$-bounded.
This contradicts the definition of $\tT_{j+1}$
and so proves the claim. 

When applying Lemma \ref{bddint:avoid}
to $J^* \in \bgen{\gG^*(\Phi^*[V'])}$ as above,
for $j \le p \le q-r+j$
we let $B^*_p$ be the set of $Im(\ups)$
with $\ups \in \Phi^*_p[V']$
such that $\ups \cup \psi^*{}'$ is full 
for some $\psi^*{}' \sub \psi^*$, 
and $\ups' \cup \psi^*{}'$ is not full 
for any $\ups' \subn \ups$, $\psi^*{}' \sub \psi^*$.
Then $|B^*_p(Im(\ups))| < \eta n$ 
for all $\ups \in \Phi^*_{p-1}[V']$,
by definition for $p>j$, and by the claim for $p=j$. 
Thus we can apply Lemma \ref{bddint:avoid} to obtain
some  $M^{IH}_j\tT_j$-bounded
$\Psi^* = \Psi^{\psi^*} \in \mb{Z}^{\mc{A}^*(\Phi^*[V'])}$ 
with $\pl^{\gG^*} \Psi^* = J^*$ such that
\begin{enumerate}
\item if $p>j$ then $U(\Psi^*)_\psi \le 1$
for all $\psi \in \Phi_p[V']$
and $U(\Psi^*)_\psi = 0$ if $Im(\psi) \in B^*_p$,
\item $U(\Psi^*)_\psi \le U(J^*)_\psi + C_0+1$ 
for all $\psi \in \Phi_j[V']$,
and $U(\Psi^*)_\psi = U(J^*)_\psi$ if $Im(\psi) \in B^*_j$.
\end{enumerate}

As described above, 
$\Psi^j$ is the sum over all such $\psi^*$
and $\phi' \in \mc{A}^*(\Phi^*[V'])$ of
$\Psi^{\psi^*}_{\phi'}\{ \phi \}$, where
$\phi'=(\phi \circ (\sS^*)^{-1})/\psi^*$ for some 
$\phi \in A^{\ov{\tT}}(\Phi)$, $\tT^* \in A^{\ov{\tT}}_{B^O}$, 
$\sS^*=\sS^*(\tT^*)$.
We claim that $\Psi^j$ is $\tT_{j+1}/2$-bounded.
To see this, we fix $\psi \in \Phi_{r-1}$,
let $\psi^a=\psi[V']$, $\psi^b=\psi \sm \psi^a$
and consider cases according to $p=|Im(\psi^a)|$.
\begin{enumerate}
\item If $p=j-1$ then $U(\Psi^j)_\psi = U(\Psi^{\psi^b})_{\psi^a} 
< M^{IH}_j\tT_j n < \tT_{j+1} n/2$,
as $\Psi^{\psi^b}$ is $M^{IH}_j\tT_j$-bounded,
\item If $p>j$ then 
$U(\Psi^j)_\psi \le \tT_{j+1} n/4 - C_0-1$ 
before $\psi$ is full, after which 
$U(\Psi^{\psi^*})_{\psi^a} = 0$ whenever $\psi^b \sub \psi^*$,
except for at most one such $\psi^*$ with 
$U(\Psi^{\psi^*})_{\psi^a} \le C_0+1$, so $U(\Psi^j)_\psi
< \tT_{j+1} n/4 \le \tT_{j+1} n/2$.
\item If $p=j$ then 
$U(\Psi^j)_\psi \le \tT_{j+1} n/4 - C_0-1$ 
before $\psi$ is full, after which 
$U(\Psi^{\psi^*})_{\psi^a} = U(J^j)_{\psi^a \cup \psi^*}$ 
whenever $\psi^b \sub \psi^*$,
except for at most one such $\psi^*$ with 
$U(\Psi^{\psi^*})_{\psi^a} \le U(J^j)_{\psi^a \cup \psi^*} + C_0+1$,
so $U(\Psi^j)_\psi \le \tT_{j+1} n/4 
+ U(J^j)_\psi \le \tT_{j+1} n/2$, as $U(J^j)_\psi \le \tT_j n$.
\end{enumerate}
Thus $U(\Psi^j)_\psi \le \tT_{j+1} n/2$
in all cases, so the claim holds.

It follows that $J^{j+1}=J^j-\pl\Psi^j$ is $\tT_{j+1}$-bounded,
so the construction can proceed to the next step.
We conclude with $J' := J^r \in (\mb{Z}^D)^{\Phi[V']_r}$ 
and $\Psi = \sum_j \Psi^j \in \mb{Z}^{\mc{A}(\Phi)}$,
such that $J'$ and $\Psi$ are $\oO_q^{-2qh}\tT$-bounded
with $\pl \Psi = J-J'$. \qed

\section{Applications} \label{sec:app}

In this section we give several applications 
of our main theorem, including the theorems 
stated in the introduction of the paper.
Most of the applications will follow from 
a decomposition theorem for hypergraphs
in various partite settings.
We also give some results on 
coloured hypergraph decomposition,
and a simple illustration
(the Tryst Table Problem) of other applications
that are not equivalent to hypergraph decomposition,
but for the sake of brevity we leave a detailed study
of these applications for future research.

Our first theorem in this section can be
viewed as a simplified form of Theorem \ref{main},
in which various general definitions 
are specialised to the setting of hypergraph decompositions.
To state it we require two definitions.

\begin{defn}
Let $\Phi$ be a $[q]$-complex 
and $H$ be an $r$-graph on $[q]$.
We say $G \in \mb{N}^{\Phi^\circ_r}$ 
is $(H,c,\oO)$-regular in $\Phi$ if there are
$y_\phi \in [\oO n^{r-q},\oO^{-1} n^{r-q}]$ 
for each $\phi \in \Phi_q$ with $\phi(H) \sub G$
so that $\sum_\phi y_\phi \phi(H) = (1 \pm c)G$.
\end{defn}

\begin{defn}
We say that an $R$-complex $\Phi$ is exactly $\Ss$-adapted 
if whenever $\phi \in \Phi_B$ and $\tau \in Bij(B',B)$ we have
$\phi \circ \tau\in \Phi_{B'}$ iff $\sS \in \Ss^B_{B'}$.
We say $\Phi$ is exactly adapted if 
$\Phi$ is exactly $\Ss$-adapted for some $\Ss$.
\end{defn}
 
\begin{theo} \label{Hdecomp:ext}
Let $H$ be an $r$-graph on $[q]$ and $\Phi$ be 
an $(\oO,h)$-extendable exactly adapted $[q]$-complex
where $n=|V(\Phi)|>n_0(q)$ is large,
$h=2^{50q^3}$, $\dD = 2^{-10^3 q^5}$, 
$n^{-\dD}<\oO<\oO_0(q)$ is small and $c=\oO^{h^{20}}$.
Suppose $G \in \bgen{H(\Phi)}$ is $(H,c,\oO)$-regular in $\Phi$
and $(\Phi,G)$ is $(\oO,h)$-extendable.
Then $G$ has an $H$-decomposition in $\Phi_q$.
\end{theo}

\nib{Proof.}
Suppose $\Phi$ is exactly $\Ss$-adapted,
let $\mc{A}=\{A\}$ with $A=\Ss^\le$,
and $\gG \in \{0,1\}^{A_r}$ with
each $\gG_\tT = 1_{Im(\tT) \in H}$.
Let $G^* \in \mb{N}^{\Phi_r}$ with
$G^*_\psi = G_{Im(\psi)}$ for $\psi \in \Phi_r$.
For any $\phi \in A(\Phi) = \Phi_q$ 
and $\tT \in A_r$ we have
$\gG(\phi)_{\phi\tT} = \gG_\tT = 1_{Im(\tT) \in H}$.
As $\Phi$ is exactly $\Ss$-adapted, we deduce
$G \in \bgen{H(\Phi)}$ iff $G^* \in \bgen{\gG(\Phi)}$,
and that an $H$-decomposition of $G$ is equivalent 
to a $\gG(\Phi)$-decomposition of $G^*$.
 
There are two types in $\gG$ for each $B \in [q]_r$:
the edge type $\{\tT \in A_B: Im(\tT) \in H\}$ and
the nonedge type $\{\tT \in A_B: Im(\tT) \notin H\}$.
Each $\gG^\tT$ is the all-1 vector for $\tT$ in an edge
type or the all-0 vector for $\tT$ in a nonedge type,
so $\gG$ is elementary.
The atom decomposition of $G^*$ is
$G^* = \sum_{e \in \Phi^\circ_r} G_e e^*$, where $e^*_\psi=1$ 
for all $\psi \in \Phi_r$ with $Im(\psi)=e$,
i.e.\ $e^*$ contains all edge types at $e$.

As $G$ is $(H,c,\oO)$-regular in $\Phi$, we have
$\sum_\phi y_\phi \phi(H) = (1 \pm c)G$ for some
$y_\phi \in [\oO n^{r-q},\oO^{-1} n^{r-q}]$ 
for each $\phi \in \Phi_q$ with $\phi(H) \sub G$.
For any such $\phi$ we have $\gG(\phi) \le_\gG G$,
so $\phi \in \mc{A}(\Phi,G)$. Also, for any
$B \in [q]_r$ and $\psi \in \Phi_B$,
writing $1_B \in T_B$ for the edge type
we have $\pl^{1_B} y_\psi 
= \sum_{\phi: t_\phi(\psi)=1_B} y_\phi
= \sum \{ y_\phi: Im(\psi) \in \phi(H) \}
= (1 \pm c)(G^*)^{1_B}_\psi$,
so $G^*$ is $(\gG,c,\oO)$-regular.

To apply Theorem \ref{main}, it remains to show that 
$(\Phi,\gG[G])$ is $(\oO,h)$-extendable.
We have $\gG[G]=(\gG[G]_B: B \in Q)$
where if $B \notin H$ then $\gG[G]_B = \Phi_B$
and if $B \in H$ then 
$\gG[G]_B = \{ \psi \in \Phi_B : G_{Im(\psi)}>0 \}$.
Let $E=(J,F,\phi)$ be any $\Phi$-extension of rank $s$
and $J' \sub J_r \sm J[F]$. 
As $(\Phi,G)$ is $(\oO,h)$-extendable
we have $X_{E,J'}(\Phi,G) > \oO n^{v_E}$.
Consider any $\phi^+ \in X_{E,J'}(\Phi,G)$.
For any $\psi \in J'$ we have $\phi^+\psi \in \Phi$ 
and $Im(\phi^+\psi) \in G$,
so $\phi^+\psi \in \gG[G]$.
Thus $\phi^+ \in X_{E,J'}(\Phi,\gG[G])$,
so $(\Phi,\gG[G])$ is $(\oO,h)$-extendable.
Now $G^*$ has a $\gG(\Phi)$-decomposition,
so $G$ has an $H$-decomposition. \qed

\medskip

The following theorem solves the
$H$-decomposition problem in the nonpartite setting
(so is similar in spirit to \cite{GKLO2});
it is an immediate corollary of 
Theorem \ref{Hdecomp:ext} with $\Ss=S_q$
and Theorem \ref{Hlattice:nonpartite}.

\begin{theo} \label{Hdecomp}
Let $H$ be an $r$-graph on $[q]$ and $\Phi$ be 
an $(\oO,h)$-extendable $S_q$-adapted $[q]$-complex
where $n=|V(\Phi)|>n_0(q)$ is large,
$h=2^{50q^3}$, $\dD = 2^{-10^3 q^5}$, 
$n^{-\dD}<\oO<\oO_0(q)$ is small and $c=\oO^{h^{20}}$.
Suppose $G$ is $H$-divisible and $(H,c,\oO)$-regular in $\Phi$
and $(\Phi,G)$ is $(\oO,h)$-extendable.
Then $G$ has an $H$-decomposition in $\Phi_q$.
\end{theo}

In the introduction we stated a simplified form
of this result (using typicality rather than
extendability and regularity); we now give the deduction.
(It will also follow from a later more general result,
but we include the proof for the sake of exposition.)

\medskip

\nib{Proof of Theorem \ref{Hdecomp:typ}.} 
Let $H$ be an $r$-graph on $[q]$ and $G$ be an $H$-divisible
$(c,h^q)$-typical $r$-graph, where $n=|V(\Phi)|>n_0(q)$,
$h=2^{50q^3}$, $\dD = 2^{-10^3 q^5}$, $d(G) > 2n^{-\dD/h^q}$
and $c < c_0 d(G)^{h^{30q}}$. We need to show 
that $G$ has an $H$-decomposition.
Let $\Phi$ be the complete $[q]$-complex on $V(G)$,
i.e.\ each $\Phi_B = Inj(B,V(G))$.
Then $\Phi$ is exactly $S_q$-adapted.

We claim that $(\Phi,G)$ is 
$(\oO,h)$-extendable with $\oO = \frac{1}{2} d(G)^{h^q}$.
To see this, consider any $\Phi$-extension 
$E = (J,F,\phi)$ with $J \sub [q](h)$
and any $J' \sub J^\circ_r \sm J^\circ[F]$.
Write $V(J) \sm F = \{x_1,\dots,x_{v_E}\}$
and suppose for $i \in [v_E]$ that there
are $m_i$ edges of $J'$ that
use $x_i$ but no $x_j$ with $j>i$.
The number of choices for the embedding
of each $x_i$ given any previous choices 
is $(1 \pm m_i c)d(G)^{m_i} n$.
Thus $X_{E,J'}(\Phi,G) > \tfrac{1}{2} d(G)^{e_E} n^{v_E}
> \oO n^{v_E}$, as claimed.

Furthermore, if $e \in G$, $f \in H$,
$J=\ova{q}$, $F=f$, $\psi \in Bij(f,e)$ 
we see that there are
$(1 \pm |H|c) d(G)^{|H|-1} n^{q-r}$
extensions of $\psi$ to $\phi \in \Phi_q$ 
with $e \in \phi(H) \sub G$.
Defining $y_\phi = d(G)^{1-|H|} n^{r-q}$
for all $\phi \in \Phi_q$ with $\phi(H) \sub G$
we see that $G$ is $(H,|H|c,\oO)$-regular in $\Phi$.
The theorem now follows from Theorem \ref{Hdecomp}. \qed

\medskip

Our next definition sets up notation for hypergraph 
decompositions in a generalised partite setting
that incorporates several earlier examples in the paper.
It is followed by a theorem that solves 
the corresponding hypergraph decomposition problem.

\begin{defn} \label{HPblowup}
Let $H$ be an $r$-graph on $[q]$ and 
$\mc{P}=(P_1,\dots,P_t)$ be a partition of $[q]$.
Let $\Ss$ be the group of all $\sS \in S_q$ 
with all $\sS(P_i)=P_i$.
Let $\Phi$ be an exactly $\Ss$-adapted $[q]$-complex
with parts $\mc{Q}=(Q_1,\dots,Q_t)$, where each 
$Q_i = \{ \psi(j): j \in P_i, \psi \in \Phi_j \}$.
Let $G \in \mb{N}^{\Phi^\circ_r}$.

For $S \sub [q]$ the $\mc{P}$-index of $S$ is 
$i_{\mc{P}}(S) = (|S \cap P_1|,\dots,|S \cap P_t|)$;
similarly, we define the $\mc{Q}$-index of subsets
of $V(\Phi)$, and also refer to both as the `index'.
For $i \in \mb{N}^t$ we let $H_i$ and $G_i$ be the
(multi)sets of edges in $H$ and $G$ of index $i$.
Let $I = \{i: H_i \ne \es \}$.
We call $G$ an $(H,\mc{P})$-blowup 
if $G_i \ne \es \Ra i \in I$. 

We say $G$ has a $\mc{P}$-partite $H$-decomposition if 
it has an $H$-decomposition using copies $\phi(H)$ 
of $H$ with all $\phi(P_i) \sub Q_i$.

For $e \sub V(\Phi)$ we define 
the degree vector $G_I(e) \in \mb{N}^I$
by $G_I(e)_i=|G_i(e)|$ for $i \in I$.
Similarly, for $f \sub [q]$ we define
$H_I(f)$ by $H_I(f)_i=|H_i(f)|$.
For $i' \in \mb{N}^t$ let $H^I_{i'}$
be the subgroup of $\mb{N}^I$ generated 
by $\{ H_I(f): i_{\mc{P}}(f) = i' \}$.
We say $G$ is $(H,\mc{P})$-divisible 
if $G_I(e) \in H^I_{i'}$ whenever $i_{\mc{P}}(e)=i'$.

Let $G^* \in \mb{N}^{\Phi_r}$ with $G^*_\psi = G_e$ 
for $\psi \in \Phi_r$, $e=Im(\psi)$ with $i_{\mc{Q}}(e) \in I$,
and $G^*_\psi$ is otherwise undefined.
\end{defn}

\begin{theo} \label{HPdecomp}
With notation as in Definition \ref{HPblowup},
suppose $n/h \le |Q_i| \le n$ with $n > n_0(q)$ large,
$G$ is an $(H,\mc{P})$-divisible $(H,\mc{P})$-blowup,
$G$ is $(H,c,\oO)$-regular in $\Phi$,
and $(\Phi,G^*)$ is $(\oO,h)$-extendable,
where $h=2^{50q^3}$, $\dD = 2^{-10^3 q^5}$, 
$n^{-\dD}<\oO<\oO_0(q)$ is small and $c=\oO^{h^{20}}$.
Then $G$ has a $\mc{P}$-partite $H$-decomposition.
\end{theo}

\nib{Proof.}
Let $A=\Ss^\le$, $H^* = \{ \tT \in A_r: Im(\tT) \in H\}$
and $\gG_\tT = 1_{\tT \in H^*}$ for $\tT \in A_r$, 
Then an (integral) $\mc{P}$-partite $H$-decomposition of $G$ 
is equivalent to an (integral) $\gG(\Phi)$-decomposition of $G^*$.
By Theorem \ref{Hdecomp:ext}, it remains to show
$G \in \bgen{H(\Phi)}$, i.e.\ $G^* \in \bgen{\gG(\Phi)}
= \mc{L}_\gG(\Phi)$ (by Lemma \ref{lattice}).

Consider any $i \in I = \{i: H_i \ne \es \}$
and $i' \in \mb{N}^t$ with all $i'_j \le i_j$.
Let $m^i_{i'} = \prod_{j \in [t]} (i_j-i'_j)!$.
For any $B' \sub B \in Q$ 
with $i_{\mc{P}}(B')=i'$ and  $i_{\mc{P}}(B)=i$
and $\psi' \in \Phi_{B'}$ with $Im(\psi')=e$ 
we have $((G^*)^\sharp_{\psi'})_B) = 
\sum \{ G^*_\psi: \psi' \sub \psi \in \Phi_B \} = m^i_{i'} |G_i(e)|$.
Writing $O=\psi'\Ss$, for any $\psi \in O$ we have 
$((G^*)^\sharp_\psi)_B) = m^i_{i'} |G_i(e)|$.
Thus we obtain $(G^*)^\sharp_\psi$ from $G_I(e)$
by copying coordinates and multiplying all copies
of each $i$-coordinate by $m^i_{i'}$.

Similarly, for any $\tT' \in A_{B'}$ with $Im(\tT')=f$ 
we have $(\gG^\sharp_{\tT'})_B) = 
\sum \{ \gG_\tT: \tT' \sub \tT \in A_B \} = m^i_{i'} |H_i(f)|$,
so $\sgen{\gG^\sharp[O]}$ is generated by vectors $v^f \in (\mb{Z}^Q)^O$
where $f \sub [q]$ with $i_{\mc{P}}(f)=i'$
and for each $\psi' \in O$, $B \in Q$ we have 
$(v^f_{\psi'})_B = m^i_{i'} |H_i(f)|$ where $i=i_{\mc{P}}(B)$.
Thus all vectors in $\sgen{\gG^\sharp[O]}$ are obtained
from vectors in $H^I_{i'}$ by the same transformation
that maps $G_I(e)$ to $((G^*)^\sharp)^O$.
As $G$ is $(H,\mc{P})$-divisible we deduce
$((G^*)^\sharp)^O \in \sgen{\gG^\sharp[O]}$ 
for any $O \in \Phi/\Ss$, as required. \qed

\medskip

Similarly to the nonpartite setting,
which also give a simplified form of the previous
theorem in which we replace extendability 
and regularity by typicality
(which we will generalise here to multigraphs).

\begin{defn} \label{HPblowup:typ}
With notation as in Definition \ref{HPblowup},
suppose that $\Phi$ is the $[q]$-complex where each $\Phi_B$ 
consists of all maps $\psi:B \to V(G)$ with 
$\psi(B \cap P_i) \sub Q_i$ for all $i \in [t]$;
we call $\Phi$ the complete $[q]$-complex wrt $(\mc{P},\mc{Q})$,
and note that $\Phi$ is exactly $\Ss$-adapted.
For $f \in Q$ let
$|G^* \cap \Phi_f| = \sum \{ G^*_\psi: \psi \in \Phi_f \}$ 
and $d_f(G^*) = |G^* \cap \Phi_f| / |\Phi_f|$.

For any $\Phi$-extension $E=(J,F,\phi)$ 
and $J' \sub J^\circ_r \sm J^\circ[F]$ let 
$X_{E,J'}(\Phi,G) = \sum_{\phi' \in X_E(\Phi)} 
\prod_{e \in J'} G_{\phi'(e)}$.

We call $G$ a $(c,s)$-typical $(H,\mc{P})$-blowup 
if for every $\Phi$-extension $E=(J,F,\phi)$ of rank $s$
and $J' \sub J^\circ_r \sm J^\circ[F]$ we have $X_{E,J'}(\Phi,G) 
= (1 \pm c) X_E(\Phi) \prod_{e \in J'} d_e(G^*)$,
where for $e \in J^\circ_f$ we write $d_e(G^*)=d_f(G^*)$.
\end{defn}

Now we give our theorem on decompositions of typical
multigraphs in the generalised partite setting.
We note that the case $\mc{P}=([q])$
implies Theorem \ref{Hdecomp:typ} 
and the case $\mc{P}=(\{1\},\dots,\{q\})$ 
implies Theorem \ref{Hdecomp:partite1}.

\begin{theo} \label{HPdecomp:typ}
Let $H$ be an $r$-graph on $[q]$ and 
$\mc{P}=(P_1,\dots,P_t)$ be a partition of $[q]$.
Suppose each $n/h \le |Q_i| \le n$
with $n > n_0(q)$ and $h=2^{50q^3}$,
that $\dD = 2^{-10^3 q^5}$, $d > 2n^{-\dD/h^q}$ 
and $c < c_0 d^{h^{30q}}$, where $c_0=c_0(q)$ is small.
Let $G$ be an $(H,\mc{P})$-divisible
$(c,h)$-typical $(H,\mc{P})$-blowup 
wrt $\mc{Q}=(Q_1,\dots,Q_t)$,
such that $d_f(G)>d$ for all $f \in H$
and $G_e < d^{-1}$ for all $e \in [n]_r$.
Then $G$ has a $\mc{P}$-partite $H$-decomposition.
\end{theo}

\nib{Proof.}
With notation as in Definition \ref{HPblowup},
it follows (as in the proof of Theorem \ref{Hdecomp:typ})
from the definition of $(c,h)$-typical $(H,\mc{P})$-blowup 
that $(\Phi,G^*)$ is $(\oO,h)$-extendable
with $\oO = \frac{1}{2} d^{h^q} > n^{-\dD}$. 
By Theorem \ref{HPdecomp} it remains to show
that $G$ is $(H,2c,\oO)$-regular in $\Phi$.

To see this, first note that 
as $G$ is $(H,\mc{P})$-divisible
we have $G_I(\es) \in H^I_0 = \sgen{H_I(\es)}$,
so there is some integer $Y$ such that
$|G_i| = Y |H_i|$ for all $i \in I$.
For each $\phi \in \Phi_q$ with $\phi(H) \sub G$
we let $y_\phi = YZ^{-1} \prod_{f \in H} G_{\phi(f)}$, 
where $Z = \prod_{j \in [t]} |Q_j|^{|P_j|}
\prod_{f \in H} d_f(G^*)$.
Then $y_\phi \in [\oO n^{r-q}, \oO^{-1} n^{r-q}]$
for each such $\phi$.

We need to show for any $e \in \Phi^\circ_r$ that
$\sum \{ y_\phi : e \in \phi(H) \} = (1 \pm 2c)G_e$.
We can suppose $G_e \ne 0$, so $i = i(e) \in I$.
Let $m_i = \prod_{j \in [t]} i_j!$.
For any $B \in H_i$ there are $m_i$ choices
of $\psi \in \Phi_B$ with $\psi(B)=e$.
It suffices to show for any such $B$ and $\psi$
that $\sum \{ y_\phi: \psi \sub \phi \}
= (1 \pm 2c)G_e / (m_i|H_i|)$.

Let $E=(\ova{q},B,\psi)$ and $J'= H \sm \{B\}$.
As $G$ is a $(c,h)$-typical $(H,\mc{P})$-blowup 
we have $X_{E,J'}(\Phi,G) 
= (1 \pm c) X_E(\Phi) \prod_{f \in J'} d_f(G^*)
= (1 \pm 2c) Z/(m_i |G_i|)$,
as $X_E(\Phi) = (1+O(n^{-1})) 
\prod_{j \in [t]} |Q_j|^{|P_j|-i_j}$
and $m_i |G_i| = |G^* \cap \Phi_B| = d_B(G^*) |\Phi_B|
= (1+O(n^{-1})) d_B(G^*) \prod_{j \in [t]} |Q_j|^{i_j}$.
Therefore 
\begin{equation} \tag*{$\Box$}
\sum \{ y_\phi: \psi \sub \phi \}
= YZ^{-1} \sum_{\phi \in X_E(\Phi)} \prod_{f \in H} G_{\phi(f)}
= G_e YZ^{-1} X_{E,J'}(\Phi,G) 
= (1 \pm 2c)G_e / (m_i|H_i|).
\end{equation}

\medskip

Now we will prove several other theorems stated 
in the introduction for which we gave an equivalent
reformulation in terms of hypergraph decompositions
in partite settings as above.
We start with the existence of resolvable designs,
or more generally, resolvable 
hypergraph decompositions of multigraphs.\footnote{
For the sake of brevity we just consider the case
that $H$ is vertex-regular, and leave the computation
of the lattice for general $H$ to the reader
(if $H$ is not vertex-regular than we need $G$
the difference between any two vertex degrees in $G$
to be divisible by the gcd of all differences
of vertex degrees in $H$).}
We deduce Theorem \ref{resolvable}
by applying Theorem \ref{resolve:G} 
with $H=K^r_q$ and $G = \lL K^r_n$.

\begin{theo} \label{resolve:G}
Let $H$ be a vertex-regular $r$-graph on $[q]$
and $G$ be a vertex-regular $H$-divisible 
$r$-multigraph on $[n]$ where
$n > n_0(q)$ is large and $q \mid n$.
Let $\Phi$ be a $S_q$-adapted $[q]$-complex on $V(G)$.
Suppose $G$ is $(H,c,\oO)$-regular in $\Phi$
and $(\Phi,G)$ is $(\oO,h)$-extendable, where 
$h=2^{50q^3}$, $\dD = 2^{-10^3 q^5}$,
$n^{-\dD/2h}<\oO<\oO_0$, $c=\oO^{h^{22}}$.
Then $G$ has a resolvable $H$-decomposition.
\end{theo}

\nib{Proof.}
We start by recalling the equivalent
partite hypergraph decomposition problem.
Let $Y$ be a set of $m$ vertices disjoint from $X$, 
where $m$ is the least integer with 
$\tbinom{m}{r-1} \ge q|G|/|H|n$.
Let $J$ be a random
$(r-1)$-graph on $Y$ with $|J|=q|G|/|H|n$.
Let $G'$ be the $r$-multigraph obtained from $G$
by adding as edges (with multiplicity one)
all $r$-sets of the form $f \cup \{x\}$
where $f \in J$ and $x \in X$. 
Let $H'$ be the $r$-graph whose vertex set is the disjoint
union of a $q$-set $A=[q]$ and an $(r-1)$-set $B$,
and whose edges consist of all $r$-sets in $A \cup B$
that are contained in $H$ or have exactly one vertex in $A$.
To adopt the notation of Definition \ref{HPblowup}
we let $\mc{P}=(P_1,P_2)$ with $P_1=A$, $P_2=B$
and $\mc{Q}=(Q_1,Q_2)$ with $Q_1=X$, $Q_2=Y$.
Then $G'$ is an $(H',\mc{P})$-blowup and we wish
to find a $\mc{P}$-partite $H'$-decomposition of $G'$.

First we check $(H',\mc{P})$-divisibility.
The set of edge indices is $I = \{ (r,0), (1,r-1) \}$.
We identify $\mb{N}^I$ with $\mb{N}^2$ 
by assigning $(r,0)$ to the first coordinate
and $(1,r-1)$ to the second. Let $i' \in \mb{N}^2$.
Suppose $i'_2>0$. Then $H'{}^I_{i'}$ is $\{(0,0)\}$
unless $i'_1 \le 1$ and $i'_2 \le r-1$,
in which case $H'{}^I_{i'}$ is generated by
$(0,1)$ if $i'_1=1$ or $(0,q)$ if $i'_1=0$.
The resulting $(H',\mc{P})$-divisibility conditions
($G_I(e) \in H'{}^I_{i'}$ whenever $i_{\mc{P}}(e)=i'$)
are satisfied trivially when $i'_1=1$
and as $q \mid n$ when $i'_1=0$.
Now suppose $i'_2=0$. Then $H'{}^I_{i'}$ is generated 
by all $(|H(f)|,0)$ with $f \in [q]_{i'_1}$ if $i'_1>1$
or by $(|H(f)|,1)$ with $f \in [q]_{i'_1}$ if $i'_1 \le 1$.
If $i'_1>1$ then the $(H',\mc{P})$-divisibility condition
is equivalent to the $H$-divisibility condition
that $gcd_{i'_1}(H)$ divides $|G(e)|$.
For $i'_1=0$ we need $(|G|,n|J|)$ to be an integer multiple
of $(|H|,q)$, so we require the $H$-divisibility condition
$|H| \mid |G|$ and also $|J| = q|G|/|H|n$.
For $i'_1=1$ we need $(|G(x)|,|J|)$ to be an integer multiple
of $(rq^{-1}|H|,1)$ for any $x \in X$
(recall that $H$ is vertex-regular),
so we require the $K^r_q$-divisibility condition
that $rq^{-1}|H|$ divides $|G(x)|$
and also that $G$ is vertex-regular.
Therefore $G'$ is $(H',\mc{P})$-divisible.

To apply Theorem \ref{HPdecomp}, it remains to check
extendability and regularity. We let $\Phi'$ 
be the $(A \cup B)$-complex where each
$\Phi'_{A' \cup B'}$ for $A' \sub A$, $B' \sub B$
consists of all $\phi \in Inj(A' \cup B',X \cup Y)$
with $\phi\mid_{A'} \in \Phi$ and $\phi(B') \sub Y$.
Consider any $\Phi'$-extension $E=(J,F,\phi)$ with 
$J \sub (A \cup B)(h)$ and any $J' \sub J_r \sm J[F]$
with $J'_{B'} \ne \es \Ra i_{\mc{P}}(B') \in I$.
Let $E(A)$ and $J'(A)$ be obtained 
by restricting to $A(h)$,
and define $E(B)$ similarly.
Then whp $X_{E,J'}(\Phi',G') = (1+o(1)) 
X_{E(A),J'(A)}(\Phi,G) m^{v_{E(B)}}$,
where $X_{E(A),J'(A)}(\Phi,G)
> \oO n^{v_{E(A)}}$ as $(\Phi,G)$ is $(\oO,h)$-extendable.
As $|G| > \oO n^r$ we have
$m \ge (q\oO/Q)^{1/(r-1)} n$,
so $X_{E,J'}(\Phi',G') > \oO^{2h} (n+m)^{v_E}$ (say),
i.e.\ $(\Phi',G')$ is $(\oO^{2h},h)$-extendable.

For regularity, as $G$ is $(H,c,\oO)$-regular in $\Phi$
we can choose $y_\phi \in [\oO n^{r-q}, \oO^{-1} n^{r-q}]$
for each $\phi \in \Phi_A$ with $\phi(H) \sub G$
with $\sum \{ y_\phi : e \in \phi(H) \} = (1 \pm c)G_e$
for any $e \in X_r$. We define
$y'_{\phi'} = y_{\phi'\mid_A} ((r-1)!|J|)^{-1}$ 
for each $\phi' \in \Phi'_{A \cup B}$ with $\phi'(B) \in J$.
Then $y'_{\phi'} \in [\oO^{2h} (n+m)^{1-q}, \oO^{-2h} (n+m)^{1-q}]$
for all $\phi' \in \Phi'_{A \cup B}$ with $\phi'(H') \sub G'$.
For any $e \in X_r$ we have
$\sum \{ y'_{\phi'} : e \in \phi'(H') \}
= \sum \{ y_\phi : e \in \phi(H) \}
= (1 \pm c) G'_e$.

Also, for any $e = f \cup \{x\}$
with $x \in X$ and $f \in J$, we have
\begin{align*}
\sum \{ y'_{\phi'} : e \in \phi'(H') \}
& = |J|^{-1} \sum \{ y_\phi: x \in Im(\phi) \} \\
& = (|J|rq^{-1}|H|)^{-1} \sum_{x \in e' \in X_r}  
\sum \{ y_\phi : e' \in \phi(H) \} \\
& = (|J|rq^{-1}|H|)^{-1} (1 \pm c) |G(x)|.
\end{align*}
As $G$ is vertex-regular,
$|G(x)| = r|G|/n = r|H||J|/q$,
so $\sum \{ y'_{\phi'} : e \in \phi'(H') \}
= 1 \pm c = (1 \pm c) G'_e$.
Therefore $G'$ is $(H',c,\oO^{2h})$-regular in $\Phi'$. \qed

\medskip

Next we show the existence of large sets of designs.
Again, we first consider the more general setting
of decompositions of multidesigns into designs.

\begin{theo} \label{large:G}
Let $\Phi$ be a $S_q$-adapted $[q]$-complex with 
$V(\Phi)=[n]$ where $n > n_0(q,\lL)$ is large.
Suppose $G \in \mb{N}^{\Phi^\circ_q}$ 
is an $r$-multidesign with all $G_e<\oO^{-1}$.
For all $0 \le i \le r$ suppose
$Z_i := \lL \tbinom{q-i}{r-i}^{-1} \tbinom{n-i}{r-i} \in \mb{Z}$
and $Z_i \mid |G(f)|$ for all $f \in [n]_i$.
Suppose also that $(\Phi,G)$ is $(\oO,h)$-extendable, where 
$h=2^{50q^3}$, $\dD = 2^{-10^3 q^5}$, $n^{-\dD/2h}<\oO<\oO_0(q,\lL)$. 
Then $G$ has a decomposition into $(n,q,r,\lL)$-designs.
\end{theo}

\nib{Proof.}
We start by recalling the equivalent
partite hypergraph decomposition problem.
Let $Y$ be a set of $m$ vertices disjoint from $X$,
where $m$ is the least integer with 
$\tbinom{m}{q-r} \ge |G| / Z_0$.
Let $J$ be a random $(q-r)$-graph on $Y$ with $|J|=|G|/Z_0$.
Let $G'$ be the $q$-multigraph obtained from $G$
by adding as edges with multiplicity $\lL$
all $q$-sets of the form $e \cup f$
with $e \sub X$ and $f \in J$.
Let $H$ be the $q$-graph whose vertex set is the disjoint
union of a $q$-set $A$ and a $(q-r)$-set $B$,
and whose edges consist of $A$ 
and all $q$-sets in $A \cup B$ that contain $B$.
To adopt the notation of Definition \ref{HPblowup}
we let $\mc{P}=(P_1,P_2)$ with $P_1=A$, $P_2=B$
and $\mc{Q}=(Q_1,Q_2)$ with $Q_1=X$, $Q_2=Y$.
Then $G'$ is an $(H,\mc{P})$-blowup and we wish
to find a $\mc{P}$-partite $H$-decomposition of $G'$.

First we check $(H,\mc{P})$-divisibility.
We identify $\mb{N}^I$ with $\mb{N}^2$ 
by assigning $(q,0)$ to the first coordinate
and $(r,q-r)$ to the second. Let $i' \in \mb{N}^2$.
Suppose $i'_2>0$. Then $H^I_{i'}$ is $\{(0,0)\}$
unless $i'_1 \le r$ and $i'_2 \le q-r$,
in which case $H^I_{i'}$ is generated 
by $(0, \tbinom{q-i'_1}{r-i'_1})$.
The corresponding $(H,\mc{P})$-divisibility 
condition is $Z_{i'_1} \in \mb{Z}$.
Now suppose $i'_2=0$. Then $H^I_{i'}$ is generated 
by $(1,0)$ if $i'_1>r$ or by 
$(1,\tbinom{q-i'_1}{r-i'_1})$ if $i'_1 \le r$.
The corresponding $(H,\mc{P})$-divisibility 
condition is trivial if $i'_1>r$.
If $i'_1 \le r$, for each $f \in [n]_{i'_1}$ we need
$(|G(f)|, \lL |J| \tbinom{n-i'_1}{r-i'_1} )$
to be an integer multiple of $(1,\tbinom{q-i'_1}{r-i'_1})$,
i.e.\ $|G(f)| =  |J| Z_{i'_1}$; this is equivalent 
to $G$ being an $r$-multidesign.
Therefore $G'$ is $(H,\mc{P})$-divisible.

To apply Theorem \ref{HPdecomp}, it remains to check
extendability and regularity. We let $\Phi'$ 
be the $(A \cup B)$-complex where each
$\Phi'_{A' \cup B'}$ for $A' \sub A$, $B' \sub B$
consists of all $\phi \in Inj(A' \cup B',X \cup Y)$
with $\phi\mid_{A'} \in \Phi$ and $\phi(B') \sub Y$.
Then $(\Phi',G')$ is $(\oO^{2h},h)$-extendable
as in the proof of Theorem \ref{resolve:G}.
For regularity, we define 
$y_\phi = G_{\phi(A)} (q!(q-r)!|J|)^{-1}$ 
for each $\phi \in \Phi'_{A \cup B}$ with $\phi(B) \in J$.
Then $y_\phi \in [\oO^{2h} (n+m)^{r-q}, \oO^{-2h} (n+m)^{r-q}]$
for all $\phi \in \Phi'_{A \cup B}$ with $\phi(H) \sub G'$.
For any $e \in X_q$  we have
$\sum \{ y_\phi : e \in \phi(H) \} =  G_e = G'_e$.
For any $e = f \cup f'$ with $f \in X_r$ and $f' \in J$,
as $G$ is an $r$-multidesign we have
$|G(f)| = Q \tbinom{n}{r}^{-1} |G| = \lL |G|/Z_0$, 
so $\sum \{ y_\phi : e \in \phi(H) \}
=  (q!(q-r)!|J|)^{-1} q! |G(f)| (q-r)!
= |G(f)| Z_0 / |G| = \lL = G'_e$. \qed

\medskip

Now we prove the existence conjecture for large sets of designs.

\medskip

\nib{Proof of Theorem \ref{large}.} 
First we note that the case that
$\lL$ is fixed and $n>n_0(q,\lL)$ is large
follows from Theorem \ref{large:G} applied to $G=K^q_n$. 
Now we can assume $\lL > \lL'(q)$ is large.
We let $\lL_0 = \prod_{i=1}^r \tbinom{q-i}{r-i}$ 
and write $\lL = \mu \lL_0 + \lL_1$ for some integers
$\mu, \lL_1$ with $0 \le \lL_1 < \lL_0$.
Write $Z_i := \lL \tbinom{q-i}{r-i}^{-1} \tbinom{n-i}{r-i}$
and note that by assumption all $Z_i \in \mb{Z}$.
Let $\ell = \tbinom{n}{q}/Z_0 = \lL^{-1} \tbinom{n-r}{q-r}$.
It suffices to decompose $K^q_n$ into
$\ell \mu$ designs with parameters $(n,q,r,\lL_0)$ 
and $\ell$ with parameters $(n,q,r,\lL_1)$.
Indeed, these can then be combined
into $\ell$ designs with parameters $(n,q,r,\lL)$.

We start by choosing the $\ell$ designs with parameters 
$(n,q,r,\lL_1)$. We do this by a greedy process,
where we start with $K^q_n$ and 
repeatedly delete some $(n,q,r,\lL_1)$-design.
Note that the divisibility conditions for
the existence of an $(n,q,r,\lL_1)$-design 
are satisfied, namely all
$\tbinom{q-i}{r-i} \mid \lL_1 \tbinom{n-i}{r-i}$.
At each step of the process we have some $q$-graph $G$.
We say that a $q$-set $e$ is full if $e \in K^q_n \sm G$,
and for $i<q$ that an $i$-set $f$ is full if
it is contained in at least $c(q) n$ full $(i+1)$-sets,
where we choose $1/\lL'(q) \ll c(q) \ll 1/h$.
Once a set is full we will avoid using it.

There can be no full $r$-set, as this would belong
to at least $(c(q) n)^{q-r}/(q-r)!$ full $q$-sets,
but we are only choosing $\ell \lL_1 
< \lL_0  \lL'(q)^{-1} \tbinom{n-r}{q-r}$ such $q$-sets.
Let $\Phi$ be the $[q]$-complex on $V(G)$
where each $\Phi_B$ consists of all 
$\phi \in Inj(B,V(G))$ such that 
all subsets of $Im(\phi)$ are not full.
Then $\Phi$ is $(1/2,h)$-extendable (say),
so by Theorem \ref{Hdecomp:ext} we can find
a $K^r_q$-decomposition of $\lL_1 K^r_n$ in $\Phi_q$,
i.e.\ an $(n,q,r,\lL_1)$-design avoiding full sets.
Thus the algorithm can be completed to choose 
$\ell$ designs with parameters $(n,q,r,\lL_1)$.

Finally, letting $\Phi$ and $G$ be as above
after the final step of the algorithm,
$(\Phi,G)$ is $(1/2,h)$-extendable,
$G$ is an $r$-multidesign,
$Z^0_i := \lL_0 \tbinom{q-i}{r-i}^{-1} 
\tbinom{n-i}{r-i} \in \mb{Z}$
and $Z^0_i \mid |G(f)|$ for all $f \in [n]_i$.
By Theorem \ref{large:G} we can decompose $G$ into $\ell \mu$ 
designs with parameters $(n,q,r,\lL_0)$, as required. \qed

\medskip

Next we prove the existence of complete resolutions.

\medskip

\nib{Proof of Theorem \ref{resolution}.}
Suppose $q$ is fixed and $n>n_0(q)$ is large
with $n = q$ mod $lcm([q])$. We start by recalling
the reformulation of complete resolution as
a partite hypergraph decomposition problem.
We consider disjoint sets of vertices $X$ and $Y$ where $|X|=n$ 
and $Y$ is partitioned into $Y_j$, $2 \le j \le q+1$ 
with $|Y_j| = \tfrac{n-j+2}{q-j+2}$. We let $G'$ be the $q$-graph 
whose edges are all $q$-sets $e \sub X \cup Y$ 
such that $|e \cap Y_j| \le 1$ for all $2 \le j \le q+1$,
and if $e \cap Y_j \ne\es$ then $e \cap Y_i \ne\es$ 
for all $i>j$. Let $H$ be the $q$-graph whose vertex set 
is the disjoint union of two $q$-sets $A$ and 
$B = \{b_2,\dots,b_{q+1}\}$, whose edges are all $q$-sets 
$e \sub A \cup B$ such that if $b_j \in e$ 
then $b_i \in e$ for all $i>j$. 
To adopt the notation of Definition \ref{HPblowup}
we let $\mc{P}=(P_1,\dots,P_{q+1})$
and $\mc{Q}=(Q_1,\dots,Q_{q+1})$ 
with $P_1=A$, $Q_1=X$ and
$P_j=\{b_j\}$, $Q_j=Y_j$ for $2 \le j \le q+1$.
Then $G'$ is an $(H,\mc{P})$-blowup and we wish
to find a $\mc{P}$-partite $H$-decomposition of $G'$.

To apply Theorem \ref{HPdecomp}, we consider the
complete $(A \cup B)$-complex $\Phi$ wrt $(\mc{P},\mc{Q})$.
Then $(\Phi,G')$ is clearly $(1/2,h)$-extendable (say).
Also, every edge of $G'$ is in exactly $\tbinom{n}{q}$
copies of $H$, so $G'$ is $(H,c,\oO)$-regular in $\Phi$
for any $c>0$ and $\oO < \oO_0$.
It remains to check $(H,\mc{P})$-divisibility.

The set of index vectors of edges is
$I = \{ i^j : 1 \le j \in q+1 \} \sub \mb{N}^{q+1}$
where $i^j_{j'}$ is $1$ for $j+1 \le j' \le q+1$,
$i^j_1 = j-1$ and $i^j_{j'}=0$ otherwise.
We identify $\mb{N}^I$ with $\mb{N}^{q+1}$
by assigning $i^j$ to coordinate $j$.
Consider $i' \in \mb{N}^{q+1}$.
We can assume $i'_{j'} \le 1$ for $j'>1$.
If there is any $j'>1$ with $i'_{j'} \ne 0$,
we let $j^0$ be the least such $j'$, 
otherwise we let $j^0=q+2$.
We can assume $i'_1 \le j^0-2$,
otherwise $H^I_{i'} = 0$.
Then $H^I_{i'}$ is generated by $v^{i'} \in \mb{N}^q$
where each $v^{i'}_j = 1_{i'_1+1 \le j \le j^0-1} 
 \tbinom{q-i'_1}{j-1-i'_1}$. 

Let $u^{i'}  \in \mb{N}^q$ where each $u^{i'}_j$ is 
$1_{i'_1+1 \le j \le j^0-1} \tbinom{n-i'_1}{j-1-i'_1} 
\prod_{j^*=j+1}^{q+1} |Y_{j^*}|^{1-i'_{j^*}}$.
Then the corresponding $(H,\mc{P})$-divisibility 
condition is that $u^{i'}$ is an integer multiple of $v^{i'}$.
It suffices to consider the case that $i'_j=1$
for all $j \ge j^0$, and so each $u^{i'}_j$ is 
$1_{i'_1+1 \le j \le j^0-1} \tbinom{n-i'_1}{j-1-i'_1} 
\prod_{j^*=j+1}^{j^0-1} |Y_{j^*}|$. 
Then for $i'_1+1 \le j \le j^0-1$ we have
$u^{i'}_j/v^{i'}_j 
= \tbinom{q-i'_1}{j-1-i'_1}^{-1}
 \tbinom{n-i'_1}{j-1-i'_1} 
\prod_{j^*=j+1}^{j^0-1} |Y_{j^*}|
= \prod_{j^*=i'_1+2}^{j^0-1} 
\tfrac{n-j^*+2}{q-j^*+2}$.
This is an integer constant independent of $j$, as required. \qed

\medskip

Next we solve the Tryst Table Problem.

\medskip

\nib{Proof of Theorem \ref{tryst}.}
Let $\Phi$ be the complete $[9]$-complex 
on an $n$-set $V$ where $n$ is large.
Let $G^* \in (\mb{Z}^2)^{\Phi_3}$ with all $G^*_\phi = (1,1)$.
Let $\mc{A}=\{A\}$ with $A=S_9^\le$.
Let $\gG \in (\mb{Z}^2)^{A_3}$ where
\begin{itemize}
\item $\gG_\tT = (1,0)$ if $Im(\tT)=\{1,4,7\}$,
\item $\gG_\tT = (0,1)$ if $Im(\tT)=\{3i-2,3i-1,3i\}$
for some $i \in [3]$ and $\tT(\min(Dom(\tT)))=3i-2$,
\item $\gG_\tT = (0,0)$ otherwise.
\end{itemize}
The Tryst Table Problem is equivalent
to finding a $\gG(\Phi)$-decomposition of $G^*$.
 
There are three types in $\gG$ for each $B \in [9]_3$,
where the type of $\tT$ is determined by $\gG_\tT$ as above,
so $\gG$ is elementary. The atom decomposition of $G^*$ 
is $G^* = \sum_{e=abc \in [n]_3} (e^1+e^a+e^b+e^c)$, 
where $e^1_\psi$ is $(1,0)$ for all $\psi$ with $Im(\psi)=e$
otherwise $0$, and each $e^x_\psi$ for $x \in e$
is $(0,1)$ for all $\psi$ with $Im(\psi)=e$
and $\psi(\min(Dom(\psi)))=x$, otherwise $0$.
(The interpretation of $e^1$ is that $e$ is the set of captains,
and of $e^x$ is that $e$ is a team with captain $x$.)
As all nonzero $\gG$-atoms at $e$ appear in $G^*$
we have $\gG[G]=(\gG[G]_B: B \in [9]_3)$
with each $\gG[G]_B = \Phi_B$, so $(\Phi,\gG[G])$ 
is $(\oO,h)$-extendable for any $\oO<1$ and $n>n_0(\oO,h)$.
To show regularity of $G^*$ we let
$y_\phi = 1/6(n-3)_6$ for all $\phi \in A(\Phi)=\Phi_9$.
Then for any $B \in [9]_3$, $\psi \in \Phi_B$, $t \in T_B$ we have 
$\pl^t y_\phi = 6(n-3)_6 |\{\phi: t_\phi(\psi)=t\}|=1$,
where the factor of $6$ either represents 
all bijections from $\{1,4,7\}$ to $e=Im(\psi)$,
or all bijections from $\{3i-2,3i-1,3i\}$ to $e$
mapping $3i-2$ to some $x \in e$,
where $x$ is fixed and $i$ ranges over $[3]$.
Therefore $G^*$ is $(\gG,c,1/7)$-regular in $\Phi$ for any $c>0$.

To apply Theorem \ref{main}, 
it remains to show that $G^* \in \bgen{\gG(\Phi)}
= \mc{L}_\gG(\Phi)$ (by Lemma \ref{lattice}).
Fix any $O \in \Phi/S_9$,
write $e=Im(O) \in \Phi^\circ$ and $i=|e|$.
Then $((G^*)^\sharp)^O \in (\mb{Z}^2)^{[9]_3 \times O}$
is a vector supported on the coordinates
$(B,\psi')$ with $B' \sub B \in [9]_3$ 
and $\psi' \in O \cap \Phi_{B'}$
in which every nonzero coordinate is equal:
we have $((G^*)^\sharp_{\psi'})_B) = 
\sum \{ G^*_\psi: \psi' \sub \psi \in \Phi_B \}  
= (n-i)_{3-i}  (1,1) $.

We need to show $((G^*)^\sharp)^O \in \sgen{\gG^\sharp[O]}$.
First consider the case $i=3$.
Then it is clear that $((G^*)^\sharp)^O$
is the sum of the $\gG^\sharp$-atoms at $O$,
as these are obtained from the $\gG$-atoms
$e^1$, $e^a$, $e^b$, $e^c$ described above by identifying 
each $\psi$ with $(B,\psi)$ where $\psi \in \Phi_B$. 
 
Now suppose $i=2$, say $Im(O)=e=ab$. 
There are four
$\gG^\sharp$-atoms at $O$, which we label
$e^{ab} = \gG^\sharp(1 \to a, 4 \to b)$
($a$ and $b$ are captains),
$e^a = \gG^\sharp(1 \to a, 2 \to b)$
($a$ captains a team containing $b$),
$e^b = \gG^\sharp(1 \to b, 2 \to a)$
($b$ captains a team containing $a$),
$e^0 = \gG^\sharp(2 \to a, 3 \to b)$
($a$ and $b$ are in the same team,
neither is the captain).

To calculate $e^{ab}$, consider 
any $\tT' \in A_2$
with $\tT'(x)=1$, $\tT'(y)=4$
and $\psi' \in \Phi_2$ with
$\psi'(x)=a$, $\psi'(y)=b$.
Then $e^{ab}_{\psi'} = \gG^\sharp_{\tT'}$,
so each $(e^{ab}_{\psi'})_{xyz}  = 
\sum \{ \gG_\tT: \tT' \sub \tT \in A_{xyz} \}
= \gG_{x \to 1, y \to 4, z \to 7} = (1,0)$.

Next, if $\tT' \in A_2$
with $\tT'(x)=1$, $\tT'(y)=2$
and $\psi' \in \Phi_2$ with
$\psi'(x)=a$, $\psi'(y)=b$
then each $(e^a_{\psi'})_{xyz}  = 
\sum \{ \gG_\tT: \tT' \sub \tT \in A_{xyz} \}
= \gG_{x \to 1, y \to 2, z \to 3}$
is $(0,1)$ if $x = \min\{x,y,z\}$
or $(0,0)$ otherwise. 
Similarly, $(e^b_{\psi'})_{xyz}$ is $(0,1)$ 
if $y = \min\{x,y,z\}$ or $(0,0)$ otherwise.

Finally, if $\tT' \in A_2$
with $\tT'(x)=2$, $\tT'(y)=3$
and $\psi' \in \Phi_2$ with
$\psi'(x)=a$, $\psi'(y)=b$
then each $(e^a_{\psi'})_{xyz}  = 
\sum \{ \gG_\tT: \tT' \sub \tT \in A_{xyz} \}
= \gG_{z \to 1, x \to 2, y \to 3}$
is $(0,1)$ if $z = \min\{x,y,z\}$
or $(0,0)$ otherwise. 

Therefore $((e^{ab}+e^a+e^b+e^0)_{\psi'})_{xyz} = (1,1)$
for every $\psi' \in O$ and $xyz \in [9]_3$,
so $((G^*)^\sharp)^O \in \sgen{\gG^\sharp[O]}$.

Now suppose $i=1$, say $Im(O)=a$. 
There are two
$\gG^\sharp$-atoms at $O$, which we label
$a^1 = \gG^\sharp(1 \to a)$ ($a$ is a captain),
$a^0 = \gG^\sharp(2 \to a)$ ($a$ is not a captain).

Consider any $\tT' \in A_1$ with $\tT'(x)=1$
and $\psi' \in \Phi_1$ with $\psi'(x)=a$.
Then each $(a^1_{\psi'})_{xyz}  = 
\sum \{ \gG_\tT: \tT' \sub \tT \in A_{xyz} \}$
is $(2,2)$ if $x=\min\{x,y,z\}$
or $(2,0)$ otherwise.

Next consider any $\tT' \in A_1$ with $\tT'(x)=2$
and $\psi' \in \Phi_1$ with $\psi'(x)=a$.
Then each $(a^0_{\psi'})_{xyz}  = 
\sum \{ \gG_\tT: \tT' \sub \tT \in A_{xyz} \}$
is $(0,0)$ if $x=\min\{x,y,z\}$
or $(0,2)$ otherwise.

Therefore $((a^1+a^0)_{\psi'})_{xyz} = (2,2)$
for every $\psi' \in O$ and $xyz \in [9]_3$.
As each $((G^*)^\sharp_{\psi'})_{xyz} = (n(n-1),n(n-1))$
we have $((G^*)^\sharp)^O \in \sgen{\gG^\sharp[O]}$.

Finally, $\gG^\sharp[\es]$ is generated by a vector
$v$ with all $(v_\es)_B = (6,6)$.
As $((G^*)^\sharp_\es)_B = (n(n-1)(n-2),n(n-1)(n-2))$
we have $(G^*)^\sharp_\es \in \gG^\sharp[\es]$. \qed

\medskip

Now we consider the more general setting of
coloured hypergraph decompositions.
We require some definitions.

\begin{defn}
Suppose $H$ is an $r$-graph on $[q]$,
edge-coloured as $H = \cup_{d \in [D]} H^d$.
We identify $H$ with a vector $H \in (\mb{N}^D)^Q$,
where each $(H_f)_d = 1_{f \in H^d}$.

Let $\Phi$ be a $[q]$-complex.
For $\phi \in \Phi_q$ we define 
$\phi(H) \in (\mb{N}^D)^{\Phi^\circ_r}$
by $\phi(H)_{\phi(f)} = H_f$.
Let $\mc{H}$ be an family of 
$[D]$-edge-coloured $r$-graphs on $[q]$.
Let $\mc{H}(\Phi) = \{ \phi(H):
 \phi \in \Phi_q, H \in \mc{H} \}$.

We say $G \in (\mb{N}^D)^{\Phi^\circ_r}$ 
is $(\mc{H},c,\oO)$-regular in $\Phi$ if there are
$y^H_\phi \in [\oO n^{r-q},\oO^{-1} n^{r-q}]$ for each 
$H \in \mc{H}$, $\phi \in \Phi_q$ with $\phi(H) \le G$
so that $\sum \{ y^H_\phi \phi(H) \} = (1 \pm c)G$.

We say that $(\Phi,G)$ is $(\oO,h)$-extendable
if $(\Phi,G')$ is $(\oO,h)$-extendable,
where $G' = (G^1,\dots,G^D)$ with each 
$G^d = \{e \in \Phi^\circ_r: (G_e)_d> 0 \}$.
\end{defn}

The following generalises Theorem \ref{Hdecomp:ext}
by allowing colours and also families of hypergraphs.

\begin{theo} \label{colHdecomp:ext}
Let $\mc{H}$ be an family of 
$[D]$-edge-coloured $r$-graphs on $[q]$.
Let $\Phi$ be an $(\oO,h)$-extendable 
exactly adapted $[q]$-complex
where $n=|V(\Phi)|>n_0(q,D)$ is large,
$h=2^{50q^3}$, $\dD = 2^{-10^3 q^5}$, 
$n^{-\dD}<\oO<\oO_0(q,D)$ is small and $c=\oO^{h^{20}}$.
Suppose $G \in \bgen{\mc{H}(\Phi)}$ 
is $(\mc{H},c,\oO)$-regular in $\Phi$
and $(\Phi,G)$ is $(\oO,h)$-extendable.
Then $G$ has an $H$-decomposition in $\Phi_q$.
\end{theo}

\nib{Proof.}
We follow the proof of Theorem \ref{Hdecomp:ext},
with appropriate modifications for the more general setting.
Suppose $\Phi$ is exactly $\Ss$-adapted and let 
$\mc{A}=\{A^H: H \in \mc{H}\}$ with each $A^H=\Ss^\le$.
Let $\gG \in (\mb{Z}^D)^{\mc{A}_r}$ with $\gG_\tT = e_d$
if $\tT \in A^H_r$, $H \in \mc{H}$, $d \in [D]$
with $Im(\tT) \in H^d$ or $\gG_\tT=0$ otherwise.
Let $G^* \in (\mb{N}^D)^{\Phi_r}$ with
$G^*_\psi = G_{Im(\psi)}$ for $\psi \in \Phi_r$.
Then $G \in \bgen{\mc{H}(\Phi)}$ iff $G^* \in \bgen{\gG(\Phi)}$,
and an $\mc{H}$-decomposition of $G$ is equivalent 
to a $\gG(\Phi)$-decomposition of $G^*$.

There are $D+1$ types in $\gG$ for each $B \in [q]_r$:
the colour $d$ type 
$\{\tT \in A^H_B: Im(\tT) \in H^d, H \in \mc{H} \}$
for each $d \in [D]$, and the nonedge type 
$\{\tT \in A^H_B: Im(\tT) \notin H \in \mc{H} \}$.
Each $\gG^\tT$ is $e_d$ in all coordinates
for $\tT$ in a colour $d$ type
or $0$ in all coordinates
for $\tT$ in a nonedge type,
so $\gG$ is elementary.
The atom decomposition of $G^*$ is
$G^* = \sum_{f \in \Phi^\circ_r} \sum_{d \in [D]} 
(G_f)_d f^d$, where $f^d_\psi=e_d$ 
for all $\psi \in \Phi_r$ with $Im(\psi)=f$,
i.e.\ $f^d$ contains all colour $d$ types at $f$.

As $G$ is $(\mc{H},c,\oO)$-regular in $\Phi$ we have
$\sum \{ y^H_\phi \phi(H) \} = (1 \pm c)G$ for some
$y^H_\phi \in [\oO n^{r-q},\oO^{-1} n^{r-q}]$ for each 
$H \in \mc{H}$, $\phi \in \Phi_q$ with $\phi(H) \le G$.
For any such $\phi \in H(\Phi)$
we have $\gG(\phi) \le_\gG G$, so $\phi \in \mc{A}(\Phi,G)$. 
We define $y_\phi = y^H_\phi$ for $\phi \in A^H(\Phi)$.
Then for any $B \in [q]_r$, $\psi \in \Phi_B$ and $d \in [D]$,
writing $t_d \in T_B$ for the colour $d$ type we have 
$\pl^{t_d} y_\psi = \sum_{\phi: t_\phi(\psi)=t_d} y_\phi
= \sum \{ y^H_\phi: Im(\psi) \in \phi(H^d), H \in \mc{H} \}
= (1 \pm c)(G^*)^{t_d}_\psi$,
so $G^*$ is $(\gG,c,\oO)$-regular.

To apply Theorem \ref{main}, it remains to show that 
each $(\Phi,\gG[G]^H)$ is $(\oO,h)$-extendable.
If $B \notin H$ then $\gG[G]^H_B = \Phi_B$
and if $B \in H^d$ for $d \in [D]$ then 
$\gG[G]^H_B = \{ \psi \in \Phi_B : Im(\psi) \in G^d \}$.
Let $E=(J,F,\phi)$ be any $\Phi$-extension of rank $s$
and $J' \sub J_r \sm J[F]$.
Let $J'' = (J^d: d \in [D])$ where each
$J^d = \bigcup \{ J'_B: B \in H^d \}$.
As $(\Phi,G)$ is $(\oO,h)$-extendable
we have $X_{E,J''}(\Phi,G) > \oO n^{v_E}$.
Consider any $\phi^+ \in X_{E,J''}(\Phi,G)$.
For any $\psi \in J^d$ we have $\phi^+\psi \in \Phi$ 
and $Im(\phi^+\psi) \in G^d$, so $\phi^+\psi \in \gG[G]^H$.
Thus $\phi^+ \in X_{E,J'}(\Phi,\gG[G]^H)$,
so $(\Phi,\gG[G]^H)$ is $(\oO,h)$-extendable.
Now $G^*$ has a $\gG(\Phi)$-decomposition,
so $G$ has an $\mc{H}$-decomposition. \qed

\medskip

We conclude by applying Theorem \ref{colHdecomp:ext}
to the two results on rainbow clique decompositions
stated in the introduction.

\medskip

\nib{Proof of Theorem \ref{rainbow:all}.}
We apply Theorem \ref{colHdecomp:ext}
with $G=[\tbinom{q}{r}]K^r_n$ and 
$\mc{H}$ equal to the set of all
rainbow $[\tbinom{q}{r}]$-colourings of $K^r_q$.
We let $\Phi$ be the complete $[q]$-complex on $[n]$
and note that $G$ is $(\mc{H},c,\oO)$-regular in $\Phi$
and $(\Phi,G)$ is $(\oO,h)$-extendable
for any $c>0$ and some $\oO=\oO(q)$.

It remains to check $G \in \bgen{\mc{H}(\Phi)}$.
Let $G^*$ and $\gG$ be as in the proof
of Theorem \ref{colHdecomp:ext}.
We need to show $G^* \in \bgen{\gG(\Phi)} 
= \mc{L}_\gG(\Phi)$ (by Lemma \ref{lattice}),
i.e.\ $((G^*)^\sharp)^O \in \sgen{\gG^\sharp[O]}$
for any $O \in \Phi/S_q$. 

Fix any $O \in \Phi/S_q$,
write $e=Im(O) \in \Phi^\circ$ and $i=|e|$.
Then $((G^*)^\sharp)^O \in ((\mb{Z}^Q)^Q)^O 
= (\mb{Z}^Q)^{Q \times O}$
is a vector supported on the coordinates
$(B,\psi')$ with $B' \sub B \in Q$ and $\psi' \in O \cap \Phi_{B'}$
with each $((G^*)^\sharp_{\psi'})_B) = 
\sum \{ G^*_\psi: \psi' \sub \psi \in \Phi_B \} \in \mb{Z}^Q$
equal to $(r-i)! \tbinom{n-i}{r-i}$ in each coordinate.
Also, $\sgen{\gG^\sharp[O]}$ is generated by $\gG^\sharp$-atoms
$\gG^\sharp(\ups)$ at $O$, each of which is
supported on the same coordinates $(B,\psi')$
as $((G^*)^\sharp)^O$, with each
$(\gG^\sharp(\ups)_{\psi'})_B)$ equal 
to some $(r-i)! v$ in each coordinate,
where $v \in \{0,1\}^Q$ is any vector
with $\sum_B v_B = \tbinom{q-i}{r-i}$.

To see that $((G^*)^\sharp)^O \in \sgen{\gG^\sharp[O]}$
we write $((G^*)^\sharp)^O$ as the sum of
$\tbinom{q-i}{r-i}^{-1} \tbinom{q}{r} \tbinom{n-i}{r-i}$
atoms at $O$, where we choose the support of the
vectors $v$ cyclically: to be formal, identify
$Q$ with $\mb{Z}/Q\mb{Z}$ and assign the $j$th atom
the vector in $\{0,1\}^{\mb{Z}/Q\mb{Z}}$ with support 
$\{ j \tbinom{q-i}{r-i} + x: x \in [\tbinom{q-i}{r-i}] \}$. \qed

\medskip

\nib{Proof of Theorem \ref{rainbow:fixed}.}
The proof is the same as that of Theorem \ref{rainbow:all},
except for the verification of 
$((G^*)^\sharp)^O \in \sgen{\gG^\sharp[O]}$
for any $O \in \Phi/S_q$. The generators $\gG^\sharp[O]$
are vectors of the same form as before except that
now $v$ must be a row of the inclusion matrix $M^r_i(q)$
(discussed after the statement of the theorem).
To see that $((G^*)^\sharp)^O \in \sgen{\gG^\sharp[O]}$,
we note that the sum $\sS^O$ of all $\gG^\sharp$-atoms at $O$
is supported (as before) on the coordinates $(B,\psi')$ 
with $B' \sub B \in Q$ and $\psi' \in O \cap \Phi_{B'}$,
with each coordinate equal to the vector in $\mb{Z}^Q$
that is $(r-i)! \tbinom{r}{i}$ in each coordinate.
Recalling that $((G^*)^\sharp)^O$ has the same description
with $\tbinom{r}{i}$ replaced by $\tbinom{n-i}{r-i}$,
we have $((G^*)^\sharp)^O = \tbinom{r}{i}^{-1}\tbinom{n-i}{r-i} \sS^O$. \qed

\end{document}